\documentclass[11pt]{article} 

\pdfoutput=1

\usepackage{amsmath,amscd,amssymb}
\usepackage[applemac]{inputenc}
\usepackage[english]{babel}
\usepackage{graphicx}
\begin{document}
\title{From non-K\" ahlerian surfaces to Cremona group of $\bb P^2(\bb C)$}

\author{Georges Dloussky\footnote{The author is partially supported by the ANR project MNGNK, decision N° ANR-10-BLAN-0118}}  
\date{ }
\newtheorem{Lem}{Lemma \thesection.}
\newtheorem{Th}[Lem]{Theorem \thesection.}
\newtheorem{Th and Def}[Lem]{Theorem and Definition \thesection.}
\newtheorem{Cor}[Lem]{Corollary \thesection.}
\newtheorem{Def}[Lem]{Definition \thesection.}
\newtheorem{Ex}[Lem]{Example \thesection.}
\newtheorem{Exs}[Lem]{Examples \thesection.}
\newtheorem{Prob}[Lem]{Problem \thesection.}
\newtheorem{Prop}[Lem]{Proposition \thesection.}
\newtheorem{Prop and Def}[Lem]{Proposition and Definition \thesection.}
\newtheorem{Rem}[Lem]{Remark \thesection.}
\newtheorem{Lem and Def}[Lem]{Lemma and Definition \thesection.}
\newtheorem{Exer}[Lem]{Exercise \thesection.}
\newtheorem{Not}[Lem]{Notations \thesection.}
\setcounter{Lem}{0}
\setcounter{section}{-1}
\def\cal{\mathcal}
\def\bb{\mathbb} 
\def\a{\alpha }
\def\b{\beta }
\def\d{\delta }
\def\D{\Delta }
\def\g{\gamma }
\def\G{\Gamma }
\def\d{\delta }
\def\D{\Delta }
\def\e{\epsilon }
\def\z{\zeta }
\def\k{\kappa }
\def\l{\lambda }
\def\L{\Lambda }
\def\m{\mu }
\def\o{\omega }
\def\p{\pi }
\def\P{\Pi }
\def\s{\sigma }
\def\t{\theta }
\def\T{\Theta }
\def\f{\varphi }
\def\Card{{\rm Card\, }}
\def\ker{{\rm Ker\,}}
\def\im{{\rm Im\,}}
\def\coker{{\rm Coker\,}}
\def\codim{{\rm codim\,}}
\def\part{\partial }
\def\dps{\displaystyle }
\def\ot{\otimes }
\def\tr{{\rm tr\,}}
\def\iff{ if and only if }
\def\wh{\widehat }
\def\wt{\widetilde }
\def\la{\langle }
\def\ra{\rangle }
\def\resto#1#2{{
#1\hskip 0.4ex\vline_{\hskip 0.2ex\raisebox{-0,2ex}
{{${\scriptstyle #2}$}}}}}
\def\map{\longrightarrow}
\def\textmap#1{\mathop{\vbox{\ialign{
                                 ##\crcr
     ${\scriptstyle\hfil\;\;#1\;\;\hfil}$\crcr
     \noalign{\kern 1pt\nointerlineskip}
     \rightarrowfill\crcr}}\;}}
\newcommand{\Cal}{\cal}
\def\textlmap#1{\mathop{\vbox{\ialign{
                                 ##\crcr
     ${\scriptstyle\hfil\;\;#1\;\;\hfil}$\crcr
     \noalign{\kern-1pt\nointerlineskip}
     \leftarrowfill\crcr}}\;}}





\newcommand{\tlowername}[2]%
{$\stackrel{\makebox[1pt]{#1}}%
{\begin{picture}(0,0)%
\put(0,0){\makebox(0,6)[t]{\makebox[1pt]{$#2$}}}%
\end{picture}}$}%

\newcommand{\tcase}[1]{\makebox[23pt]%
{\raisebox{2.5pt}{#1{20}}}}%

\newcommand{\Tcase}[2]{\makebox[23pt]%
{\raisebox{2.5pt}{$\stackrel{#2}{#1{20}}$}}}%

\newcommand{\tbicase}[1]{\makebox[23pt]%
{\raisebox{1pt}{#1{20}}}}%

\newcommand{\Tbicase}[3]{\makebox[23pt]{\raisebox{-7pt}%
{$\stackrel{#2}{\mbox{\tlowername{#1{20}}{\scriptstyle{#3}}}}$}}}%


\newcommand{\AR}[1]%
{\begin{picture}(#1,0)%
\put(0,0){\vector(1,0){#1}}%
\end{picture}}%

\newcommand{\DOTAR}[1]%
{\NUMBEROFDOTS=#1%
\divide\NUMBEROFDOTS by 3%
\begin{picture}(#1,0)%
\multiput(0,0)(3,0){\NUMBEROFDOTS}{\circle*{1}}%
\put(#1,0){\vector(1,0){0}}%
\end{picture}}%

\newcommand{\MONO}[1]%
{\begin{picture}(#1,0)%
\put(0,0){\vector(1,0){#1}}%
\put(2,-2){\line(0,1){4}}%
\end{picture}}%

\newcommand{\EPI}[1]%
{\begin{picture}(#1,0)(-#1,0)%
\put(-#1,0){\vector(1,0){#1}}%
\put(-6,-2){\line(0,1){4}}%
\end{picture}}%

\newcommand{\BIMO}[1]%
{\begin{picture}(#1,0)(-#1,0)%
\put(-#1,0){\vector(1,0){#1}}%
\put(-6,-2){\line(0,1){4}}%
\put(-#1,-2){\hspace{2pt}\line(0,1){4}}%
\end{picture}}%

\newcommand{\BIAR}[1]%
{\begin{picture}(#1,4)%
\put(0,0){\vector(1,0){#1}}%
\put(0,4){\vector(1,0){#1}}%
\end{picture}}%

\newcommand{\EQL}[1]%
{\begin{picture}(#1,0)%
\put(0,1){\line(1,0){#1}}%
\put(0,-1){\line(1,0){#1}}%
\end{picture}}%

\newcommand{\ADJAR}[1]%
{\begin{picture}(#1,4)%
\put(0,0){\vector(1,0){#1}}%
\put(#1,4){\vector(-1,0){#1}}%
\end{picture}}%


\newcommand{\ar}{\tcase{\AR}}%

\newcommand{\Ar}[1]{\Tcase{\AR}{#1}}%

\newcommand{\dotar}{\tcase{\DOTAR}}%

\newcommand{\Dotar}[1]{\Tcase{\DOTAR}{#1}}%

\newcommand{\mono}{\tcase{\MONO}}%

\newcommand{\Mono}[1]{\Tcase{\MONO}{#1}}%

\newcommand{\epi}{\tcase{\EPI}}%

\newcommand{\Epi}[1]{\Tcase{\EPI}{#1}}%

\newcommand{\bimo}{\tcase{\BIMO}}%

\newcommand{\Bimo}[1]{\Tcase{\BIMO}{#1}}%

\newcommand{\iso}{\Tcase{\AR}{\cong}}%

\newcommand{\Iso}[1]{\Tcase{\AR}{\cong{#1}}}%

\newcommand{\biar}{\tbicase{\BIAR}}%

\newcommand{\Biar}[2]{\Tbicase{\BIAR}{#1}{#2}}%

\newcommand{\eql}{\tcase{\EQL}}%

\newcommand{\adjar}{\tbicase{\ADJAR}}%

\newcommand{\Adjar}[2]{\Tbicase{\ADJAR}{#1}{#2}}%


\newcommand{\BKAR}[1]%
{\begin{picture}(#1,0)%
\put(#1,0){\vector(-1,0){#1}}%
\end{picture}}%

\newcommand{\BKDOTAR}[1]%
{\NUMBEROFDOTS=#1%
\divide\NUMBEROFDOTS by 3%
\begin{picture}(#1,0)%
\multiput(#1,0)(-3,0){\NUMBEROFDOTS}{\circle*{1}}%
\put(0,0){\vector(-1,0){0}}%
\end{picture}}%

\newcommand{\BKMONO}[1]%
{\begin{picture}(#1,0)(-#1,0)%
\put(0,0){\vector(-1,0){#1}}%
\put(-2,-2){\line(0,1){4}}%
\end{picture}}%

\newcommand{\BKEPI}[1]%
{\begin{picture}(#1,0)%
\put(#1,0){\vector(-1,0){#1}}%
\put(6,-2){\line(0,1){4}}%
\end{picture}}%

\newcommand{\BKBIMO}[1]%
{\begin{picture}(#1,0)%
\put(#1,0){\vector(-1,0){#1}}%
\put(6,-2){\line(0,1){4}}%
\put(#1,-2){\hspace{-2pt}\line(0,1){4}}%
\end{picture}}%

\newcommand{\BKBIAR}[1]%
{\begin{picture}(#1,4)%
\put(#1,0){\vector(-1,0){#1}}%
\put(#1,4){\vector(-1,0){#1}}%
\end{picture}}%

\newcommand{\BKADJAR}[1]%
{\begin{picture}(#1,4)%
\put(0,4){\vector(1,0){#1}}%
\put(#1,0){\vector(-1,0){#1}}%
\end{picture}}%


\newcommand{\bkar}{\tcase{\BKAR}}%

\newcommand{\Bkar}[1]{\Tcase{\BKAR}{#1}}%

\newcommand{\bkdotar}{\tcase{\BKDOTAR}}%

\newcommand{\Bkdotar}[1]{\Tcase{\BKDOTAR}{#1}}%

\newcommand{\bkmono}{\tcase{\BKMONO}}%

\newcommand{\Bkmono}[1]{\Tcase{\BKMONO}{#1}}%

\newcommand{\bkepi}{\tcase{\BKEPI}}%

\newcommand{\Bkepi}[1]{\Tcase{\BKEPI}{#1}}%

\newcommand{\bkbimo}{\tcase{\BKBIMO}}%

\newcommand{\Bkbimo}[1]{\Tcase{\BKBIMO}{\hspace{9pt}#1}}%

\newcommand{\bkiso}{\Tcase{\BKAR}{\cong}}%

\newcommand{\Bkiso}[1]{\Tcase{\BKAR}{\cong{#1}}}%

\newcommand{\bkbiar}{\tbicase{\BKBIAR}}%

\newcommand{\Bkbiar}[2]{\Tbicase{\BKBIAR}{#1}{#2}}%

\newcommand{\bkeql}{\tcase{\EQL}}%

\newcommand{\bkadjar}{\tbicase{\BKADJAR}}%

\newcommand{\Bkadjar}[2]{\Tbicase{\BKADJAR}{#1}{#2}}%


\newcommand{\lowername}[2]%
{$\stackrel{\makebox[1pt]{#1}}%
{\begin{picture}(0,0)%
\truex{600}%
\put(0,0){\makebox(0,\value{x})[t]{\makebox[1pt]{$#2$}}}%
\end{picture}}$}%

\newcommand{\hcase}[2]%
{\makebox[0pt]%
{\raisebox{-1pt}[0pt][0pt]{#1{#2}}}}%

\newcommand{\Hcase}[3]%
{\makebox[0pt]
{\raisebox{-1pt}[0pt][0pt]%
{$\stackrel{\makebox[0pt]{$\textstyle{#2}$}}{#1{#3}}$}}}%

\newcommand{\hcasE}[3]%
{\makebox[0pt]%
{\raisebox{-9pt}[0pt][0pt]%
{\lowername{#1{#3}}{#2}}}}%

\newcommand{\hbicase}[2]%
{\makebox[0pt]%
{\raisebox{-2.5pt}[0pt][0pt]{#1{#2}}}}%

\newcommand{\Hbicase}[4]%
{\makebox[0pt]
{\raisebox{-10.5pt}[0pt][0pt]%
{$\stackrel{\makebox[0pt]{$\textstyle{#2}$}}%
{\mbox{\lowername{#1{#4}}{#3}}}$}}}%


\newcommand{\EAR}[1]%
{\begin{picture}(#1,0)%
\put(0,0){\vector(1,0){#1}}%
\end{picture}}%

\newcommand{\EDOTAR}[1]%
{\truex{100}\truey{300}%
\NUMBEROFDOTS=#1%
\divide\NUMBEROFDOTS by \value{y}%
\begin{picture}(#1,0)%
\multiput(0,0)(\value{y},0){\NUMBEROFDOTS}%
{\circle*{\value{x}}}%
\put(#1,0){\vector(1,0){0}}%
\end{picture}}%

\newcommand{\EMONO}[1]%
{\begin{picture}(#1,0)%
\put(0,0){\vector(1,0){#1}}%
\truex{300}\truey{600}%
\put(\value{x},-\value{x}){\line(0,1){\value{y}}}%
\end{picture}}%

\newcommand{\EEPI}[1]%
{\begin{picture}(#1,0)(-#1,0)%
\put(-#1,0){\vector(1,0){#1}}%
\truex{300}\truey{600}\truez{800}%
\put(-\value{z},-\value{x}){\line(0,1){\value{y}}}%
\end{picture}}%

\newcommand{\EBIMO}[1]%
{\begin{picture}(#1,0)(-#1,0)%
\put(-#1,0){\vector(1,0){#1}}%
\truex{300}\truey{600}\truez{800}%
\put(-\value{z},-\value{x}){\line(0,1){\value{y}}}%
\put(-#1,-\value{x}){\hspace{3pt}\line(0,1){\value{y}}}%
\end{picture}}%

\newcommand{\EBIAR}[1]%
{\truex{400}%
\begin{picture}(#1,\value{x})%
\put(0,0){\vector(1,0){#1}}%
\put(0,\value{x}){\vector(1,0){#1}}%
\end{picture}}%

\newcommand{\EEQL}[1]%
{\begin{picture}(#1,0)%
\truex{200}%
\put(0,\value{x}){\line(1,0){#1}}%
\put(0,0){\line(1,0){#1}}%
\end{picture}}%

\newcommand{\EADJAR}[1]%
{\truex{400}%
\begin{picture}(#1,\value{x})%
\put(0,0){\vector(1,0){#1}}%
\put(#1,\value{x}){\vector(-1,0){#1}}%
\end{picture}}%


\newcommand{\earv}[1]{\hcase{\EAR}{#100}}%

\newcommand{\ear}%
{\hspace{\SOURCE\unitlength}%
\hcase{\EAR}{\ARROWLENGTH}}%

\newcommand{\Earv}[2]{\Hcase{\EAR}{#1}{#200}}%

\newcommand{\Ear}[1]%
{\hspace{\SOURCE\unitlength}%
\Hcase{\EAR}{#1}{\ARROWLENGTH}}%

\newcommand{\eaRv}[2]{\hcasE{\EAR}{#1}{#200}}%

\newcommand{\eaR}[1]%
{\hspace{\SOURCE\unitlength}%
\hcasE{\EAR}{#1}{\ARROWLENGTH}}%

\newcommand{\edotarv}[1]{\hcase{\EDOTAR}{#100}}%

\newcommand{\edotar}%
{\hspace{\SOURCE\unitlength}%
\hcase{\EDOTAR}{\ARROWLENGTH}}%

\newcommand{\Edotarv}[2]{\Hcase{\EDOTAR}{#1}{#200}}%

\newcommand{\Edotar}[1]%
{\hspace{\SOURCE\unitlength}%
\Hcase{\EDOTAR}{#1}{\ARROWLENGTH}}%

\newcommand{\edotaRv}[2]{\hcasE{\EDOTAR}{#1}{#200}}%

\newcommand{\edotaR}[1]%
{\hspace{\SOURCE\unitlength}%
\hcasE{\EDOTAR}{#1}{\ARROWLENGTH}}%

\newcommand{\emonov}[1]{\hcase{\EMONO}{#100}}%

\newcommand{\emono}%
{\hspace{\SOURCE\unitlength}%
\hcase{\EMONO}{\ARROWLENGTH}}%

\newcommand{\Emonov}[2]{\Hcase{\EMONO}{#1}{#200}}%

\newcommand{\Emono}[1]%
{\hspace{\SOURCE\unitlength}%
\Hcase{\EMONO}{#1}{\ARROWLENGTH}}%

\newcommand{\emonOv}[2]{\hcasE{\EMONO}{#1}{#200}}%

\newcommand{\emonO}[1]%
{\hspace{\SOURCE\unitlength}%
\hcasE{\EMONO}{#1}{\ARROWLENGTH}}%

\newcommand{\eepiv}[1]{\hcase{\EEPI}{#100}}%

\newcommand{\eepi}%
{\hspace{\SOURCE\unitlength}%
\hcase{\EEPI}{\ARROWLENGTH}}%

\newcommand{\Eepiv}[2]{\Hcase{\EEPI}{#1}{#200}}%

\newcommand{\Eepi}[1]%
{\hspace{\SOURCE\unitlength}%
\Hcase{\EEPI}{#1}{\ARROWLENGTH}}%

\newcommand{\eepIv}[2]{\hcasE{\EEPI}{#1}{#200}}%

\newcommand{\eepI}[1]%
{\hspace{\SOURCE\unitlength}%
\hcasE{\EEPI}{#1}{\ARROWLENGTH}}%

\newcommand{\ebimov}[1]{\hcase{\EBIMO}{#100}}%

\newcommand{\ebimo}%
{\hspace{\SOURCE\unitlength}%
\hcase{\EBIMO}{\ARROWLENGTH}}%

\newcommand{\Ebimov}[2]{\Hcase{\EBIMO}{#1}{#200}}%

\newcommand{\Ebimo}[1]%
{\hspace{\SOURCE\unitlength}%
\Hcase{\EBIMO}{#1}{\ARROWLENGTH}}%

\newcommand{\ebimOv}[2]{\hcasE{\EBIMO}{#1}{#200}}%

\newcommand{\ebimO}[1]%
{\hspace{\SOURCE\unitlength}%
\hcasE{\EBIMO}{#1}{\ARROWLENGTH}}%

\newcommand{\eisov}[1]{\Hcase{\EAR}{\cong}{#100}}%

\newcommand{\eiso}%
{\hspace{\SOURCE\unitlength}%
\Hcase{\EAR}{\cong}{\ARROWLENGTH}}%

\newcommand{\Eisov}[2]{\Hcase{\EAR}{\cong#1}{#200}}%

\newcommand{\Eiso}[1]%
{\hspace{\SOURCE\unitlength}%
\Hcase{\EAR}{\cong#1}{\ARROWLENGTH}}%

\newcommand{\eisOv}[2]{\hcasE{\EAR}{\cong#1}{#200}}%

\newcommand{\eisO}[1]%
{\hspace{\SOURCE\unitlength}%
\hcasE{\EAR}{\cong#1}{\ARROWLENGTH}}%

\newcommand{\ebiarv}[1]{\hbicase{\EBIAR}{#100}}%

\newcommand{\ebiar}%
{\hspace{\SOURCE\unitlength}%
\hbicase{\EBIAR}{\ARROWLENGTH}}%

\newcommand{\Ebiarv}[3]{\Hbicase{\EBIAR}{#1}{#2}{#300}}%

\newcommand{\Ebiar}[2]%
{\hspace{\SOURCE\unitlength}%
\Hbicase{\EBIAR}{#1}{#2}{\ARROWLENGTH}}%

\newcommand{\eeqlv}[1]{\hcase{\EEQL}{#100}}%

\newcommand{\eeql}%
{\hspace{\SOURCE\unitlength}%
\hbicase{\EEQL}{\ARROWLENGTH}}%

\newcommand{\eadjarv}[1]{\hbicase{\EADJAR}{#100}}%

\newcommand{\eadjar}%
{\hspace{\SOURCE\unitlength}%
\hbicase{\EADJAR}{\ARROWLENGTH}}%

\newcommand{\Eadjarv}[3]{\Hbicase{\EADJAR}{#1}{#2}{#300}}%

\newcommand{\Eadjar}[2]%
{\hspace{\SOURCE\unitlength}%
\Hbicase{\EADJAR}{#1}{#2}{\ARROWLENGTH}}%


\newcommand{\WAR}[1]%
{\begin{picture}(#1,0)%
\put(#1,0){\vector(-1,0){#1}}%
\end{picture}}%

\newcommand{\WDOTAR}[1]%
{\truex{100}\truey{300}%
\NUMBEROFDOTS=#1%
\divide\NUMBEROFDOTS by \value{y}%
\begin{picture}(#1,0)%
\multiput(#1,0)(-\value{y},0){\NUMBEROFDOTS}%
{\circle*{\value{x}}}%
\put(0,0){\vector(-1,0){0}}%
\end{picture}}%

\newcommand{\WMONO}[1]%
{\begin{picture}(#1,0)(-#1,0)%
\put(0,0){\vector(-1,0){#1}}%
\truex{300}\truey{600}%
\put(-\value{x},-\value{x}){\line(0,1){\value{y}}}%
\end{picture}}%

\newcommand{\WEPI}[1]%
{\begin{picture}(#1,0)%
\put(#1,0){\vector(-1,0){#1}}%
\truex{300}\truey{600}\truez{800}%
\put(\value{z},-\value{x}){\line(0,1){\value{y}}}%
\end{picture}}%

\newcommand{\WBIMO}[1]%
{\begin{picture}(#1,0)%
\put(#1,0){\vector(-1,0){#1}}%
\truex{300}\truey{600}\truez{800}%
\put(\value{z},-\value{x}){\line(0,1){\value{y}}}%
\put(#1,-\value{x}){\hspace{-3pt}\line(0,1){\value{y}}}%
\end{picture}}%

\newcommand{\WBIAR}[1]%
{\truex{400}%
\begin{picture}(#1,\value{x})%
\put(#1,0){\vector(-1,0){#1}}%
\put(#1,\value{x}){\vector(-1,0){#1}}%
\end{picture}}%

\newcommand{\WADJAR}[1]%
{\truex{400}%
\begin{picture}(#1,\value{x})%
\put(0,\value{x}){\vector(1,0){#1}}%
\put(#1,0){\vector(-1,0){#1}}%
\end{picture}}%


\newcommand{\warv}[1]{\hcase{\WAR}{#100}}%

\newcommand{\war}%
{\hspace{\SOURCE\unitlength}%
\hcase{\WAR}{\ARROWLENGTH}}%

\newcommand{\Warv}[2]{\Hcase{\WAR}{#1}{#200}}%

\newcommand{\War}[1]%
{\hspace{\SOURCE\unitlength}%
\Hcase{\WAR}{#1}{\ARROWLENGTH}}%

\newcommand{\waRv}[2]{\hcasE{\WAR}{#1}{#200}}%

\newcommand{\waR}[1]%
{\hspace{\SOURCE\unitlength}%
\hcasE{\WAR}{#1}{\ARROWLENGTH}}%

\newcommand{\wdotarv}[1]{\hcase{\WDOTAR}{#100}}%

\newcommand{\wdotar}%
{\hspace{\SOURCE\unitlength}%
\hcase{\WDOTAR}{\ARROWLENGTH}}%

\newcommand{\Wdotarv}[2]{\Hcase{\WDOTAR}{#1}{#200}}%

\newcommand{\Wdotar}[1]%
{\hspace{\SOURCE\unitlength}%
\Hcase{\WDOTAR}{#1}{\ARROWLENGTH}}%

\newcommand{\wdotaRv}[2]{\hcasE{\WDOTAR}{#1}{#200}}%

\newcommand{\wdotaR}[1]%
{\hspace{\SOURCE\unitlength}%
\hcasE{\WDOTAR}{#1}{\ARROWLENGTH}}%

\newcommand{\wmonov}[1]{\hcase{\WMONO}{#100}}%

\newcommand{\wmono}%
{\hspace{\SOURCE\unitlength}%
\hcase{\WMONO}{\ARROWLENGTH}}%

\newcommand{\Wmonov}[2]{\Hcase{\WMONO}{#1}{#200}}%

\newcommand{\Wmono}[1]%
{\hspace{\SOURCE\unitlength}%
\Hcase{\WMONO}{#1}{\ARROWLENGTH}}%

\newcommand{\wmonOv}[2]{\hcasE{\WMONO}{#1}{#200}}%

\newcommand{\wmonO}[1]%
{\hspace{\SOURCE\unitlength}%
\hcasE{\WMONO}{#1}{\ARROWLENGTH}}%

\newcommand{\wepiv}[1]{\hcase{\WEPI}{#100}}%

\newcommand{\wepi}%
{\hspace{\SOURCE\unitlength}%
\hcase{\WEPI}{\ARROWLENGTH}}%

\newcommand{\Wepiv}[2]{\Hcase{\WEPI}{#1}{#200}}%

\newcommand{\Wepi}[1]%
{\hspace{\SOURCE\unitlength}%
\Hcase{\WEPI}{#1}{\ARROWLENGTH}}%

\newcommand{\wepIv}[2]{\hcasE{\WEPI}{#1}{#200}}%

\newcommand{\wepI}[1]%
{\hspace{\SOURCE\unitlength}%
\hcasE{\WEPI}{#1}{\ARROWLENGTH}}%

\newcommand{\wbimov}[1]{\hcase{\WBIMO}{#100}}%

\newcommand{\wbimo}%
{\hspace{\SOURCE\unitlength}%
\hcase{\WBIMO}{\ARROWLENGTH}}%

\newcommand{\Wbimov}[2]{\Hcase{\WBIMO}{#1}{#200}}%

\newcommand{\Wbimo}[1]%
{\hspace{\SOURCE\unitlength}%
\Hcase{\WBIMO}{#1}{\ARROWLENGTH}}%

\newcommand{\wbimOv}[2]{\hcasE{\WBIMO}{#1}{#200}}%

\newcommand{\wbimO}[1]%
{\hspace{\SOURCE\unitlength}%
\hcasE{\WBIMO}{#1}{\ARROWLENGTH}}%

\newcommand{\wisov}[1]{\Hcase{\WAR}{\cong}{#100}}%

\newcommand{\wiso}%
{\hspace{\SOURCE\unitlength}%
\Hcase{\WAR}{\cong}{\ARROWLENGTH}}%

\newcommand{\Wisov}[2]{\Hcase{\WAR}{\cong#1}{#200}}%

\newcommand{\Wiso}[1]%
{\hspace{\SOURCE\unitlength}%
\Hcase{\WAR}{#1}{\ARROWLENGTH}}%

\newcommand{\wisOv}[2]{\hcasE{\WAR}{\cong#1}{#200}}%

\newcommand{\wisO}[1]%
{\hspace{\SOURCE\unitlength}%
\hcasE{\WAR}{#1}{\ARROWLENGTH}}%

\newcommand{\wbiarv}[1]{\hbicase{\WBIAR}{#100}}%

\newcommand{\wbiar}%
{\hspace{\SOURCE\unitlength}%
\hbicase{\WBIAR}{\ARROWLENGTH}}%

\newcommand{\Wbiarv}[3]{\Hbicase{\WBIAR}{#1}{#2}{#300}}%

\newcommand{\Wbiar}[2]%
{\hspace{\SOURCE\unitlength}%
\Hbicase{\WBIAR}{#1}{#2}{\ARROWLENGTH}}%

\newcommand{\weqlv}[1]{\hbicase{\EEQL}{#100}}%

\newcommand{\weql}%
{\hspace{\SOURCE\unitlength}%
\hbicase{\EEQL}{\ARROWLENGTH}}%

\newcommand{\wadjarv}[1]{\hbicase{\WADJAR}{#100}}%

\newcommand{\wadjar}%
{\hspace{\SOURCE\unitlength}%
\hbicase{\WADJAR}{\ARROWLENGTH}}%

\newcommand{\Wadjarv}[3]{\Hbicase{\WADJAR}{#1}{#2}{#300}}%

\newcommand{\Wadjar}[2]%
{\hspace{\SOURCE\unitlength}%
\Hbicase{\WADJAR}{#1}{#2}{\ARROWLENGTH}}%


\newcommand{\vcase}[2]{#1{#2}}%

\newcommand{\Vcase}[3]{\makebox[0pt]%
{\makebox[0pt][r]{\raisebox{0pt}[0pt][0pt]{${#2}\hspace{2pt}$}}}#1{#3}}%

\newcommand{\vcasE}[3]{\makebox[0pt]%
{#1{#3}\makebox[0pt][l]{\raisebox{0pt}[0pt][0pt]{\hspace{2pt}$#2$}}}}%

\newcommand{\vbicase}[2]{\makebox[0pt]{{#1{#2}}}}%

\newcommand{\Vbicase}[4]{\makebox[0pt]%
{\makebox[0pt][r]{\raisebox{0pt}[0pt][0pt]{$#2$\hspace{4pt}}}#1{#4}%
\makebox[0pt][l]{\raisebox{0pt}[0pt][0pt]{\hspace{5pt}$#3$}}}}%


\newcommand{\SAR}[1]%
{\begin{picture}(0,0)%
\put(0,0){\makebox(0,0)%
{\begin{picture}(0,#1)%
\put(0,#1){\vector(0,-1){#1}}%
\end{picture}}}\end{picture}}%

\newcommand{\SDOTAR}[1]%
{\truex{100}\truey{300}%
\NUMBEROFDOTS=#1%
\divide\NUMBEROFDOTS by \value{y}%
\begin{picture}(0,0)%
\put(0,0){\makebox(0,0)%
{\begin{picture}(0,#1)%
\multiput(0,#1)(0,-\value{y}){\NUMBEROFDOTS}%
{\circle*{\value{x}}}%
\put(0,0){\vector(0,-1){0}}%
\end{picture}}}\end{picture}}%

\newcommand{\SMONO}[1]%
{\begin{picture}(0,0)%
\put(0,0){\makebox(0,0)%
{\begin{picture}(0,#1)%
\put(0,#1){\vector(0,-1){#1}}%
\truex{300}\truey{600}%
\put(0,#1){\begin{picture}(0,0)%
\put(-\value{x},-\value{x}){\line(1,0){\value{y}}}\end{picture}}%
\end{picture}}}\end{picture}}%

\newcommand{\SEPI}[1]%
{\begin{picture}(0,0)%
\put(0,0){\makebox(0,0)%
{\begin{picture}(0,#1)%
\put(0,#1){\vector(0,-1){#1}}%
\truex{300}\truey{600}\truez{800}%
\put(-\value{x},\value{z}){\line(1,0){\value{y}}}%
\end{picture}}}\end{picture}}%

\newcommand{\SBIMO}[1]%
{\begin{picture}(0,0)%
\put(0,0){\makebox(0,0)%
{\begin{picture}(0,#1)%
\put(0,#1){\vector(0,-1){#1}}%
\truex{300}\truey{600}\truez{800}%
\put(0,#1){\begin{picture}(0,0)%
\put(-\value{x},-\value{x}){\line(1,0){\value{y}}}\end{picture}}%
\put(-\value{x},\value{z}){\line(1,0){\value{y}}}%
\end{picture}}}\end{picture}}%

\newcommand{\SBIAR}[1]%
{\begin{picture}(0,0)%
\truex{200}%
\put(0,0){\makebox(0,0)%
{\begin{picture}(0,#1)\put(-\value{x},#1){\vector(0,-1){#1}}%
\put(\value{x},#1){\vector(0,-1){#1}}%
\end{picture}}}\end{picture}}%

\newcommand{\SEQL}[1]%
{\begin{picture}(0,0)%
\truex{100}%
\put(0,0){\makebox(0,0)%
{\begin{picture}(0,#1)\put(-\value{x},#1){\line(0,-1){#1}}%
\put(\value{x},#1){\line(0,-1){#1}}%
\end{picture}}}\end{picture}}%

\newcommand{\SADJAR}[1]{\begin{picture}(0,0)%
\truex{200}%
\put(0,0){\makebox(0,0)%
{\begin{picture}(0,#1)\put(-\value{x},#1){\vector(0,-1){#1}}%
\put(\value{x},0){\vector(0,1){#1}}%
\end{picture}}}\end{picture}}%


\newcommand{\sarv}[1]{\vcase{\SAR}{#100}}%

\newcommand{\sar}{\sarv{50}}%

\newcommand{\Sarv}[2]{\Vcase{\SAR}{#1}{#200}}%

\newcommand{\Sar}[1]{\Sarv{#1}{50}}%

\newcommand{\saRv}[2]{\vcasE{\SAR}{#1}{#200}}%

\newcommand{\saR}[1]{\saRv{#1}{50}}%

\newcommand{\sdotarv}[1]{\vcase{\SDOTAR}{#100}}%

\newcommand{\sdotar}{\sdotarv{50}}%

\newcommand{\Sdotarv}[2]{\Vcase{\SDOTAR}{#1}{#200}}%

\newcommand{\Sdotar}[1]{\Sdotarv{#1}{50}}%

\newcommand{\sdotaRv}[2]{\vcasE{\SDOTAR}{#1}{#200}}%

\newcommand{\sdotaR}[1]{\sdotaRv{#1}{50}}%

\newcommand{\smonov}[1]{\vcase{\SMONO}{#100}}%

\newcommand{\smono}{\smonov{50}}%

\newcommand{\Smonov}[2]{\Vcase{\SMONO}{#1}{#200}}%

\newcommand{\Smono}[1]{\Smonov{#1}{50}}%

\newcommand{\smonOv}[2]{\vcasE{\SMONO}{#1}{#200}}%

\newcommand{\smonO}[1]{\smonOv{#1}{50}}%

\newcommand{\sepiv}[1]{\vcase{\SEPI}{#100}}%

\newcommand{\sepi}{\sepiv{50}}%

\newcommand{\Sepiv}[2]{\Vcase{\SEPI}{#1}{#200}}%

\newcommand{\Sepi}[1]{\Sepiv{#1}{50}}%

\newcommand{\sepIv}[2]{\vcasE{\SEPI}{#1}{#200}}%

\newcommand{\sepI}[1]{\sepIv{#1}{50}}%

\newcommand{\sbimov}[1]{\vcase{\SBIMO}{#100}}%

\newcommand{\sbimo}{\sbimov{50}}%

\newcommand{\Sbimov}[2]{\Vcase{\SBIMO}{#1}{#200}}%

\newcommand{\Sbimo}[1]{\Sbimov{#1}{50}}%

\newcommand{\sbimOv}[2]{\vcasE{\SBIMO}{#1}{#200}}%

\newcommand{\sbimO}[1]{\sbimOv{#1}{50}}%

\newcommand{\sisov}[1]{\vcasE{\SAR}{\cong}{#100}}%

\newcommand{\siso}{\sisov{50}}%

\newcommand{\Sisov}[2]%
{\Vbicase{\SAR}{#1\hspace{-2pt}}{\hspace{-2pt}\cong}{#200}}%

\newcommand{\Siso}[1]{\Sisov{#1}{50}}%

\newcommand{\sbiarv}[1]{\vbicase{\SBIAR}{#100}}%

\newcommand{\sbiar}{\sbiarv{50}}%

\newcommand{\Sbiarv}[3]{\Vbicase{\SBIAR}{#1}{#2}{#300}}%

\newcommand{\Sbiar}[2]{\Sbiarv{#1}{#2}{50}}%

\newcommand{\seqlv}[1]{\vbicase{\SEQL}{#100}}%

\newcommand{\seql}{\seqlv{50}}%

\newcommand{\sadjarv}[1]{\vbicase{\SADJAR}{#100}}%

\newcommand{\sadjar}{\sadjarv{50}}%

\newcommand{\Sadjarv}[3]{\Vbicase{\SADJAR}{#1}{#2}{#300}}%

\newcommand{\Sadjar}[2]{\Sadjarv{#1}{#2}{50}}%


\newcommand{\NAR}[1]%
{\begin{picture}(0,0)%
\put(0,0){\makebox(0,0)%
{\begin{picture}(0,#1)\put(0,0){\vector(0,1){#1}}%
\end{picture}}}\end{picture}}%

\newcommand{\NDOTAR}[1]%
{\truex{100}\truey{300}%
\NUMBEROFDOTS=#1%
\divide\NUMBEROFDOTS by \value{y}%
\begin{picture}(0,0)%
\put(0,0){\makebox(0,0)%
{\begin{picture}(0,#1)%
\multiput(0,0)(0,\value{y}){\NUMBEROFDOTS}%
{\circle*{\value{x}}}%
\put(0,#1){\vector(0,1){0}}%
\end{picture}}}\end{picture}}%

\newcommand{\NMONO}[1]%
{\begin{picture}(0,0)%
\put(0,0){\makebox(0,0)%
{\begin{picture}(0,#1)%
\put(0,0){\vector(0,1){#1}}%
\truex{300}\truey{600}%
\put(-\value{x},\value{x}){\line(1,0){\value{y}}}%
\end{picture}}}%
\end{picture}}%

\newcommand{\NEPI}[1]%
{\begin{picture}(0,0)%
\put(0,0){\makebox(0,0)%
{\begin{picture}(0,#1)%
\put(0,0){\vector(0,1){#1}}%
\truex{300}\truey{600}\truez{800}%
\put(0,#1){\begin{picture}(0,0)%
\put(-\value{x},-\value{z}){\line(1,0){\value{y}}}\end{picture}}%
\end{picture}}}\end{picture}}%

\newcommand{\NBIMO}[1]%
{\begin{picture}(0,0)%
\put(0,0){\makebox(0,0)%
{\begin{picture}(0,#1)%
\put(0,0){\vector(0,1){#1}}%
\truex{300}\truey{600}\truez{800}%
\put(-\value{x},\value{x}){\line(1,0){\value{y}}}%
\put(0,#1){\begin{picture}(0,0)%
\put(-\value{x},-\value{z}){\line(1,0){\value{y}}}\end{picture}}%
\end{picture}}}\end{picture}}%

\newcommand{\NBIAR}[1]%
{\begin{picture}(0,0)%
\truex{200}%
\put(0,0){\makebox(0,0)%
{\begin{picture}(0,#1)\put(-\value{x},0){\vector(0,1){#1}}%
\put(\value{x},0){\vector(0,1){#1}}%
\end{picture}}}\end{picture}}%

\newcommand{\NADJAR}[1]{\begin{picture}(0,0)%
\truex{200}%
\put(0,0){\makebox(0,0)%
{\begin{picture}(0,#1)\put(\value{x},#1){\vector(0,-1){#1}}%
\put(-\value{x},0){\vector(0,1){#1}}%
\end{picture}}}\end{picture}}%


\newcommand{\narv}[1]{\vcase{\NAR}{#100}}%

\newcommand{\nar}{\narv{50}}%

\newcommand{\Narv}[2]{\Vcase{\NAR}{#1}{#200}}%

\newcommand{\Nar}[1]{\Narv{#1}{50}}%

\newcommand{\naRv}[2]{\vcasE{\NAR}{#1}{#200}}%

\newcommand{\naR}[1]{\naRv{#1}{50}}%

\newcommand{\ndotarv}[1]{\vcase{\NDOTAR}{#100}}%

\newcommand{\ndotar}{\ndotarv{50}}%

\newcommand{\Ndotarv}[2]{\Vcase{\NDOTAR}{#1}{#200}}%

\newcommand{\Ndotar}[1]{\Ndotarv{#1}{50}}%

\newcommand{\ndotaRv}[2]{\vcasE{\NDOTAR}{#1}{#200}}%

\newcommand{\ndotaR}[1]{\ndotaRv{#1}{50}}%

\newcommand{\nmonov}[1]{\vcase{\NMONO}{#100}}%

\newcommand{\nmono}{\nmonov{50}}%

\newcommand{\Nmonov}[2]{\Vcase{\NMONO}{#1}{#200}}%

\newcommand{\Nmono}[1]{\Nmonov{#1}{50}}%

\newcommand{\nmonOv}[2]{\vcasE{\NMONO}{#1}{#200}}%

\newcommand{\nmonO}[1]{\nmonOv{#1}{50}}%

\newcommand{\nepiv}[1]{\vcase{\NEPI}{#100}}%

\newcommand{\nepi}{\nepiv{50}}%

\newcommand{\Nepiv}[2]{\Vcase{\NEPI}{#1}{#200}}%

\newcommand{\Nepi}[1]{\Nepiv{#1}{50}}%

\newcommand{\nepIv}[2]{\vcasE{\NEPI}{#1}{#200}}%

\newcommand{\nepI}[1]{\nepIv{#1}{50}}%

\newcommand{\nbimov}[1]{\vcase{\NBIMO}{#100}}%

\newcommand{\nbimo}{\nbimov{50}}%

\newcommand{\Nbimov}[2]{\Vcase{\NBIMO}{#1}{#200}}%

\newcommand{\Nbimo}[1]{\Nbimov{#1}{50}}%

\newcommand{\nbimOv}[2]{\vcasE{\NBIMO}{#1}{#200}}%

\newcommand{\nbimO}[1]{\nbimOv{#1}{50}}%

\newcommand{\nisov}[1]{\vcasE{\NAR}{\cong}{#100}}%

\newcommand{\niso}{\nisov{50}}%

\newcommand{\Nisov}[2]%
{\Vbicase{\NAR}{#1\hspace{-2pt}}{\hspace{-2pt}\cong}{#200}}%

\newcommand{\Niso}[1]{\Nisov{#1}{50}}%

\newcommand{\nbiarv}[1]{\vbicase{\NBIAR}{#100}}%

\newcommand{\nbiar}{\nbiarv{50}}%

\newcommand{\Nbiarv}[3]{\Vbicase{\NBIAR}{#1}{#2}{#300}}%

\newcommand{\Nbiar}[2]{\Nbiarv{#1}{#2}{50}}%

\newcommand{\neqlv}[1]{\vbicase{\SEQL}{#100}}%

\newcommand{\neql}{\neqlv{50}}%

\newcommand{\nadjarv}[1]{\vbicase{\NADJAR}{#100}}%

\newcommand{\nadjar}{\nadjarv{50}}%

\newcommand{\Nadjarv}[3]{\Vbicase{\NADJAR}{#1}{#2}{#300}}%

\newcommand{\Nadjar}[2]{\Nadjarv{#1}{#2}{50}}%


\newcommand{\fdcase}[3]{\begin{picture}(0,0)%
\put(0,-150){#1}%
\truex{200}\truey{600}\truez{600}%
\put(-\value{x},-\value{x}){\makebox(0,\value{z})[r]{${#2}$}}%
\put(\value{x},-\value{y}){\makebox(0,\value{z})[l]{${#3}$}}%
\end{picture}}%

\newcommand{\fdbicase}[3]{\begin{picture}(0,0)%
\put(0,-150){#1}%
\truex{800}\truey{50}%
\put(-\value{x},\value{y}){${#2}$}%
\truex{200}\truey{950}%
\put(\value{x},-\value{y}){${#3}$}%
\end{picture}}%


\newcommand{\NEAR}{\begin{picture}(0,0)%
\put(-2900,-2900){\vector(1,1){5800}}%
\end{picture}}%

\newcommand{\NEDOTAR}%
{\truex{100}\truey{212}%
\NUMBEROFDOTS=5800%
\divide\NUMBEROFDOTS by \value{y}%
\begin{picture}(0,0)%
\multiput(-2900,-2900)(\value{y},\value{y}){\NUMBEROFDOTS}%
{\circle*{\value{x}}}%
\put(2900,2900){\vector(1,1){0}}%
\end{picture}}%

\newcommand{\NEMONO}{\begin{picture}(0,0)%
\put(-2900,-2900){\vector(1,1){5800}}%
\put(-2900,-2900){\begin{picture}(0,0)%
\truex{141}%
\put(\value{x},\value{x}){\makebox(0,0){$\times$}}%
\end{picture}}\end{picture}}%

\newcommand{\NEEPI}{\begin{picture}(0,0)%
\put(-2900,-2900){\vector(1,1){5800}}%
\put(2900,2900){\begin{picture}(0,0)%
\truex{545}%
\put(-\value{x},-\value{x}){\makebox(0,0){$\times$}}%
\end{picture}}\end{picture}}%

\newcommand{\NEBIMO}{\begin{picture}(0,0)%
\put(-2900,-2900){\vector(1,1){5800}}%
\put(2900,2900){\begin{picture}(0,0)%
\truex{545}%
\put(-\value{x},-\value{x}){\makebox(0,0){$\times$}}%
\end{picture}}
\put(-2900,-2900){\begin{picture}(0,0)%
\truex{141}%
\put(\value{x},\value{x}){\makebox(0,0){$\times$}}%
\end{picture}}\end{picture}}%

\newcommand{\NEBIAR}{\begin{picture}(0,0)%
\put(-2900,-2900){\begin{picture}(0,0)%
\truex{141}%
\put(-\value{x},\value{x}){\vector(1,1){5800}}%
\put(\value{x},-\value{x}){\vector(1,1){5800}}%
\end{picture}}\end{picture}}%

\newcommand{\NEEQL}{\begin{picture}(0,0)%
\put(-2900,-2900){\begin{picture}(0,0)%
\truex{70}%
\put(-\value{x},\value{x}){\line(1,1){5800}}%
\put(\value{x},-\value{x}){\line(1,1){5800}}%
\end{picture}}\end{picture}}%

\newcommand{\NEADJAR}{\begin{picture}(0,0)%
\put(-2900,-2900){\begin{picture}(0,0)%
\truex{141}%
\put(\value{x},-\value{x}){\vector(1,1){5800}}%
\end{picture}}%
\put(2900,2900){\begin{picture}(0,0)%
\truex{141}%
\put(-\value{x},\value{x}){\vector(-1,-1){5800}}%
\end{picture}}\end{picture}}%

\newcommand{\NEARV}[1]{\begin{picture}(0,0)%
\put(0,0){\makebox(0,0){\begin{picture}(#1,#1)%
\put(0,0){\vector(1,1){#1}}\end{picture}}}%
\end{picture}}%


\newcommand{\near}{\fdcase{\NEAR}{}{}}%

\newcommand{\Near}[1]{\fdcase{\NEAR}{#1}{}}%

\newcommand{\neaR}[1]{\fdcase{\NEAR}{}{#1}}%

\newcommand{\nedotar}{\fdcase{\NEDOTAR}{}{}}%

\newcommand{\Nedotar}[1]{\fdcase{\NEDOTAR}{#1}{}}%

\newcommand{\nedotaR}[1]{\fdcase{\NEDOTAR}{}{#1}}%

\newcommand{\nemono}{\fdcase{\NEMONO}{}{}}%

\newcommand{\Nemono}[1]{\fdcase{\NEMONO}{#1}{}}%

\newcommand{\nemonO}[1]{\fdcase{\NEMONO}{}{#1}}%

\newcommand{\neepi}{\fdcase{\NEEPI}{}{}}%

\newcommand{\Neepi}[1]{\fdcase{\NEEPI}{#1}{}}%

\newcommand{\neepI}[1]{\fdcase{\NEEPI}{}{#1}}%

\newcommand{\nebimo}{\fdcase{\NEBIMO}{}{}}%

\newcommand{\Nebimo}[1]{\fdcase{\NEBIMO}{#1}{}}%

\newcommand{\nebimO}[1]{\fdcase{\NEBIMO}{}{#1}}%

\newcommand{\neiso}{\fdcase{\NEAR}{\hspace{-2pt}\cong}{}}%

\newcommand{\Neiso}[1]{\fdcase{\NEAR}{\hspace{-2pt}\cong}{#1}}%

\newcommand{\nebiar}{\fdbicase{\NEBIAR}{}{}}%

\newcommand{\Nebiar}[2]{\fdbicase{\NEBIAR}{#1}{#2}}%

\newcommand{\neeql}{\fdbicase{\NEEQL}{}{}}%

\newcommand{\neadjar}{\fdbicase{\NEADJAR}{}{}}%

\newcommand{\Neadjar}[2]{\fdbicase{\NEADJAR}{#1}{#2}}%


\newcommand{\nearv}[1]{\fdcase{\NEARV{#100}}{}{}}%

\newcommand{\Nearv}[2]{\fdcase{\NEARV{#200}}{#1}{}}%

\newcommand{\neaRv}[2]{\fdcase{\NEARV{#200}}{}{#1}}%


\newcommand{\SWAR}{\begin{picture}(0,0)%
\put(2900,2900){\vector(-1,-1){5800}}%
\end{picture}}%

\newcommand{\SWDOTAR}%
{\truex{100}\truey{212}%
\NUMBEROFDOTS=5800%
\divide\NUMBEROFDOTS by \value{y}%
\begin{picture}(0,0)%
\multiput(2900,2900)(-\value{y},-\value{y}){\NUMBEROFDOTS}%
{\circle*{\value{x}}}%
\put(-2900,-2900){\vector(-1,-1){0}}%
\end{picture}}%

\newcommand{\SWMONO}{\begin{picture}(0,0)%
\put(2900,2900){\vector(-1,-1){5800}}%
\put(2900,2900){\begin{picture}(0,0)%
\truex{141}%
\put(-\value{x},-\value{x}){\makebox(0,0){$\times$}}%
\end{picture}}\end{picture}}%

\newcommand{\SWEPI}{\begin{picture}(0,0)%
\put(2900,2900){\vector(-1,-1){5800}}%
\put(-2900,-2900){\begin{picture}(0,0)%
\truex{525}%
\put(\value{x},\value{x}){\makebox(0,0){$\times$}}%
\end{picture}}\end{picture}}%

\newcommand{\SWBIMO}{\begin{picture}(0,0)%
\put(2900,2900){\vector(-1,-1){5800}}%
\put(2900,2900){\begin{picture}(0,0)%
\truex{141}%
\put(-\value{x},-\value{x}){\makebox(0,0){$\times$}}%
\end{picture}}%
\put(-2900,-2900){\begin{picture}(0,0)%
\truex{525}%
\put(\value{x},\value{x}){\makebox(0,0){$\times$}}%
\end{picture}}\end{picture}}%

\newcommand{\SWBIAR}{\begin{picture}(0,0)%
\put(2900,2900){\begin{picture}(0,0)%
\truex{141}%
\put(\value{x},-\value{x}){\vector(-1,-1){5800}}%
\put(-\value{x},\value{x}){\vector(-1,-1){5800}}%
\end{picture}}\end{picture}}%

\newcommand{\SWADJAR}{\begin{picture}(0,0)%
\put(-2900,-2900){\begin{picture}(0,0)%
\truex{141}%
\put(-\value{x},\value{x}){\vector(1,1){5800}}%
\end{picture}}%
\put(2900,2900){\begin{picture}(0,0)%
\truex{141}%
\put(\value{x},-\value{x}){\vector(-1,-1){5800}}%
\end{picture}}\end{picture}}%

\newcommand{\SWARV}[1]{\begin{picture}(0,0)%
\put(0,0){\makebox(0,0){\begin{picture}(#1,#1)%
\put(#1,#1){\vector(-1,-1){#1}}\end{picture}}}%
\end{picture}}%


\newcommand{\swar}{\fdcase{\SWAR}{}{}}%

\newcommand{\Swar}[1]{\fdcase{\SWAR}{#1}{}}%

\newcommand{\swaR}[1]{\fdcase{\SWAR}{}{#1}}%

\newcommand{\swdotar}{\fdcase{\SWDOTAR}{}{}}%

\newcommand{\Swdotar}[1]{\fdcase{\SWDOTAR}{#1}{}}%

\newcommand{\swdotaR}[1]{\fdcase{\SWDOTAR}{}{#1}}%

\newcommand{\swmono}{\fdcase{\SWMONO}{}{}}%

\newcommand{\Swmono}[1]{\fdcase{\SWMONO}{#1}{}}%

\newcommand{\swmonO}[1]{\fdcase{\SWMONO}{}{#1}}%

\newcommand{\swepi}{\fdcase{\SWEPI}{}{}}%

\newcommand{\Swepi}[1]{\fdcase{\SWEPI}{#1}{}}%

\newcommand{\swepI}[1]{\fdcase{\SWEPI}{}{#1}}%

\newcommand{\swbimo}{\fdcase{\SWBIMO}{}{}}%

\newcommand{\Swbimo}[1]{\fdcase{\SWBIMO}{#1}{}}%

\newcommand{\swbimO}[1]{\fdcase{\SWBIMO}{}{#1}}%

\newcommand{\swiso}{\fdcase{\SWAR}{\hspace{-2pt}\cong}{}}%

\newcommand{\Swiso}[1]{\fdcase{\SWAR}{\hspace{-2pt}\cong}{#1}}%

\newcommand{\swbiar}{\fdbicase{\SWBIAR}{}{}}%

\newcommand{\Swbiar}[2]{\fdbicase{\SWBIAR}{#1}{#2}}%

\newcommand{\sweql}{\fdbicase{\NEEQL}{}{}}%

\newcommand{\swadjar}{\fdbicase{\SWADJAR}{}{}}%

\newcommand{\Swadjar}[2]{\fdbicase{\SWADJAR}{#1}{#2}}%


\newcommand{\swarv}[1]{\fdcase{\SWARV{#100}}{}{}}%

\newcommand{\Swarv}[2]{\fdcase{\SWARV{#200}}{#1}{}}%

\newcommand{\swaRv}[2]{\fdcase{\SWARV{#200}}{}{#1}}%


\newcommand{\sdcase}[3]{\begin{picture}(0,0)%
\put(0,-150){#1}%
\truex{100}\truez{600}%
\put(\value{x},\value{x}){\makebox(0,\value{z})[l]{${#2}$}}%
\truex{300}\truey{800}%
\put(-\value{x},-\value{y}){\makebox(0,\value{z})[r]{${#3}$}}%
\end{picture}}%

\newcommand{\sdbicase}[3]{\begin{picture}(0,0)%
\put(0,-150){#1}%
\truex{250}\truey{600}\truez{850}%
\put(\value{x},\value{x}){\makebox(0,\value{y})[l]{${#2}$}}%
\put(-\value{x},-\value{z}){\makebox(0,\value{y})[r]{${#3}$}}%
\end{picture}}%


\newcommand{\SEAR}{\begin{picture}(0,0)%
\put(-2900,2900){\vector(1,-1){5800}}%
\end{picture}}%

\newcommand{\SEDOTAR}%
{\truex{100}\truey{212}%
\NUMBEROFDOTS=5800%
\divide\NUMBEROFDOTS by \value{y}%
\begin{picture}(0,0)%
\multiput(-2900,2900)(\value{y},-\value{y}){\NUMBEROFDOTS}%
{\circle*{\value{x}}}%
\put(2900,-2900){\vector(1,-1){0}}%
\end{picture}}%

\newcommand{\SEMONO}{\begin{picture}(0,0)%
\put(-2900,2900){\vector(1,-1){5800}}%
\put(-2900,2900){\begin{picture}(0,0)%
\truex{141}%
\put(\value{x},-\value{x}){\makebox(0,0){$\times$}}%
\end{picture}}\end{picture}}%

\newcommand{\SEEPI}{\begin{picture}(0,0)%
\put(-2900,2900){\vector(1,-1){5800}}%
\put(2900,-2900){\begin{picture}(0,0)%
\truex{525}%
\put(-\value{x},\value{x}){\makebox(0,0){$\times$}}%
\end{picture}}\end{picture}}%

\newcommand{\SEBIMO}{\begin{picture}(0,0)%
\put(-2900,2900){\vector(1,-1){5800}}%
\put(-2900,2900){\begin{picture}(0,0)%
\truex{141}%
\put(\value{x},-\value{x}){\makebox(0,0){$\times$}}%
\end{picture}}%
\put(2900,-2900){\begin{picture}(0,0)%
\truex{525}%
\put(-\value{x},\value{x}){\makebox(0,0){$\times$}}%
\end{picture}}\end{picture}}%

\newcommand{\SEBIAR}{\begin{picture}(0,0)%
\put(-2900,2900){\begin{picture}(0,0)%
\truex{141}
\put(-\value{x},-\value{x}){\vector(1,-1){5800}}%
\put(\value{x},\value{x}){\vector(1,-1){5800}}%
\end{picture}}\end{picture}}%

\newcommand{\SEEQL}{\begin{picture}(0,0)%
\put(-2900,2900){\begin{picture}(0,0)%
\truex{70}%
\put(-\value{x},-\value{x}){\line(1,-1){5800}}%
\put(\value{x},\value{x}){\line(1,-1){5800}}%
\end{picture}}\end{picture}}%

\newcommand{\SEADJAR}{\begin{picture}(0,0)%
\put(-2900,2900){\begin{picture}(0,0)%
\truex{141}%
\put(-\value{x},-\value{x}){\vector(1,-1){5800}}%
\end{picture}}%
\put(2900,-2900){\begin{picture}(0,0)%
\truex{141}%
\put(\value{x},\value{x}){\vector(-1,1){5800}}%
\end{picture}}\end{picture}}%

\newcommand{\SEARV}[1]{\begin{picture}(0,0)%
\put(0,0){\makebox(0,0){\begin{picture}(#1,#1)%
\put(0,#1){\vector(1,-1){#1}}\end{picture}}}%
\end{picture}}%


\newcommand{\sear}{\sdcase{\SEAR}{}{}}%

\newcommand{\Sear}[1]{\sdcase{\SEAR}{#1}{}}%

\newcommand{\seaR}[1]{\sdcase{\SEAR}{}{#1}}%

\newcommand{\sedotar}{\sdcase{\SEDOTAR}{}{}}%

\newcommand{\Sedotar}[1]{\sdcase{\SEDOTAR}{#1}{}}%

\newcommand{\sedotaR}[1]{\sdcase{\SEDOTAR}{}{#1}}%

\newcommand{\semono}{\sdcase{\SEMONO}{}{}}%

\newcommand{\Semono}[1]{\sdcase{\SEMONO}{#1}{}}%

\newcommand{\semonO}[1]{\sdcase{\SEMONO}{}{#1}}%

\newcommand{\seepi}{\sdcase{\SEEPI}{}{}}%

\newcommand{\Seepi}[1]{\sdcase{\SEEPI}{#1}{}}%

\newcommand{\seepI}[1]{\sdcase{\SEEPI}{}{#1}}%

\newcommand{\sebimo}{\sdcase{\SEBIMO}{}{}}%

\newcommand{\Sebimo}[1]{\sdcase{\SEBIMO}{#1}{}}%

\newcommand{\sebimO}[1]{\sdcase{\SEBIMO}{}{#1}}%

\newcommand{\seiso}{\sdcase{\SEAR}{\hspace{-2pt}\cong}{}}%

\newcommand{\Seiso}[1]{\sdcase{\SEAR}{\hspace{-2pt}\cong}{#1}}%

\newcommand{\sebiar}{\sdbicase{\SEBIAR}{}{}}%

\newcommand{\Sebiar}[2]{\sdbicase{\SEBIAR}{#1}{#2}}%

\newcommand{\seeql}{\sdbicase{\SEEQL}{}{}}%

\newcommand{\seadjar}{\sdbicase{\SEADJAR}{}{}}%

\newcommand{\Seadjar}[2]{\sdbicase{\SEADJAR}{#1}{#2}}%


\newcommand{\searv}[1]{\sdcase{\SEARV{#100}}{}{}}%

\newcommand{\Searv}[2]{\sdcase{\SEARV{#200}}{#1}{}}%

\newcommand{\seaRv}[2]{\sdcase{\SEARV{#200}}{}{#1}}%


\newcommand{\NWAR}{\begin{picture}(0,0)%
\put(2900,-2900){\vector(-1,1){5800}}%
\end{picture}}%

\newcommand{\NWDOTAR}%
{\truex{100}\truey{212}%
\NUMBEROFDOTS=5800%
\divide\NUMBEROFDOTS by \value{y}%
\begin{picture}(0,0)%
\multiput(2900,-2900)(-\value{y},\value{y}){\NUMBEROFDOTS}%
{\circle*{\value{x}}}%
\put(-2900,2900){\vector(-1,1){0}}%
\end{picture}}%

\newcommand{\NWMONO}{\begin{picture}(0,0)%
\put(2900,-2900){\vector(-1,1){5800}}%
\put(2900,-2900){\begin{picture}(0,0)%
\truex{141}%
\put(-\value{x},\value{x}){\makebox(0,0){$\times$}}%
\end{picture}}\end{picture}}%

\newcommand{\NWEPI}{\begin{picture}(0,0)%
\put(2900,-2900){\vector(-1,1){5800}}%
\put(-2900,2900){\begin{picture}(0,0)%
\truex{525}%
\put(\value{x},-\value{x}){\makebox(0,0){$\times$}}%
\end{picture}}\end{picture}}%

\newcommand{\NWBIMO}{\begin{picture}(0,0)%
\put(2900,-2900){\vector(-1,1){5800}}%
\put(2900,-2900){\begin{picture}(0,0)%
\truex{141}%
\put(-\value{x},\value{x}){\makebox(0,0){$\times$}}%
\end{picture}}%
\put(-2900,2900){\begin{picture}(0,0)%
\truex{525}%
\put(\value{x},-\value{x}){\makebox(0,0){$\times$}}%
\end{picture}}\end{picture}}%

\newcommand{\NWBIAR}{\begin{picture}(0,0)%
\put(2900,-2900){\begin{picture}(0,0)%
\truex{141}%
\put(\value{x},\value{x}){\vector(-1,1){5800}}%
\end{picture}}%
\put(2900,-2900){\begin{picture}(0,0)%
\truex{141}
\put(-\value{x},-\value{x}){\vector(-1,1){5800}}%
\end{picture}}\end{picture}}%

\newcommand{\NWADJAR}{\begin{picture}(0,0)%
\put(-2900,2900){\begin{picture}(0,0)%
\truex{141}%
\put(\value{x},\value{x}){\vector(1,-1){5800}}%
\end{picture}}%
\put(2900,-2900){\begin{picture}(0,0)%
\truex{141}%
\put(-\value{x},-\value{x}){\vector(-1,1){5800}}%
\end{picture}}\end{picture}}%

\newcommand{\NWARV}[1]{\begin{picture}(0,0)%
\put(0,0){\makebox(0,0){\begin{picture}(#1,#1)%
\put(#1,0){\vector(-1,1){#1}}\end{picture}}}%
\end{picture}}%


\newcommand{\nwar}{\sdcase{\NWAR}{}{}}%

\newcommand{\Nwar}[1]{\sdcase{\NWAR}{#1}{}}%

\newcommand{\nwaR}[1]{\sdcase{\NWAR}{}{#1}}%

\newcommand{\nwdotar}{\sdcase{\NWDOTAR}{}{}}%

\newcommand{\Nwdotar}[1]{\sdcase{\NWDOTAR}{#1}{}}%

\newcommand{\nwdotaR}[1]{\sdcase{\NWDOTAR}{}{#1}}%

\newcommand{\nwmono}{\sdcase{\NWMONO}{}{}}%

\newcommand{\Nwmono}[1]{\sdcase{\NWMONO}{#1}{}}%

\newcommand{\nwmonO}[1]{\sdcase{\NWMONO}{}{#1}}%

\newcommand{\nwepi}{\sdcase{\NWEPI}{}{}}%

\newcommand{\Nwepi}[1]{\sdcase{\NWEPI}{#1}{}}%

\newcommand{\nwepI}[1]{\sdcase{\NWEPI}{}{#1}}%

\newcommand{\nwbimo}{\sdcase{\NWBIMO}{}{}}%

\newcommand{\Nwbimo}[1]{\sdcase{\NWBIMO}{#1}{}}%

\newcommand{\nwbimO}[1]{\sdcase{\NWBIMO}{}{#1}}%

\newcommand{\nwiso}{\sdcase{\NWAR}{\hspace{-2pt}\cong}{}}%

\newcommand{\Nwiso}[1]{\sdcase{\NWAR}{\hspace{-2pt}\cong}{#1}}%

\newcommand{\nwbiar}{\sdbicase{\NWBIAR}{}{}}%

\newcommand{\Nwbiar}[2]{\sdbicase{\NWBIAR}{#1}{#2}}%

\newcommand{\nweql}{\sdbicase{\SEEQL}{}{}}%

\newcommand{\nwadjar}{\sdbicase{\NWADJAR}{}{}}%

\newcommand{\Nwadjar}[2]{\sdbicase{\NWADJAR}{#1}{#2}}%


\newcommand{\nwarv}[1]{\sdcase{\NWARV{#100}}{}{}}%

\newcommand{\Nwarv}[2]{\sdcase{\NWARV{#200}}{#1}{}}%

\newcommand{\nwaRv}[2]{\sdcase{\NWARV{#200}}{}{#1}}%



\newcommand{\ENEAR}[2]%
{\makebox[0pt]{\begin{picture}(0,0)%
\put(0,-150){\makebox(0,0){\begin{picture}(0,0)%
\put(-6600,-3300){\vector(2,1){13200}}%
\truex{200}\truey{800}\truez{600}%
\put(-\value{x},\value{x}){\makebox(0,\value{z})[r]{${#1}$}}%
\put(\value{x},-\value{y}){\makebox(0,\value{z})[l]{${#2}$}}%
\end{picture}}}\end{picture}}}%

\newcommand{\enear}{\ENEAR{}{}}%

\newcommand{\Enear}[1]{\ENEAR{#1}{}}%

\newcommand{\eneaR}[1]{\ENEAR{}{#1}}%

\newcommand{\ESEAR}[2]%
{\makebox[0pt]{\begin{picture}(0,0)%
\put(0,-150){\makebox(0,0){\begin{picture}(0,0)%
\put(-6600,3300){\vector(2,-1){13200}}%
\truex{200}\truey{800}\truez{600}%
\put(\value{x},\value{x}){\makebox(0,\value{z})[l]{${#1}$}}%
\put(-\value{x},-\value{y}){\makebox(0,\value{z})[r]{${#2}$}}%
\end{picture}}}\end{picture}}}%

\newcommand{\esear}{\ESEAR{}{}}%

\newcommand{\Esear}[1]{\ESEAR{#1}{}}%

\newcommand{\eseaR}[1]{\ESEAR{}{#1}}%

\newcommand{\WNWAR}[2]%
{\makebox[0pt]{\begin{picture}(0,0)%
\put(0,-150){\makebox(0,0){\begin{picture}(0,0)%
\put(6600,-3300){\vector(-2,1){13200}}%
\truex{200}\truey{800}\truez{600}%
\put(\value{x},\value{x}){\makebox(0,\value{z})[l]{${#1}$}}%
\put(-\value{x},-\value{y}){\makebox(0,\value{z})[r]{${#2}$}}%
\end{picture}}}\end{picture}}}%

\newcommand{\wnwar}{\WNWAR{}{}}%

\newcommand{\Wnwar}[1]{\WNWAR{#1}{}}%

\newcommand{\wnwaR}[1]{\WNWAR{}{#1}}%

\newcommand{\WSWAR}[2]%
{\makebox[0pt]{\begin{picture}(0,0)%
\put(0,-150){\makebox(0,0){\begin{picture}(0,0)%
\put(6600,3300){\vector(-2,-1){13200}}%
\truex{200}\truey{800}\truez{600}%
\put(-\value{x},\value{x}){\makebox(0,\value{z})[r]{${#1}$}}%
\put(\value{x},-\value{y}){\makebox(0,\value{z})[l]{${#2}$}}%
\end{picture}}}\end{picture}}}%

\newcommand{\wswar}{\WSWAR{}{}}%

\newcommand{\Wswar}[1]{\WSWAR{#1}{}}%

\newcommand{\wswaR}[1]{\WSWAR{}{#1}}%



\newcommand{\NNEAR}[2]%
{\raisebox{-1pt}[0pt][0pt]{\begin{picture}(0,0)%
\put(0,0){\makebox(0,0){\begin{picture}(0,0)%
\put(-3300,-6600){\vector(1,2){6600}}%
\truex{100}\truez{600}%
\put(-\value{x},\value{x}){\makebox(0,\value{z})[r]{${#1}$}}%
\put(\value{x},-\value{z}){\makebox(0,\value{z})[l]{${#2}$}}%
\end{picture}}}\end{picture}}}%

\newcommand{\nnear}{\NNEAR{}{}}%

\newcommand{\Nnear}[1]{\NNEAR{#1}{}}%

\newcommand{\nneaR}[1]{\NNEAR{}{#1}}%

\newcommand{\SSWAR}[2]%
{\raisebox{-1pt}[0pt][0pt]{\begin{picture}(0,0)%
\put(0,0){\makebox(0,0){\begin{picture}(0,0)%
\put(3300,6600){\vector(-1,-2){6600}}%
\truex{100}\truez{600}%
\put(-\value{x},\value{x}){\makebox(0,\value{z})[r]{${#1}$}}%
\put(\value{x},-\value{z}){\makebox(0,\value{z})[l]{${#2}$}}%
\end{picture}}}\end{picture}}}%

\newcommand{\sswar}{\SSWAR{}{}}%

\newcommand{\Sswar}[1]{\SSWAR{#1}{}}%

\newcommand{\sswaR}[1]{\SSWAR{}{#1}}%

\newcommand{\SSEAR}[2]%
{\raisebox{-1pt}[0pt][0pt]{\begin{picture}(0,0)%
\put(0,0){\makebox(0,0){\begin{picture}(0,0)%
\put(-3300,6600){\vector(1,-2){6600}}%
\truex{200}\truez{600}%
\put(\value{x},\value{x}){\makebox(0,\value{z})[l]{${#1}$}}%
\put(-\value{x},-\value{z}){\makebox(0,\value{z})[r]{${#2}$}}%
\end{picture}}}\end{picture}}}%

\newcommand{\ssear}{\SSEAR{}{}}%

\newcommand{\Ssear}[1]{\SSEAR{#1}{}}%

\newcommand{\sseaR}[1]{\SSEAR{}{#1}}%

\newcommand{\NNWAR}[2]%
{\raisebox{-1pt}[0pt][0pt]{\begin{picture}(0,0)%
\put(0,0){\makebox(0,0){\begin{picture}(0,0)%
\put(3300,-6600){\vector(-1,2){6600}}%
\truex{200}\truez{600}%
\put(\value{x},\value{x}){\makebox(0,\value{z})[l]{${#1}$}}%
\put(-\value{x},-\value{z}){\makebox(0,\value{z})[r]{${#2}$}}%
\end{picture}}}\end{picture}}}%

\newcommand{\nnwar}{\NNWAR{}{}}%

\newcommand{\Nnwar}[1]{\NNWAR{#1}{}}%

\newcommand{\nnwaR}[1]{\NNWAR{}{#1}}%



\newcommand{\Necurve}[2]%
{\begin{picture}(0,0)%
\truex{1300}\truey{2000}\truez{200}%
\put(0,\value{x}){\oval(#200,\value{y})[t]}%
\put(0,\value{x}){\makebox(0,0){\begin{picture}(#200,0)%
\put(#200,0){\vector(0,-1){\value{z}}}%
\put(0,0){\line(0,-1){\value{z}}}\end{picture}}}%
\truex{2500}%
\put(0,\value{x}){\makebox(0,0)[b]{${#1}$}}%
\end{picture}}%

\newcommand{\necurve}[1]{\Necurve{}{#1}}%

\newcommand{\Nwcurve}[2]%
{\begin{picture}(0,0)%
\truex{1300}\truey{2000}\truez{200}%
\put(0,\value{x}){\oval(#200,\value{y})[t]}%
\put(0,\value{x}){\makebox(0,0){\begin{picture}(#200,0)%
\put(#200,0){\line(0,-1){\value{z}}}%
\put(0,0){\vector(0,-1){\value{z}}}\end{picture}}}%
\truex{2500}%
\put(0,\value{x}){\makebox(0,0)[b]{${#1}$}}%
\end{picture}}%

\newcommand{\nwcurve}[1]{\Nwcurve{}{#1}}%

\newcommand{\Securve}[2]%
{\begin{picture}(0,0)%
\truex{1300}\truey{2000}\truez{200}%
\put(0,-\value{x}){\oval(#200,\value{y})[b]}%
\put(0,-\value{x}){\makebox(0,0){\begin{picture}(#200,0)%
\put(#200,0){\vector(0,1){\value{z}}}%
\put(0,0){\line(0,1){\value{z}}}\end{picture}}}%
\truex{2500}%
\put(0,-\value{x}){\makebox(0,0)[t]{${#1}$}}%
\end{picture}}%

\newcommand{\securve}[1]{\Securve{}{#1}}%

\newcommand{\Swcurve}[2]%
{\begin{picture}(0,0)%
\truex{1300}\truey{2000}\truez{200}%
\put(0,-\value{x}){\oval(#200,\value{y})[b]}%
\put(0,-\value{x}){\makebox(0,0){\begin{picture}(#200,0)%
\put(#200,0){\line(0,1){\value{z}}}%
\put(0,0){\vector(0,1){\value{z}}}\end{picture}}}%
\truex{2500}%
\put(0,-\value{x}){\makebox(0,0)[t]{${#1}$}}%
\end{picture}}%

\newcommand{\swcurve}[1]{\Swcurve{}{#1}}%



\newcommand{\Escurve}[2]%
{\begin{picture}(0,0)%
\truex{1400}\truey{2000}\truez{200}%
\put(\value{x},0){\oval(\value{y},#200)[r]}%
\put(\value{x},0){\makebox(0,0){\begin{picture}(0,#200)%
\put(0,0){\vector(-1,0){\value{z}}}%
\put(0,#200){\line(-1,0){\value{z}}}\end{picture}}}%
\truex{2500}%
\put(\value{x},0){\makebox(0,0)[l]{${#1}$}}%
\end{picture}}%

\newcommand{\escurve}[1]{\Escurve{}{#1}}%

\newcommand{\Encurve}[2]%
{\begin{picture}(0,0)%
\truex{1400}\truey{2000}\truez{200}%
\put(\value{x},0){\oval(\value{y},#200)[r]}%
\put(\value{x},0){\makebox(0,0){\begin{picture}(0,#200)%
\put(0,0){\line(-1,0){\value{z}}}%
\put(0,#200){\vector(-1,0){\value{z}}}\end{picture}}}%
\truex{2500}%
\put(\value{x},0){\makebox(0,0)[l]{${#1}$}}%
\end{picture}}%

\newcommand{\encurve}[1]{\Encurve{}{#1}}%

\newcommand{\Wscurve}[2]%
{\begin{picture}(0,0)%
\truex{1300}\truey{2000}\truez{200}%
\put(-\value{x},0){\oval(\value{y},#200)[l]}%
\put(-\value{x},0){\makebox(0,0){\begin{picture}(0,#200)%
\put(0,0){\vector(1,0){\value{z}}}%
\put(0,#200){\line(1,0){\value{z}}}\end{picture}}}%
\truex{2400}%
\put(-\value{x},0){\makebox(0,0)[r]{${#1}$}}%
\end{picture}}%

\newcommand{\wscurve}[1]{\Wscurve{}{#1}}%

\newcommand{\Wncurve}[2]%
{\begin{picture}(0,0)%
\truex{1300}\truey{2000}\truez{200}%
\put(-\value{x},0){\oval(\value{y},#200)[l]}%
\put(-\value{x},0){\makebox(0,0){\begin{picture}(0,#200)%
\put(0,0){\line(1,0){\value{z}}}%
\put(0,#200){\vector(1,0){\value{z}}}\end{picture}}}%
\truex{2400}%
\put(-\value{x},0){\makebox(0,0)[r]{${#1}$}}%
\end{picture}}%

\newcommand{\wncurve}[1]{\Wncurve{}{#1}}%



\newcommand{\Freear}[8]{\begin{picture}(0,0)%
\put(#400,#500){\vector(#6,#7){#800}}%
\truex{#200}\truey{#300}%
\put(\value{x},\value{y}){$#1$}\end{picture}}%

\newcommand{\freear}[5]{\Freear{}{0}{0}{#1}{#2}{#3}{#4}{#5}}%


\newcount\SCALE%

\newcount\NUMBER%

\newcount\LINE%

\newcount\COLUMN%

\newcount\WIDTH%

\newcount\SOURCE%

\newcount\ARROW%

\newcount\TARGET%

\newcount\ARROWLENGTH%

\newcount\NUMBEROFDOTS%

\newcounter{x}%

\newcounter{y}%

\newcounter{z}%

\newcounter{horizontal}%

\newcounter{vertical}%

\newskip\itemlength%

\newskip\firstitem%

\newskip\seconditem%

\newcommand{\printarrow}{}%


\newcommand{\truex}[1]{%
\NUMBER=#1%
\multiply\NUMBER by 100%
\divide\NUMBER by \SCALE%
\setcounter{x}{\NUMBER}}%

\newcommand{\truey}[1]{%
\NUMBER=#1%
\multiply\NUMBER by 100%
\divide\NUMBER by \SCALE%
\setcounter{y}{\NUMBER}}%

\newcommand{\truez}[1]{%
\NUMBER=#1%
\multiply\NUMBER by 100%
\divide\NUMBER by \SCALE%
\setcounter{z}{\NUMBER}}%

\newcommand{\changecounters}[1]{%
\SOURCE=\ARROW%
\ARROW=\TARGET%
\settowidth{\itemlength}{#1}%
\ifdim \itemlength > 2800\unitlength%
\addtolength{\itemlength}{-2800\unitlength}%
\TARGET=\itemlength%
\divide\TARGET by 1310%
\multiply\TARGET by 100%
\divide\TARGET by \SCALE%
\else%
\TARGET=0%
\fi%
\ARROWLENGTH=5000%
\advance\ARROWLENGTH by -\SOURCE%
\advance\ARROWLENGTH by -\TARGET%
\advance\SOURCE by -\TARGET}%

\newcommand{\initialize}[1]{%
\LINE=0%
\COLUMN=0%
\WIDTH=0%
\ARROW=0%
\TARGET=0%
\changecounters{#1}%
\renewcommand{\printarrow}{#1}%
\begin{center}%
\vspace{10pt}%
\begin{picture}(0,0)}%

\newcommand{\DIAG}[1]{%
\SCALE=100%
\setlength{\unitlength}{655sp}%
\initialize{\mbox{$#1$}}}%

\newcommand{\DIAGV}[2]{%
\SCALE=#1%
\setlength{\unitlength}{655sp}%
\multiply\unitlength by \SCALE%
\divide\unitlength by 100%
\initialize{\mbox{$#2$}}}%

\newcommand{\n}[1]{%
\changecounters{\mbox{$#1$}}%
\put(\COLUMN,\LINE){\makebox(0,0){\printarrow}}%
\thinlines%
\renewcommand{\printarrow}{\mbox{$#1$}}%
\advance\COLUMN by 4000}%

\newcommand{\nn}[1]{%
\put(\COLUMN,\LINE){\makebox(0,0){\printarrow}}%
\thinlines%
\ifnum \WIDTH < \COLUMN%
\WIDTH=\COLUMN%
\else%
\fi%
\advance\LINE by -4000%
\COLUMN=0%
\ARROW=0%
\TARGET=0%
\changecounters{\mbox{$#1$}}%
\renewcommand{\printarrow}{\mbox{$#1$}}}%

\newcommand{\conclude}{%
\put(\COLUMN,\LINE){\makebox(0,0){\printarrow}}%
\thinlines%
\ifnum \WIDTH < \COLUMN%
\WIDTH=\COLUMN%
\else%
\fi%
\setcounter{horizontal}{\WIDTH}%
\setcounter{vertical}{-\LINE}%
\end{picture}}%

\newcommand{\diag}{%
\conclude%
\raisebox{0pt}[0pt][\value{vertical}\unitlength]{}%
\hspace*{\value{horizontal}\unitlength}%
\vspace{10pt}%
\end{center}%
\setlength{\unitlength}{1pt}}%

\newcommand{\diagv}[3]{%
\conclude%
\NUMBER=#1%
\rule{0pt}{\NUMBER pt}%
\hspace*{-#2pt}%
\raisebox{0pt}[0pt][\value{vertical}\unitlength]{}%
\hspace*{\value{horizontal}\unitlength}
\NUMBER=#3%
\advance\NUMBER by 10%
\vspace*{\NUMBER pt}%
\end{center}%
\setlength{\unitlength}{1pt}}%

\newcommand{\N}[1]%
{\raisebox{0pt}[7pt][0pt]{$#1$}}%

\newcommand{\movename}[3]{%
\hspace{#2pt}%
\raisebox{#3pt}[5pt][2pt]{\raisebox{#3pt}{$#1$}}%
\hspace{-#2pt}}%

\newcommand{\movearrow}[3]{%
\makebox[0pt]{%
\hspace{#2pt}\hspace{#2pt}%
\raisebox{#3pt}[0pt][0pt]{\raisebox{#3pt}{$#1$}}}}%

\newcommand{\movevertex}[3]{%
\mbox{\hspace{#2pt}%
\raisebox{#3pt}{\raisebox{#3pt}{$#1$}}%
\hspace{-#2pt}}}%

\newcommand{\crosslength}[2]{%
\settowidth{\firstitem}{#1}%
\settowidth{\seconditem}{#2}%
\ifdim\firstitem < \seconditem%
\itemlength=\seconditem%
\else%
\itemlength=\firstitem%
\fi%
\divide\itemlength by 2%
\hspace{\itemlength}}%

\newcommand{\cross}[2]{%
\crosslength{\mbox{$#1$}}{\mbox{$#2$}}%
\begin{picture}(0,0)%
\put(0,0){\makebox(0,0){$#1$}}%
\thinlines%
\put(0,0){\makebox(0,0){$#2$}}%
\thinlines%
\end{picture}%
\crosslength{\mbox{$#1$}}{\mbox{$#2$}}}%

\newcommand{\B}{\thicklines}


\newcommand{\adj}{\begin{picture}(9,6)%
\put(1,3){\line(1,0){6}}\put(7,0){\line(0,1){6}}%
\end{picture}}%

\newcommand{\com}{\begin{picture}(12,8)%
\put(6,4){\oval(8,8)[b]}\put(6,4){\oval(8,8)[r]}%
\put(6,8){\vector(-1,0){2}}\end{picture}}%

\newcommand{\Nat}[3]{\raisebox{-2pt}%
{\begin{picture}(34,15)%
\put(2,10){\vector(1,0){30}}%
\put(2,0){\vector(1,0){30}}%
\put(13,2){$\Downarrow$}%
\put(20,3){$\scriptstyle{#2}$}%
\put(4,11){$\scriptstyle{#1}$}%
\put(4,1){$\scriptstyle{#3}$}%
\end{picture}}}%

\newcommand{\nat}{\raisebox{-2pt}%
{\begin{picture}(34,10)%
\put(2,10){\vector(1,0){30}}%
\put(2,0){\vector(1,0){30}}%
\put(13,2){$\Downarrow$}%
\end{picture}}}%

\newcommand{\Binat}[5]{\raisebox{-7.5pt}%
{\begin{picture}(34,25)%
\put(2,20){\vector(1,0){30}}%
\put(2,10){\vector(1,0){30}}%
\put(2,0){\vector(1,0){30}}%
\put(13,12){$\Downarrow$}%
\put(13,2){$\Downarrow$}%
\put(20,13){$\scriptstyle{#2}$}%
\put(20,3){$\scriptstyle{#4}$}%
\put(4,21){$\scriptstyle{#1}$}%
\put(4,11){$\scriptstyle{#3}$}%
\put(4,1){$\scriptstyle{#5}$}%
\end{picture}}}%

\newcommand{\binat}{\raisebox{-7.5pt}%
{\begin{picture}(34,20)%
\put(2,20){\vector(1,0){30}}%
\put(2,10){\vector(1,0){30}}%
\put(2,0){\vector(1,0){30}}%
\put(13,12){$\Downarrow$}%
\put(13,2){$\Downarrow$}%
\end{picture}}}%
 
\setcounter{section}{0}
\noindent
\maketitle
\abstract{For any minimal compact complex surface $S$ with $n=b_2(S)>0$ containing global spherical shells (GSS) we study the effectiveness of the $2n$ parameters given by the $n$ blown up points. There exists a family of surfaces $\cal S\to B$ with GSS which contains as fibers $S$, some Inoue-Hirzebruch surface and non minimal surfaces, such that blown up points are generically effective parameters. These families are versal outside a non empty hypersurface $T\subset B$. We deduce that, for any configuration of rational curves, there is a non empty open set in the Oeljeklaus-Toma moduli space such that the corresponding surfaces are defined by a contracting germ in Cremona group, in particular admit a birational structure.}
\tableofcontents

\section{Introduction}
Hopf surfaces are defined by contracting invertible germs $F:(\bb C^2,0)\to (\bb C^2,0)$. There are well-known normal forms
$$F(z_1,z_2)=(az_1+tz_2^m,bz_2), \quad 0<|a|\le |b|<1,\ (a-b^m)t=0,\ m\in\bb N^\star,$$
which give effective parameters of the versal deformation and give charts with transition mappings in the group $Aut(\bb C^2)$ of polynomial automorphisms of $\bb C^2$, in particular in the Cremona group $Bir(\bb P^2(\bb C))$ of birational mappings of $\bb P^2(\bb C)$. Hopf surfaces are particular cases of a larger family of compact complex surfaces in the VII$_0$ class of Kodaira, namely surfaces $S$ containing global spherical shells (GSS). When $b_2(S)\ge 1$, these surfaces  also called Kato surfaces admit neither affine nor projective stuctures \cite{IKO, KO,Kl}. Their explicit construction  consists in the composition $\Pi$ of $n=b_2(S)$ blowing-ups (depending on $2n$ parameters)  followed by a special glueing by a germ of isomorphism $\s$ (depending on an infinite number of parameters). These surfaces are not almost homogeneous \cite{P} hence $0\le \dim H^0(S,\T)\le 1$ and Chern classes  of surfaces in class VII$_0$ satisfy the conditions $b_2(S)=c_2(S)=-c_1^2(S)$. By Riemann-Roch formula, we obtain the dimension of the base of the versal deformation of $S$,
$$2n\le \dim H^1(S,\T)=2b_2(S)+\dim H^0(S,\T)\le 2n+1,$$
where $\T$ is the sheaf of holomorphic vector fields.\\
Some questions are raised
\begin{description}
\item{(1)} Are the $2n$ parameters of the blown up points effective parameters ?
\item{(2)} If there are non trivial global vector fields, there is at least one missing parameter. How to choose it ?
\item{(3)} Do compact surfaces with GSS admit a birational structure, i.e. is there an atlas with transition mappings in Cremona group $Bir(\bb P^2(\bb C))$. More precisely is there in each conjugation class of contracting germs of the form $\Pi\s$ (or of strict germs, following Favre terminology \cite{Fav}) a birational representative ?
\end{description}
{\bf Known results}: 
\begin{itemize}
\item If $S$ is a Enoki surface (see \cite{DK})  known normal forms, namely
$$F(z_1,z_2)=\bigl(t^nz_1z_2^n+\sum_{i=0}^{n-1}a_it^{i+1}z_2^{i+1}, tz_2\bigr),\quad 0<|t|<1,$$
 are birational. The parameters $t$ and $a_i$, $i=0,\ldots,n-1$ are effective at $a=(a_0,\ldots,a_{n-1})=0$ (i.e. if there are global vector fields or $S$ is a Inoue surface) and give the versal deformation. If $a\neq 0$, there is no global vector fields, $t$ and $a_i$, $i=0,\ldots,n-1$ but one $a_j\neq 0$ give the versal deformation. The complex numbers $a_i$ are the coordinates of the blown up points $O_i$ in the successive exceptional curves $C_i$.
 \item If $S$ is a Inoue-Hirzebruch surface (see \cite{D2})  
$$N(z_1,z_2)=(z_1^pz_2^q,z_1^rz_2^s),$$
is the composition of blowing-ups hence is birational.
 Here  $\left(\begin{array}{cc}p&q\\r&s\end{array}\right)\in Gl(2,\bb Z)$ is the composition of matrices 
$$\left(\begin{array}{cc}1&1\\0&1\end{array}\right) \quad {\rm or}\quad \left(\begin{array}{cc}0&1\\1&1\end{array}\right)$$
 with at least one of the second type. There is no parameters because these surfaces are logarithmically rigid.
 \item If $S$ is of intermediate type (see definition in section 2), there are normal forms due to C.Favre \cite{Fav} 
$$F(z_1,z_2)=(\l z_1z_2^{\frak s}+P(z_2),z_2^k), \qquad \l\in\bb C^\star,\ \frak s\in\bb N^\star,\  k\ge 2,$$
where $P$ is a special polynomial. These normal forms are adapted to logarithmic deformations and show the existence of a foliation, however {\it are not birational}.  In \cite{OT} K.Oeljeklaus and M.Toma explain how to recover second Betti number which is now hidden and give coarse moduli spaces of surfaces with fixed intersection matrix,
\item Some special cases of intermediate surfaces are obtained from Hénon mappings $H$ or composition of Hénon mappings. More precisely,  the germ of $H$ at the fixed point at infinity is strict, hence yields a surface with a GSS \cite{HO, DO2}. These germs are birational.
\end{itemize}
\vspace{2mm}

{\bf Motivation}: A.Teleman \cite{Te2,Te3} proved  that for $b_2(S)=1,2$, any minimal surface in class VII$_0^+$ contains a cycle of rational curves, therefore has a deformation into a surface with GSS. In order to prove that any surface in class VII$_0^+$ contains a GSS, we should solve the following

{\bf Problem}: Let $\cal S\to \D$ be a family of compact surfaces over the disc such that for every $u\in\D^\star$, $S_u$ contains a GSS. Does $S_0$ contain a GSS ? In other words, are the surface with GSS  closed in families ?\\
To solve this problem we have to study families of surfaces in which curves do not fit into flat families, the volume of some curves in these families may be not uniformly bounded (see \cite{DT1}) and configurations of curves change. Favre normal forms of polynomial germs associated to surfaces with GSS, cannot be used because the discriminant of the intersection form is fixed. Moreover, if using the algorithm in \cite{OT} we put $F$ under the form $\Pi\s$, $\s$ is not fixed in the logarithmic family, depends on the blown up points and degenerates when a generic blown up point approaches the intersection of two curves.\\
Therefore this article and \cite{D6}, section 5, focus on the problem of finding effective parameters and new normal forms of contracting germs in intermediate cases of surfaces with fixed $\s$, such that surfaces are minimal or not and intersection matrices are not fixed. Since usual holomorphic objects, curves or foliations, do not fit in global family, it turns out that birational structures could be the adapted notion. Clearly the problem of their unicity raises. \\
\vspace{2mm}\\
{\bf Main results}: A marked surface $(S,C_0)$ is a surface $S$ with a fixed rational curve $C_0$. In section 2, we define large families $\Phi_{J,\s}:\cal S_{J,\s}\to B_J$ of marked surfaces with GSS with fixed second Betti number $n=b_2(S)$ which use the same $n$ charts of blowing-ups identified by a subset $J\subset\{0,\ldots,n-1\}$. The base admits a stratification by strata over which the intersection matrix of the $n$ rational curves is fixed. With these fixed charts, we construct explicit global sections of the direct image sheaf of the vertical vector fields $R^1\Phi_{J,\s,\star}\T$ over $B_J$, which express the dependence on the parameters of the blown up points: $[\t_i]$ are the infinitesimal deformations along the rational curves and $[\mu_i]$, $i=0,\ldots,n-1$ the infinitesimal deformations transversaly to the rational curves. Surfaces with non trivial global vector fields exist over an analytic set of codimension at least 2. The choice of a rational curve $C_0$ (the first one in the construction) fixes the conjugacy class of a contracting germ. Using a result by A.Teleman \cite{Te4}, we obtain in section 3,
 \begin{Th} Let $(S,C_0)$ be a minimal marked surface containing a GSS of intermediate type,  with $n=b_2(S)$.  Let $J=I_\infty(C_0)$ be the indices of the blown up points at infinity and let $\Phi_{J,\s}:\cal S_{J,\s}\to B_J$ be the family of surfaces with GSS associated to $J$ and $\s$. Then, there exists a non empty hypersurface $T_{J,\s}\subset B_J$ containing $Z=\{u\in B\mid h^0(S_u,\T_u)>0\}$ such that for $u\in B_J\setminus T_{J,\s}$,
\begin{description}
\item{a)}  $\{[\t^i_u],[\mu^i_u] \mid 0\le i\le n-1\}$ is a base of $H^1(S_u,\T_u)$,
\item{b)} $\{[\t^i_u] \mid O_i\ {\rm is\ generic}\}$ is a base of   the space of infinitesimal logarithmic deformations $H^1(S_u,\Theta_u(-Log\ D_u))$, where $D_u$ is the maximal divisor in $S_u$.
\end{description}
Moreover 
\begin{description}
\item{i)} If $T_{J,\s}$ intersects a stratum $B_{J,M}$ then $T_{J,\s}\cap B_{J,M}$ is a hypersurface in $B_{J,M}$,
\item{ii)} $T_{J,\s}$  intersect each stratum $B_{J,M}$ such that the corresponding surfaces admit twisted vector fields and $Z\cap B_{J,M}\subset T_{J,\s}$,
\end{description}
\end{Th}
\begin{Cor} Any marked surface $(S,C_0)$ belongs to a large family $\Phi_{J,\s}:\cal S_{J,\s}\to B_J$ and there is a non empty hypersurface $T_{J,\s}$ such that over $B_J\setminus T_{J,\s}$ this family is versal.
\end{Cor}
This answers to the first question and the result is the best possible because $T_{J,\s}$ is never empty. What happens on the hypersurface $T_{J,\s}$ ? Is it possible that there is a curve of isomorphic surfaces ? Is the canonical image of a stratum $B_{J,M}$ in the Oeljeklaus-Toma coarse moduli space open ?  Do we obtain all possible surfaces ?\\
We give a partial answer to the question (3):
\begin{Cor} Let $M$ be any intersection matrix of a minimal compact complex surface containing a GSS (i.e. of a Kato surface) then the O-T moduli space of such surfaces contains a non empty open set of surfaces admitting a birational structure.
\end{Cor}

 A complete answer is given in \cite{D6} section 5, if there is only one branch attached to the cycle.\\

This article stems from discussions with Karl Oeljeklaus and Matei Toma at the university of Osnabrück about the case $b_2(S)=2$, I thank them for their relevant remarks. I thank Andrei Teleman for fruitful discussions in particular to have pointed out  that thanks to his results \cite{Te4} the cocycles $\t_i$ and $\mu_i$ {\it cannot} be everywhere independent.

\section{Surfaces with Global Spherical Shells}
\subsection{Basic constructions}
\begin{Def} Let $S$ be a compact complex surface. We say that $S$ contains a global spherical shell, if there is a biholomorphic map $\f:U\to S$ from a neighbourhood $U\subset \bb C^2\setminus\{0\}$ of the sphere $S^3$ into $S$ such that $S\setminus \f(S^3)$ is connected.
\end{Def}
Hopf surfaces are the simplest examples of surfaces with GSS.\\

Let $S$ be a surface containing a GSS with $n=b_2(S)$. It is known that $S$ contains $n$ rational curves and to each curve it is possible to associate a contracting germ of mapping $F=\Pi\s=\Pi_0\cdots\Pi_{n-1}\s:(\bb C^2,0) \to (\bb C^2,0)$ where $\Pi=\Pi_0\cdots\Pi_{n-1}:B^\Pi\to B$ is a sequence of $n$ blowing-ups  and $\s$ is a germ of isomorphism (see \cite{D1}). The surface is obtained by gluing  two open shells as explained by the following picture
\begin{center}
\includegraphics[width=12cm]{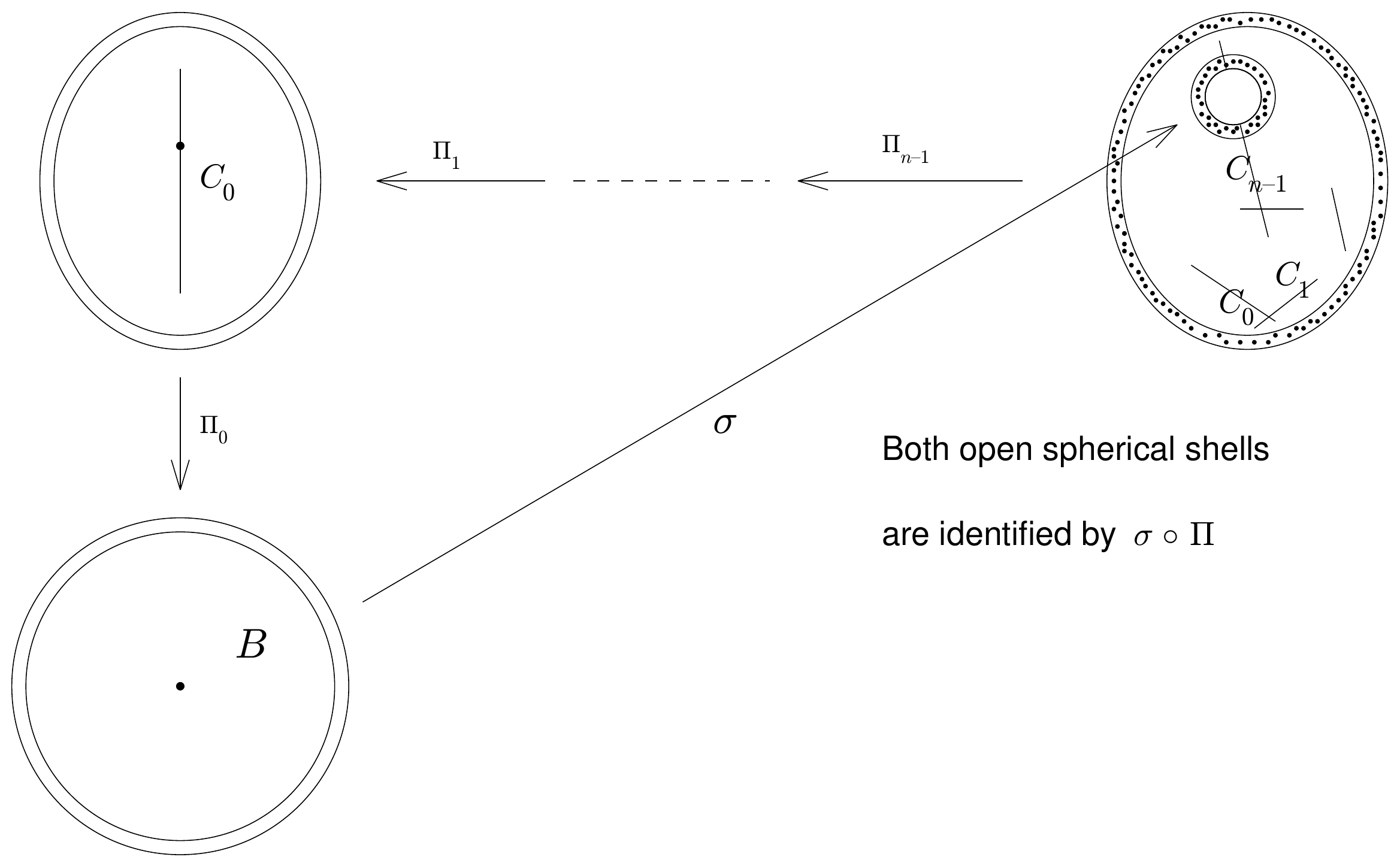}
\end{center}

\begin{Def}\label{Un-arbre-RecEnoki} Let $S$ be a  surface containing a GSS, with $n=b_2(S)$. A {\bf Enoki covering} of $S$ is an open covering $\cal U=(U_i)_{0\le i\le n-1}$ obtained in the following way:
\begin{itemize}
\item $W_0$ is the  ball of radius $1+\e$  blown up at the origin, $C_0=\Pi_0^{-1}(0)$, $B'_0\subset \subset B_0$ are small balls centered at $O_0=(a_0,0)\in W_0$, $U_0=W_0\setminus B'_0$,
\item For $1\le i\le n-1$, $W_i$ is the  ball $B_{i-1}$  blown up at $O_{i-1}$, $C_i=\Pi_i^{-1}(O_{i-1})$, $B'_{i}\subset \subset B_{i}$  are small balls centered at $O_{i}\in W_i$, $U_i=W_i\setminus B'_i$.
\end{itemize}
The  pseudoconcave boundary of $U_i$ is patched with the pseudoconvex boundary of $U_{i+1}$ by $\Pi_i$, for $i=0,\ldots,n-2$ and the pseudoconcave boundary of $U_{n-1}$ is patched with the pseudoconvex boundary of $U_0$ by $\s\Pi_0$, where
$$\begin{array}{cccc}
\s:&B(1+\e)&\to&W_{n-1}\\
&z=(z_1,z_2)&\mapsto&\s(z)
\end{array}$$
is biholomorphic on its image, satisfying $\s(0)=O_{n-1}$. 
\end{Def}

If we want to obtain a minimal surface, the sequence of blowing-ups has to be made in the following way:
\begin{itemize}
\item $\Pi_0 $ blows up the origin of the two dimensional unit ball $B$, 
\item $\Pi_1$ blows up a point $O_0\in C_0=\Pi_0^{-1}(0)$,\ldots 
\item $\Pi_{i+1}$ blows up a point $O_{i}\in C_{i}=\Pi_{i}^{-1}(O_{i-1})$, for $i=0,\ldots,n-2$, and 
\item $\s:\bar B\to B^\Pi$ sends isomorphically a neighbourhood of $\bar B$ onto a small ball in $B^\Pi$ in such a way that $\s(0)\in C_{n-1}$. 
\end{itemize}
Each $W_i$ is covered by two charts with coordinates $(u_i,v_i)$ and $(u'_i,v'_i)$ in which $\Pi_i$ writes $\Pi_i(u_i,v_i)=(u_iv_i+a_{i-1},v_i)$ and $\Pi_i(u'_i,v'_i)=(v'_i+a_{i-1},u'_iv'_i)$. In these charts the exceptional curves has always the equations $v_i=0$ and $v'_i=0$.\\
A  blown up point $O_i\in C_i$ will be called {\bf generic} if it is not at the intersection of two curves.
The data $(S,C)$ of a surface $S$ and of a rational curve in $S$ will be called a {\bf marked surface}.

If we assume that $S$ is minimal and that we are in the intermediate case,  there is at least one blowing-up at a generic point, and one at the intersection of two curves (hence $n\ge 2$). If there is only one tree i.e. one regular sequence and if we choose $C_0$ as being the curve which induces the root of the tree, we suppose that
\begin{itemize}
\item
 $\Pi_1$ is a generic blowing-up, 
 \item $\Pi_{n-1}$ blows-up the intersection of $C_{n-2}$ with another rational curve and
 \item $\s(0)$ is one of the two intersection points of $C_{n-1}$ with the previous curves.
 \end{itemize}
The  Enoki covering is obtained as in the following picture:
\begin{center}\label{Enoki}
\includegraphics[width=12cm]{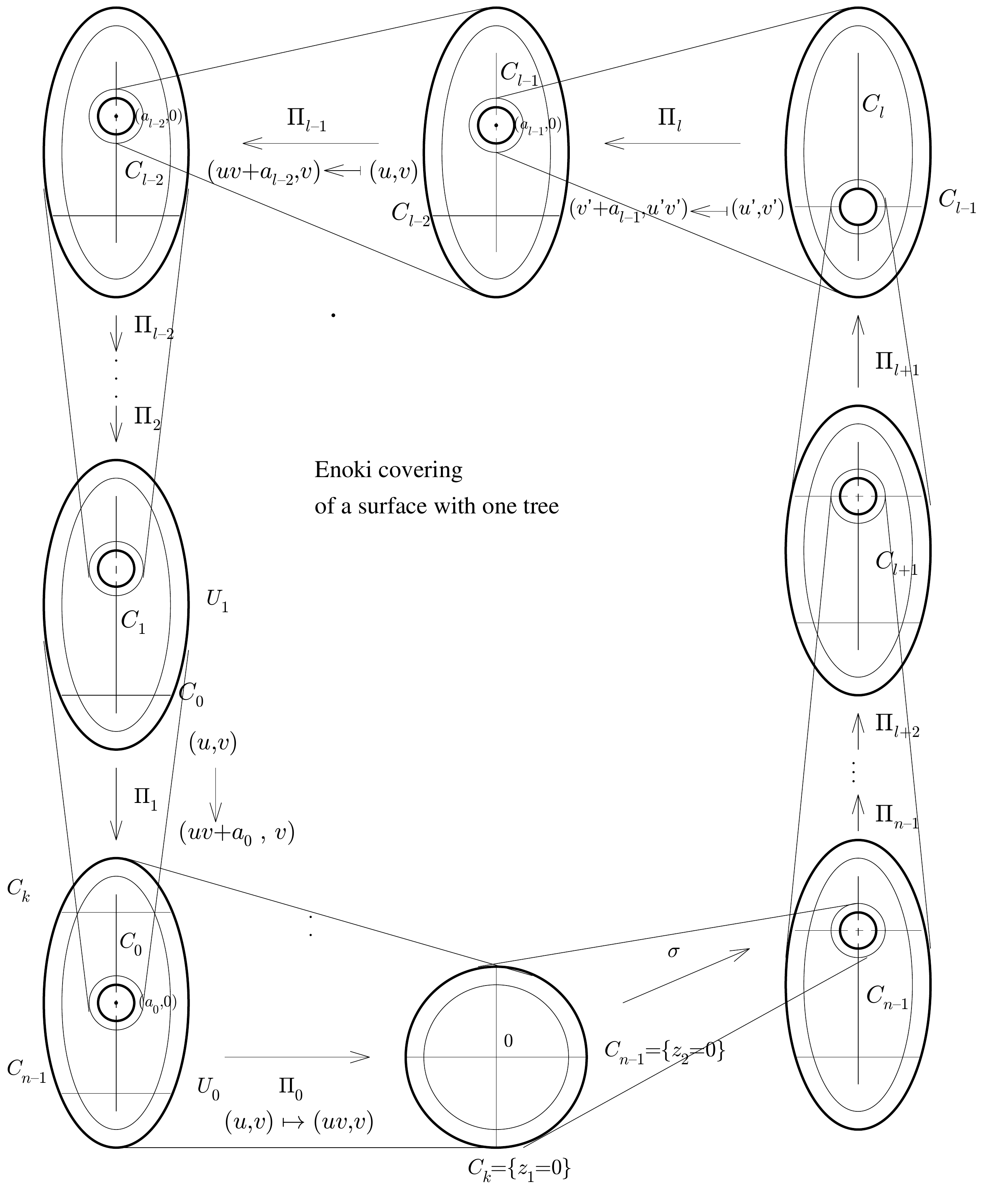}
\end{center}
where 
\begin{itemize}
\item $1\le l\le n-1$ and $n\ge 2$. If all, but one, blowing-ups are generic, then $l=n-1$ 
\item For $i=1,\ldots,l-1$, $\Pi_i(u_i,v_i)=(u_iv_i+a_{i-1},v_i)$ are generic blowing-ups,
\item $\Pi_{l}(u'_{l},v'_{l})=(v'_{l}+a_{l-1},u'_{l}v'_{l})$ is also generic, but $O_l$ is the origin of the chart $(u'_l,v'_l)$,
\item For $i=l+1,\ldots,n-1$, $\Pi_i(u_i,v_i)=(u_iv_i,v_i)$ or $\Pi_i(u'_i,v'_i)=(v'_i,u'_iv'_i)$ are blowing-ups at the intersection of two curves.
\end{itemize}

The general case of $\rho\geq 1$ trees is obtained by joining $\rho$ sequences similar to the previous one, i.e.,
$$\begin{array}{cll}
F&= &\Pi\s\\
&&\\
& = &(\Pi_0\cdots\Pi_{l_0-1}\Pi_{l_0}\cdots\Pi_{n_1-1})\cdots \\
&&\\
&&(\Pi_{n_1+\cdots+n_\k}\cdots\Pi_{n_1+\cdots+n_\k+l_{\k}-1}\Pi_{n_1+\cdots+n_\k+l_{\k}}\cdots\Pi_{n_1+\cdots+n_\k+n_{\k+1}-1}) \cdots\\ 
&&\\
&&(\Pi_{n_1+\cdots+n_{\rho -1}}\cdots\Pi_{n_1+\cdots+n_{\rho -1}+l_{\rho-1}-1}\Pi_{n_1+\cdots+n_{\rho -1}+l_{\rho-1}}\cdots\Pi_{n_1+\cdots+n_{\rho}-1}) \s.
\end{array}$$
where $n_1+\cdots+n_{\rho}=n$.

\begin{center}
\includegraphics[width=8cm]{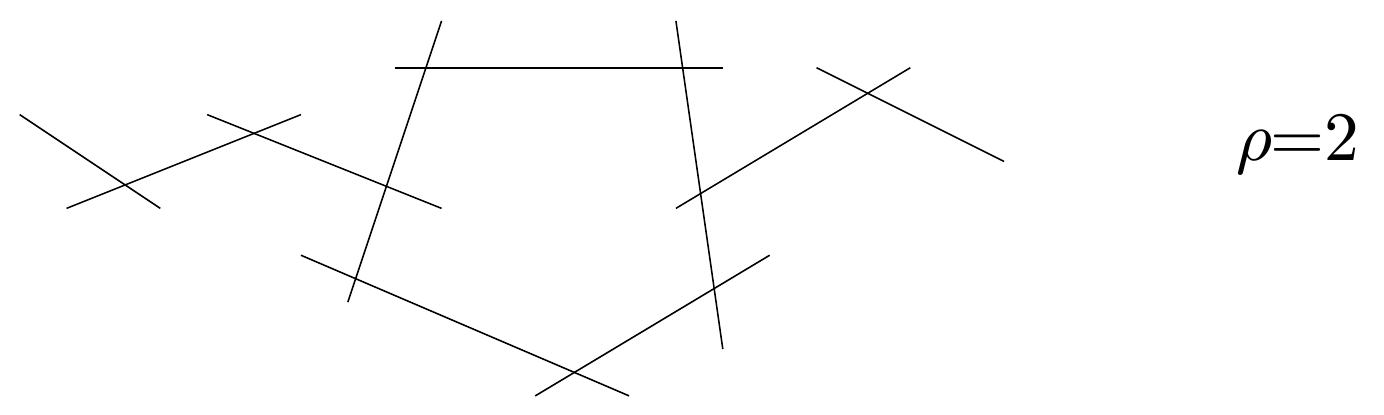}
\end{center}
We may suppose, up to a conjugation of $F$ by a linear map, that 
$$\partial_1\s_2(0)=\frac{\partial \s_2}{\partial z_1}(0)=0$$
 it means that the strict transform of the curve $\s^{-1}(C_{n-1})$ intersects $C_0$  at the infinite point of the chart $(u,v)$, i.e. the origin of $(u',v')$. This condition is convenient for computations.\\
When $n=2$, we denote by $U_{01}= U_0\cap \Pi_1(U_1)\subset U_0$ and $U_{10}= U_1\cap \s\Pi_0(U_0)\subset U_1$ the two connected components of the intersection  $U_0\cap U_1$ of the images in $S$ of $U_0$ and $U_1$, denoted in the same way.\\
If $n\geq 3$, $U_{i,i+1}= U_i\cap \Pi_{i+1}(U_{i+1})$, $i=0,\ldots,n-2$, $U_{n-1,0}= U_{n-1}\cap \s\Pi_0(U_0)$.\\
We refer to \cite{D1} for the description of configurations of curves. We index the curves $(C_i)_{i\in\bb Z}$ in the universal covering space following the canonical order  (see \cite{D1}). Let $a(S)=(a_i)_{i\in\bb Z}$ be the family of positive integers defined by $a_i=-C_i^2$. By \cite{D1} p104, this family is periodic of period $n$ and for any index $i\in\bb Z$ we define a positive integer independant of $i$,
$$2n\le \s_n(S):=\sum_{j=i}^{i+n-1}a_i\le 3n.$$
The family $(a_i)_{i\in\bb Z}$ splits into sequences
$$s_p=(p+2,2,\ldots,2) \quad{\rm and} \quad r_m=(2,\ldots,2)$$
 of length respectively $p$ and $m$, where $p\ge 1$ and $m\ge 1$. We call $s_p$ (resp. $r_m$), $p\ge 1$ ($m\ge 1$) the singular (resp. regular) sequence of length $p$ (resp. $m$). We have
$$\rho:= \sharp\{\rm trees\}=\sharp\{\rm regular\ sequences\}.$$

\subsection{Large families  of marked surfaces\label{logfamilies}}
With the previous notations, we consider global families of minimal compact surfaces with the same charts, parameterized by the coordinates of the blown up points on the successive exceptional curves obtained in the construction of the surfaces and such that any marked surface with GSS $(S,C_0)$ belongs to at least one of these families.  More precisely,  let $F(z)=\Pi_0\cdots\Pi_{n-1}\s(z)$ be a germ associated to any marked surface $(S,C_0)$ with $tr(S)=0$. In order to fix the notations we suppose that $C_0=\Pi_0^{-1}(0)$ meets two other curves (see the picture after definition \thesection.\ref{Un-arbre-RecEnoki}), hence $\s(0)$ is the intersection of $C_{n-1}$ with another curve. We suppose that
 $$\part_1\s_2(0)=0.$$
  We denote by $I_\infty(C_0)\subset \{0,\ldots,n-1\}$ the subset of indices which correspond to blown up points at infinity, that is to say,
$$I_\infty(C_0):=\bigl\{i\mid O_i\ {\rm is\ the\ origin\ of\ the\ chart}\ (u'_i,v'_i)\bigr\}.$$
Each generic blow-up
$$\Pi_i(u_i,v_i)=(u_iv_i+a_{i-1},v_i)\quad {\rm or}\quad \Pi_i(u'_i,v'_i)=(v'_i+a_{i-1},u'_iv'_i)$$
 may be deformed moving the blown up point $(a_{i-1},0)$. If we do not want to change the configuration  we take
 $$\begin{array}{l}
 {\rm for\  all}\  \k=0,\ldots,\rho-1\quad {\rm (with }\ n_0=0{\rm)},\\
 \\
 \hspace{20mm}\left\{
 \begin{array}{l}
 a_{n_1+\cdots+n_\kappa}\in \bb C^\star,\\
 \\
 \forall i,\  1\le i\le l_{\k} -1,\qquad  a_{n_1+\cdots+n_\kappa +i}\in \bb C, \\
 \\
  \forall j,\  0\le j\le n_{\k+1}-l_{\k}-1,\qquad a_{n_1+\cdots+n_\kappa+l_{\k}+j}=0.
  \end{array}\right.
  \end{array}$$
  The mapping $\s$ is supposed to be fixed.  We obtain a large family of compact surfaces which contains $S$ such that all the surfaces $S_a$ have the same intersection matrix $$M=M(S_a)=M(S),$$
   therefore is a logarithmic deformation. For $J=I_\infty(C_0)$ we denote 
  this family
$$\Phi_{J,M,\s}:\cal S_{J,M,\s}\to B_{J,M}$$
where
$$\begin{array}{l}
B_{J,M}\\
\\
:=\bb C^\star\times \bb C^{l_0-1}\times\{0\}^{n_1-l_0}\times\cdots\times \bb C^\star\times \bb C^{l_\kappa -1}\times\{0\}^{n_{\k+1}-l_\k}\times\cdots\times \bb C^\star\times\bb C^{l_{\rho-1}-1}\times\{0\}^{n_\rho-l_{\rho-1}}\\
\\
\simeq\bb C^\star\times \bb C^{l_0-1}\times\cdots\times \bb C^\star\times \bb C^{l_\kappa -1}\times\cdots\times \bb C^\star\times\bb C^{l_{\rho-1}-1}\end{array}$$
and $n_1+\cdots+n_\rho=n$.\\

In $\cal S_{J,M,\s}$ there is a flat family of divisors $\cal D\subset \cal S$ with irreducible components 
$$\cal D_i, \quad i=0,\ldots,n-1,$$
such that for every $a\in B_{J,M}$, $M=(D_{i,a}.D_{j,a})_{0\le i,j\le n-1}$.
  We may extend this family towards smaller or larger strata which produce minimal surfaces:
  \begin{itemize}
  \item On one hand, {\bf towards a unique Inoue-Hirzebruch surface}:
  Over 
  $$\bb C^{l_0}\times\{0\}^{n_1-l_0}\times\cdots\times \bb C^{l_\kappa}\times\{0\}^{n_{\k+1}-l_\k}\times\cdots\times \bb C^{l_{\rho-1}}\times\{0\}^{n_\rho-l_{\rho-1}}\simeq \bb C^{l_0}\times\cdots\times \bb C^{l_\kappa}\times\cdots\times \bb C^{l_{\rho-1}},$$ 
$$\Phi_{J,\s}:\cal S_{J,\s}\to \bb C^{l_0}\times\cdots\times \bb C^{l_\kappa}\times\cdots\times \bb C^{l_{\rho-1}}.$$
If for an index  $\kappa$, $a_{n_1+\cdots+n_\kappa}=0$, there is a jump in the configuration of the curves. For instance, if for all $\kappa$, $\kappa=0,\ldots,\rho-1$
 $$a_{n_1+\cdots+n_\kappa}=\cdots=a_{n_1+\cdots+n_\kappa+l_\kappa-1}=0$$
 we obtain a Inoue-Hirzebruch surface.  To be more precise the base 
 $$\bb C^{l_0}\times\cdots\times \bb C^{l_\kappa}\times\cdots\times \bb C^{l_{\rho-1}}$$
 splits into locally closed submanifolds called {\bf strata}
\begin{itemize}
\item the Zariski open set $\bb C^\star\times\bb C^{l_0-1}\times\cdots\times \bb C^\star\times\bb C^{l_\k-1}\times\cdots\bb C^\star\times\bb C^{l_{\rho-1}-1}$,
\item $\rho=C^1_\rho$ codimension one strata 
$$\bb C^\star\times\bb C^{l_0-1}\times\cdots\times \{0\}\times \bb C^\star\times\bb C^{l_\k-2}\times\cdots\times \bb C^\star\times\bb C^{l_{\rho-1}-1}, \quad 0\le \k\le \rho-1,$$
\item $C^p_{\rho+p-1}$ codimension $p$ strata, $1\le p:=p_0+\cdots +p_{\rho-1}\le l_0+\cdots+l_{\rho-1}$,
$$\{0\}^{p_0}\times \bb  C^\star\times \bb C^{l_0-p_0-1}\times\cdots\times\{0\}^{p_\k}\times\bb C^\star \times \bb C^{l_\k-p_\k-1}\times\cdots\times\{0\}^{p_{\rho-1}}\times \bb C^\star\times \bb C^{l_{\rho-1}-p_{\rho-1}-1}$$
\end{itemize}
 \item On second hand, {\bf towards Enoki surfaces}. If for all indices such that $O_i$ is at the intersection of two rational curves, in particular for $i\in J$, the blown up point $O_i$  is moved to $O_i=(a_i,0)$ with $a_i\neq 0$, all the blown up points become generic, the trace of the contracting germ is different from $0$. We obtain also all the intermediate configurations.\\
\begin{Prop} \label{towardsEnoki} There is a monomial holomorphic function $t:\bb C^{\Card J}\to \bb C$ depending on the variables $a_j$, $j\in J$ such that over $B_J:=\{|t(a)|<1\}\subset \bb C^n$, the family $\Phi_{J,\s}:\cal S_{J,\s}\to B_{J}$ may be extended and for every $a\in B_J$, $t(a)=\tr(S_a)$. 
\end{Prop}

Proof: The trace of a surface does not depend on the germs associated to this surface therefore we may suppose that $O_0=(a_0,0)$ is in the chart $(u'_0,v'_0)$, i.e. $0\in J$.\\
Suppose that $\Card J=1$, then for $i\neq 0$, $\Pi_i(u_i,v_i)=(u_iv_i+a_{i-1},v_i)$ and 
$$\s(z)=(\s_1(z)+a_{n-1},\s_2(z)).$$
 We have
$$\begin{array}{lcl}
F(z)&=&\dps\Pi\s(z)=\Pi_0\bigl(\s_1(z)\s_2(z)^{n-1}+\sum_{j=0}^{n-1}a_j\s_2(z)^j,\s_2(z)\bigr)\\
&&\\
&=&\dps\Bigl(\s_2(z), \s_1(z)\s_2(z)^{n}+\sum_{j=0}^{n-1}a_j\s_2(z)^{j+1}\Bigr),\end{array}\leqno{(\spadesuit)}$$
and with our convention on $\s$,
$$\tr DF(0)=\tr \left(\begin{array}{cc}\part_1\s_2(0)&\part_2\s_2(0)\\a_0\part_1\s_2(0)&a_0\part_2\s_2(0)\end{array}\right)= \tr \left(\begin{array}{cc}0&\part_2\s_2(0)\\0&a_0\part_2\s_2(0)\end{array}\right)=a_0\part_2\s_2(0).$$
The general case is obtained by the composition $F=F_1\circ\cdots\circ F_N$, where $N=\Card J$,  $F_N$  of the type of $(\spadesuit)$ 
$$F_N(z)=\Bigl(\s_2(z), \s_1(z)\s_2(z)^{m_N}+\sum_{j=0}^{m_N-1}a^N_j\s_2(z)^{j+1}\Bigr), \quad m_N\ge 1$$
and other $F_k$ have similar expressions with $\s=Id$ and $m_k\ge 1$, i.e.
$$F_k(u,v)=\Bigl(v,uv^{m_k}+\sum_{j=0}^{m_k-1} a^k_jv^{j+1}\Bigr)$$
with $m_1+\cdots+m_N=n$. Therefore
$$F(z)=\bigl(\star\ ,\ \part_2\s_2(0)a^1_0a^2_0\cdots a^N_0 z_2\bigr)\qquad {\rm mod}\  (z_1,z_2)^2$$
and $\tr DF(0)=\part_2\s_2(0) a^1_0a^2_0\cdots a^N_0$.\hfill$\Box$\\

Now, $B_{J}\subset \bb C^n$ is an open neighbourhood of 
$$ \bb C^{l_0}\times\{0\}^{n_1-l_0}\times\cdots\times  \bb C^{l_\kappa}\times\{0\}^{n_{\k+1}-l_\k}\times\cdots\times \bb C^{l_{\rho-1}}\times\{0\}^{n_\rho-l_{\rho-1}}.$$
 and we extend the family
 $$\Phi_{J,\s}:\cal S_{J,\s}\to B_J.$$
 thanks to proposition \thesection.\ref{towardsEnoki}. We obtain larger strata of minimal surfaces, from dimension $l+1$ to dimension $n$.
 \end{itemize}
 
 \begin{Ex} 1) {\bf \boldmath Example with 2 curves}: For $(3,2)=-(C_0^2,C_1^2)$, $J=\{1\}$,  $O_0=(a_0,0)$ and $O_1=(a_1,0)$ with $a_0\in\bb C^\star$, $a_1=0$.
 The stratum of Inoue-Hirzebruch surface $(4,2)$ is obtained for $a_0=0$, and generic surfaces are obtained for $a_0\in\bb C$, $a_1\neq 0$. If $\s(z_1,z_2)=(z_1+a_1,z_2)$, 
 $$F(z)=\Pi\s(z)=\bigl(z_2(z_1+a_1)(z_2+a_0), z_2(z_1+a_1)\bigr)$$
 $tr DF(0)=a_1$, hence $B_J=\bb C\times\D$.\\
 
 2)  {\bf \boldmath Example with 6 curves}: If we start with the sequence 
$$(42\ 2\ 3\ 3\ 2)=(s_2r_1s_1s_1r_1)=-(C_0^2,C_1^2,C_2^2,C_3^2,C_4^2,C_5^2)$$
$J=\{1,4,5\}$, and the blown up points are
$O_i=(a_i,0)$, $i=0,\ldots,5$ with 
$$a_0\in\bb C^\star,\quad a_1=0,\quad a_2=0,\quad a_3\in\bb C^\star,  \quad a_4=0, \quad a_5=0.$$
 Strata towards Inoue-Hirzebruch surfaces are
\begin{itemize}
\item $(522\ 3\ 3\ 2)$ when $a_3=0$, 
\item $(42\ 2\ 3\ 42)$, when $a_0=0$,
\item $(522\ 3\ 42)$, when $a_0=a_3=0$,
 which is a Inoue-Hirzebruch surface with one cycle.
\end{itemize}
Towards Enoki surfaces, we move each non generic point into generic one:
\begin{itemize}
\item $(3\ 22\ 3\ 3\ 2)$ with $a_1=0$, $a_2\in\bb C^\star$,
\item $(222\ 3\ 3\ 2)$ with $a_1\in\bb C^\star$,
\item $(42\ 22\ 3\ 2)$ with $a_4\in\bb C^\star$,
\item $(42\ 2\ 3\ 22)$ with $a_5\in\bb C^\star$,
\item $(3\ 222\ 3\ 2)$ with $a_1=0$, $a_2\in\bb C^\star$, $a_4\in\bb C^\star$,
\item $(2222\ 3\ 2)$ with  $a_1\in\bb C^\star$, $a_4\in\bb C^\star$,
\item $(3\ 22\ 3\ 22)$ with $a_1=0$, $a_2\in\bb C^\star$, $a_5\in\bb C^\star$,
\item $(222\ 3\ 22)$ with  $a_1\in\bb C^\star$, $a_5\in\bb C^\star$,
\item $(42\ 2222)$ with $a_4\in\bb C^\star$, $a_5\in\bb C^\star$,
\item $(3\ 22222)$ with $a_1=0$, $a_2\in\bb C^\star$, $a_4\in\bb C^\star$, $a_5\in\bb C^\star$,
\item $(222222)$ with  $a_1\in\bb C^\star$, $a_4\in\bb C^\star$, $a_5\in\bb C^\star$
\end{itemize}
\end{Ex}

\begin{Prop} \label{minimalcodim2} For any $J\neq\emptyset$, any invertible germ $\s:(\bb C^2,0)\to (\bb C^2,0)$, any point $u$ in  the stratum $\{t\neq 0\}\subset B_J$,
$$H^0(S_u,\T_u)=0,$$
i.e. there is no Inoue surface in the stratum $\{t\neq 0\}\subset B_J$.
\end{Prop}
Proof: 1) Let $N=\Card J\ge 1$, we may suppose that the numbering is chosen so that $0\in J$, then $J=\{j_1=0,j_2,\ldots,j_N\}$, and $G=\Pi\circ \s$ is the composition $G=G_1\circ\cdots\circ G_N$, where
$$G_i(z)=\Pi_{j_i}\circ\cdots\circ\Pi_{j_{i+1}-1},\quad i=1,\ldots, N-1$$
$$G_N(z)=\Pi_{j_N}\circ\cdots\circ\Pi_{n-1}\circ\s$$
For each index $j_i$, 
$$\Pi_{j_i}(u',v')=(v'+a_{j_i-1},u'v'), \quad a_{j_1-1}=a_{-1}=0,$$
and for all the other indices $k$
$$\Pi_k(u,v)=(uv+a_{k-1},v)$$
therefore
$$G_i(z)=\left(v,uv^{l_i+1}+\sum_{k=0}^{l_i} a_k^iv^{k+1}\right),\quad i=1,\ldots, N-1$$
$$G_N(z)=\left(\s_2(z),\s_1(z)\s_2(z)^{l_N+1}+\sum_{k=0}^{l_N}a_k^N\s_2(z)^{k+1}\right)$$
$$\s(z)=\bigl(\s_1(z)+a_{l_N}^N,\s_2(z)\bigr).$$
A simple proof by induction shows that
$$G(z)=G_1\circ\cdots\circ G_N(z)=\Bigl(a_0^2\cdots a_0^N\s_2(z)+\cdots, a_0^1\cdots a_0^N\s_2(z)+\cdots\Bigr)$$
2)  We suppose that $G$ is a germ associated to a Inoue surface therefore there is an invertible germ $\f:(\bb C^2,0)\to (\bb C^2,0)$ such that $F\circ \f=\f\circ G$ where
$$F(z_1,z_2)=(t^n z_1z_2^n,tz_2),$$
is the canonical germ associated to Inoue surfaces and  $n=\sum_{i=1}^N l_i +N$. Recall that a germ associated to a Inoue surface admits an invariant curve $\G$, i.e. there is a unique germ of curve such that $F_{\mid\G}:\G\to\G$ is a contraction and the curve $\G$ induces the elliptic curve of the associated surface $S(G)\simeq S(F)$ \cite{D1}. The curves $\G$ and $\s^{-1}(C_{n-1})$ are transversal, therefore replacing if necessary $G$ by a conjugate, we may suppose that $\part_1\s_2(0)=\part_2\s_1(0)=0$. Since the invariant curve of $F$ is $\{z_1=0\}$ and its contracted curve is $\{z_2=0\}$, the germ $\f(z)=(\f_1(z),\f_2(z))$ has a diagonal linear part, i.e. $\part_2\f_1(0)=\part_1\f_2(0)=0$.
 
We have now
$$F\circ\f(z)=\left(t^n\f_1(z)\f_2(z)^{n},t\f_2(z)\right)$$
$$\f\circ G(z)=\Bigl(\f_1\bigl(a_0^2\cdots a_0^N\s_2(z)+\cdots, a_0^1\cdots a_0^N\s_2(z)+\cdots\bigr), \ast\Bigr)$$
Since 
$$\f_1(z)=A_{10}z_1\quad {\rm mod}\quad (z_1,z_2)^2$$ 
with $A_{10}\neq 0$, the left member of the equality $\f\circ G=F\circ\f$ is, 
$$t^{n}\f_1(z)\f_2(z)^{n}= A_{10}a_0^2\cdots a_0^N\part_2\s_2(0)z_2+\cdots$$
which is impossible since there is no linear term in the left member.\hfill$\Box$

\vspace{5mm}

 Remain non minimal surfaces: we still extend the previous family on a small neighbourhood $\wh{B_J}$ of $B_J$, moving the blown up point transversally to the exceptional curves $C_i=\{v_i=0\}\cup\{v'_i=0\}$, introducing $n$ new parameters
 $$\Pi_i(u_i,v_i)=(u_iv_i+a_{i-1},v_i+b_{i-1}),  \quad {\rm or}\quad \Pi_i(u'_i,v'_i)=(v'_i+a_{i-1},u'_iv'_i+b_{i-1}), \quad |b_{i-1}|<<1,$$
 we obtain 
 $$\wh{\Phi}_{J,\s}:\wh{\cal S}_{J,\s}\to \wh{B}_J,$$
 with $\dim \wh{B_J}=2n=2b_2$.
 Since for any $(a,b)\in \wh{B_J}$, $h^1(S_{a,b},\T_{a,b})=2b_2(S_{a,b})+h^0(S_{a,b},\T_{a,b})$, there are some questions:
 \begin{itemize}
 \item Are the parameters $a_i,b_i$, $i=0,\ldots,n-1$, effective ?
 \item Which parameter to add when $h^1(S_{a,b},\T_{a,b})=2b_2(S_{a,b})+1$ in order to obtain a complete family ?
\item If we choose $\s=Id$ or more generally an invertible polynomial mapping, we obtain a birational polynomial germs. Does this families contain all the isomorphy classes of surfaces with fixed intersection matrix $M$ ?  
\end{itemize}

\begin{Rem} It is difficult to determine the maximal domain $\wh{B}_J$ over which $\wh{\Phi}_{J,\s}$ may be defined. When the surface is minimal, i.e. when $b=(b_0,\ldots,b_{n-1})=0$, $F_{a,b}(0)=0$. However, when $b\neq 0$, the fixed point $\z=(\zeta_1,\zeta_2)$  moves and the existence condition for the corresponding surface is that the eigenvalues $\l_1$ and $\l_2$ of $DF_{a,b}(\z)$ satisfy $|\l_i|<1$, $i=1,2$.
\end{Rem}

\subsection{Minimal and non minimal deformations}

Let $S=S(F)$ be a minimal surface with GSS and $\cal U=(U_{i,i+1})$ a Enoki covering of $S$. We denote by $(e_i)_{0\le i\le n-1}$ the base of the free $\bb Z$-module $H_2(S,\bb Z)$ which trivializes the intersection form, i.e. $e_i.e_j=-\d_{ij}$. Here a {\bf simply minimal divisor} is a connected divisor which may be blown down on a regular point.

\begin{Prop}\label{nonminimaliteHi} Let $\wh{\Phi}_{J,\s}:\wh{\cal S}_{J,\s}\to \wh{B}_J$ be a large family of marked surfaces with GSS. Then for any $i=0,\ldots,n-1$ there exists 
\begin{itemize}
\item A smooth hypersurface $H_{i}\subset \wh{B}_J$,
\item A flat family of divisors $\Phi_{J,\s}:\cal E_{i}\to \wh{B}_J\setminus H_{i}$,
\end{itemize}
such that
\begin{enumerate}
\item For any $(a,b)\in \wh{B}_J\setminus H_{i}$, $E_{i,(a,b)}$ is a simply exceptional divisor such that 
$$[E_{i,(a,b)}]=e_i,$$
\item  $S_{(a,b)}$ contains a simply exceptional divisor $E_{i,(a,b)}$ such that $[E_{i,(a,b)}]=e_i$ \iff $(a,b)\not\in H_{i}$,
\item Any intersection $H_{i_1} \cap\cdots\cap H_{i_p}$ of $p$ different such hypersurfaces is smooth of codimension $p$.
\end{enumerate}
\end{Prop}
Proof: The fundamental remark is that $(a,b)\in H_i$ \iff in the construction of the surface $S_{(a,b)}$ there is a sequence of indices $i,i+j_1,\ldots,i+j_p=i$ mod $n$ such that the curve $C_{i+j_k}$ is blown up by $C_{i+j_{k+1}}$. If this sequence of blow-ups ends before reaching the index $i$, say at $i+j_q$, $C_i+C_{i+j_1}+\cdots+C_{i+j_q}$ will be a simply exceptional divisor. Therefore, the total transform of $C_i$ has to check
$$O_{i-1}=(a_{i-1},b_{i-1})\in\Pi_{i-1}^{-1}\cdots\Pi_0^{-1}\s^{-1}\Pi_{n-1}^{-1}\cdots\Pi_{i+1}^{-1}(C_i),$$
or equivalently
$$\Pi_{i+1}\circ\cdots\circ\Pi_{n-1}\circ\s\circ\Pi_0\circ\cdots\circ\Pi_{i-1}(a_{i-1},b_{i-1})\in C_i=\{v_i=0\}.$$
We have 
$$\Pi_{i+1}(u_{i+1},v_{i+1})=(u_{i+1}v_{i+1}+a_i,v_{i+1}+b_i)\quad {\rm or}\quad
\Pi_{i+1}(u'_{i+1},v'_{i+1})=(v'_{i+1}+a_i,u'_{i+1}v'_{i+1}+b_i)$$
therefore the condition the equation of $H_{i}$  is
$$b_i+P(a_0,b_0,\ldots,a_{i-1},b_{i-1},a_{i+1},b_{i+1},\ldots,a_{n-1},b_{n-1},\wt{\s_1}(a_{i-1},b_{i-1}),\wt{\s_2}(a_{i-1},b_{i-1}))=0$$
where 
\begin{itemize}
\item $P$ is a polynomial,
\item $\wt\s(a_{i-1},b_{i-1})= \s\circ\Pi_0\cdots\circ\Pi_{i-1}(a_{i-1},b_{i-1})$ does not depend on $b_i$,
\end{itemize}
\ldots and this is the equation of a smooth hypersurface. The third assertion follows readily from these equations.\hfill$\Box$

\section{Infinitesimal deformations of surfaces with GSS}
 
\subsection{Infinitesimal deformations of the families $\cal S_{J,\s}$}
We define the following cocycles which are the infinitesimal deformations of the families $\wh{\cal S}_{J,\s}\to \wh{B}_J$:
\begin{itemize}
\item For  $i=0,\ldots, n-1$, the cocycles $\t^ i$ called the ``tangent cocycles'' move the blown up points $O_i$ along the curve $C_i$ and vanish only (at order two) at the point ``at infinity''  $C_i\cap C_{i-1}$,
\item For  $i=0,\ldots, n-1$, the cocycle $\mu^ i$ called the ``tranversal cocycles'' move $O_i$ transversaly to $C_i$
\end{itemize}
On a stratum where there are global twisted vector fields we need another infinitesimal deformation (see \cite{D6} for an explicit construction). \\

More precisely,
$$\begin{array}{ccll}
 \theta^i & = & \left\{\begin{array}{ccc}
\displaystyle \frac{\partial}{\partial u_i} &{\rm on}& U_{i,i+1}\\
&&\\
 0 &{\rm over}& U_{j,j+1},\;j\neq i
 \end{array}\right. & \;{\rm If\;}O_i \; {\rm belongs\; to\;the\;chart\;} (u_i,v_i) \\
&&&\\
  \theta^i & = & \left\{\begin{array}{ccc}
\displaystyle  \frac{\partial}{\partial u'_i} &{\rm on}& U_{i,i+1}\\
&&\\
 0 &{\rm over}& U_{j,j+1},\;j\neq i
 \end{array}\right. & \begin{array}{l}  {\rm If\;}O_i \; {\rm belongs\; to\;the\;chart\;} (u'_i,v'_i)\\ {\rm in\ particular\ if\ }O_i=C_i\cap C_{i-1}\end{array}
\end{array}$$

Since $\theta^i$ just moves the blown up point $O_i$ along the curve $C_i$, all surfaces in these deformations are minimal.\\
 We introduce now $n$ other cocycles which move the blown up point $O_i$ transversaly to the  exceptional curves $C_i$. They yield non minimal surfaces, for instance blown up Hopf surfaces but also surfaces with GSS blown up $k$ times, $1\le k\le n$.\\
 For  $i=0,\ldots, n-1$, 
$$\begin{array}{ccll}
\mu^ i & = & \left\{\begin{array}{ccc}
\displaystyle \frac{\partial}{\partial v_i} &{\rm on}& U_{i,i+1}\\
&&\\
 0 &{\rm over}& U_{j,j+1},\;j\neq i
 \end{array}\right. & \;{\rm If\;}O_i \; {\rm belongs\; to\;the\;chart\;} (u_i,v_i) \\
&&&\\
  \mu^ i & = & \left\{\begin{array}{ccc}
\displaystyle  \frac{\partial}{\partial v'_i} &{\rm on}& U_{i,i+1}\\
&&\\
 0 &{\rm over}& U_{j,j+1},\;j\neq i
 \end{array}\right. & \begin{array}{l} {\rm If\;}O_i \; {\rm belongs\; to\;the\;chart\;} (u'_i,v'_i)\\ {\rm in\ particular\ if\ }O_i=C_i\cap C_{i-1}\end{array}
\end{array}$$  

For any $J\subset \{0,\ldots,n-1\}$, the family $\wh{\cal S}_{J,\s}\to \wh{B}_J$  is globally endowed with a family of Enoki coverings.
 Using the family of Enoki coverings, all the cocycles $\t^i$,  $\mu^i$ are globally defined over $\wh{\cal S}_{M,\s}$ and give global sections
 $$[\t^i]\in H^0\bigl(\wh{B}_J,R^1\wh{\Phi}_{J,\s\star}\T\bigr), \quad \quad [\mu^i]\in H^0\bigl(\wh{B}_J,R^1 \wh{\Phi}_{J,\s\star}\T\bigr), \quad i=0,\ldots,n-1$$
 and more precisely for the $l$ indices $i$ such that $O_i$ is generic
 $$[\t^i]\in H^0\bigl(B_{J,M},R^1\Phi_{J,M,\s\star}(\T(-Log\ \cal D))\bigr).$$
For any $(a,b)\in \wh{B}_J$, the cocycles $[\t^i(a,b)], [\mu^i(a,b)]\in R^1\Phi_\star\T_{(a,b)}\ot\bb C=H^1(S_{(a,b)},\T_{(a,b)})$, $i=0,\ldots,n-1$ are  infinitesimal deformations at $(a,b)\in \wh{B}_{J}$ associated to the family $\wh{\cal S}_{J,\s}\to \wh{B}_J$.

\subsection{Splitting of the space of infinitesimal deformations}
 We  divide minimal deformation in two types of deformations: logarithmic deformations for which the intersection matrix of the maximal divisor $D$ does not change, in particular the surfaces remain minimal, and deformations in which the cycle may be smoothed at some singular points or disappear and surfaces  may become non minimal.
\begin{Th}\label{splitting} Let $S$ be a minimal surface containing a GSS with $b_2(S)=n\ge 1$ rational curves $D_0,\ldots,D_{n-1}$ such that $M(S)$ is negative definite. Let $U$ be a spc neighbourhood of $D$, $\rho$ the number of trees in $D$, $r_{l_0},\ldots, r_{l_{\rho-1}}$ the corresponding regular sequences and 
$$l=\sum_{i=0}^{\rho-1} l_i$$
the sum of the length of the regular sequences which is also the number of generic blow-ups.
 Then we have the exact sequence
$$0\to H^1(S,\T_S(-\log D)) \to H^1(S,\T_S)\to H^1(U,\T_{\mid U}) \to 0.\leqno{(\ast)}$$
Moreover
$$\dim H^1(S,\T(-\log D)) = l+\dim H^0(S,\T_S)=3b_2(S)-\s_n(S)+\dim H^0(S,\T_S),$$ 
$$\dim H^1(U,\T_{\mid U})=2b_2(S)-l=\s_n(S)-b_2(S).$$
\end{Th}
Proof: Consider the exact sequence on $S$
$$0\to \T_S(-\log D)\to \T_S\to J_D\to 0\leqno{(\maltese)}$$
 where 
 $$J_D:=\T_S/\T_S(-\log D)=\bigoplus_{i=0}^{n-1}N_{D_i},$$
  $Supp(J_D)=D$, and $N_{D_i}$ the normal bundle of $D_i$. The long exact sequence of cohomology gives
 $$\cdots \to H^0(D,J_D)\to H^1(S,\T_S(-\log D)) \to H^1(S,\T_S)\to H^1(D,J_D) \to H^2(S,\T_S(-\log D)) \to \cdots$$
 If $\t\in H^0(D,J_D)$ its restriction $\t_{D_i}$ to each curve $D_i$ is a section  in the normal bundle $N_{D_i}$ of $D_i$. Since $D_i^2\le -2$, $H^0(D_i,N_{D_i})=0$, hence $\t=0$ and $H^0(D,J_D)=0$. Moreover, by \cite{N2}, thm (1.3), $H^2(S,\T_S(-\log D))=0$, therefore we have
 $$0\to H^1(S,\T_S(-\log D)) \to H^1(S,\T_S)\to H^1(D,J_D) \to 0.\leqno{(\ast)}$$
 We compute now $H^1(D,J_D)$: the restriction of $(\maltese)$ to $U$ gives
 $$0\to H^1(U,\T_U(-\log D)) \to H^1(U,\T_U)\to H^1(D,J_D) \to 0$$
 since by Siu theorem $H^2(U,\T_U(-\log D))=0$.\\ 
 Besides, denoting by $C$ the cycle of rational curves and by $H=D-C$ the sum of trees which meet $C$,  we have the exact sequence
 $$0\to \T_U(-\log D)\to \T_U(-\log C)\to J_H\to 0$$
 where $J_H:=\T_U(-\log C)/\T_U(-\log D)$ and $Supp(J_H)\subset H$. \\
 By \cite{N1} lemma (4.3), $H^1(U,\T_U(-\log C))=0$, and $H^0(H,J_H)=0$, hence
 $$ H^1(U,\T_U(-\log D))=0$$
 With $(\ast)$ we conclude.\\
 By \cite{B} (see appendix \ref{Bruasse}), $h^1(S,\T(-\log D)) =3b_2(S)-\s_n(S)+h^0(S,\T)$. Moreover $3b_2(S)-\s_n(S)$ is the number of generic blown up points $O_i$  and also is equal to the sum of lengths of regular sequences. \hfill $\Box$

\subsection{Infinitesimal non logarithmic deformations \label{section4.2}}
We would like to show that $[\t^ 0],\ldots,[\t^ {n-1}]$,  $[\m^ 0],\ldots,[\m^ {n-1}]$  are generically linearly independent. We suppose that there exists a linear relation
$$\sum_{i=0}^{n-1}(\a_i[\t^ i]+\b_i[\m^ i])=0.$$

We choose the curve $C_0$ such that $O_0$ is a generic point but $O_{n-1}$ is the intersection of two curves. Hence $D_0$ the curve in $S$ induced by $C_0$ is the root of a tree. We shall use this fact later. We have the following linear system where $X_i$ is a vector field over $U_i$, $i=0,\ldots,n-1$:
$$\left\{\begin{array}{lclr}
X_0-{\Pi_1}_\star X_1& =& \displaystyle \a_0\frac{\partial}{\partial u''_0}+\b_0\frac{\partial}{\partial v''_0} &{\rm on} \quad U_{01}\subset U_0\\
\hspace{10mm}\vdots&&\hspace{20mm}\vdots&\\
X_{i}-{\Pi_{i+1}}_\star X_{i+1}&=&\displaystyle \a_{i}\frac{\partial}{\partial u''_{i}}+\b_{i}\frac{\partial}{\partial v''_{i}}  &{\rm on} \quad U_{i,i+1}\subset U_i\\
\hspace{10mm}\vdots&&\hspace{20mm}\vdots&\\
X_{n-2}-{\Pi_{n-1}}_\star X_{n-1}&=&\displaystyle \a_{n-2}\frac{\partial}{\partial u''_{n-2}}+\b_{n-2}\frac{\partial}{\partial v''_{n-2}}  &{\rm on} \quad U_{n-2,n-1}\subset U_{n-2}\\
&&&\\
X_{n-1}- (\s\Pi_0)_\star X_0 & = &\displaystyle \a_{n-1}\frac{\partial}{\partial u''_{n-1}}+\b_{n-1}\frac{\partial}{\partial v''_{n-1}}&  {\rm on} \quad U_{n-1,0}\subset U_{n-1}
\end{array}\right.\leqno{(E1)}$$
where,  $u''_i=u_i$ or $u''_i=u'_i$ (resp. $v''_i=v_i$ or $v''_i=v'_i$).\\
We notice that by Hartogs theorem, $X_i$ extends to $W_i$,  hence $X_i$ is tangent to $C_i$  for $i=0,\ldots,n-1$; moreover 
$$ {\Pi_1}_\star X_1(O_0)=\cdots = {\Pi_{n-1}}_\star X_{n-1}(O_{n-2})= (\s\Pi_0)_\star X_0 (O_{n-1})=0.$$
 Therefore the $i$-th equation at $O_{i}$ gives $\b_i=0$. 
\begin{Rem} If we replace the vector field $\frac{\part}{\part v_i}$ by any non vanishing transversal vector field, the proof works as well.
\end{Rem}

 Now, we show that if $O_i$ is the intersection point of two curves, then $\a_i=0$.
 In fact, there are two cases:
 \begin{description}
 \item {\bf First case}  \boldmath$O_i=C_i\cap C_{i-1}$\unboldmath: In the $(i-1)$-th equation, $X_{i-1}$ and $\frac{\partial}{\partial u_{i-1}}$ or $\frac{\partial}{\partial u'_{i-1}}$ are defined on whole $W_{i-1}$, therefore it is the same for ${\Pi_{i}}_\star X_{i}$, so ${\Pi_{i}}_\star X_{i}$ is tangent to $C_{i-1}$. As consequence, $X_{i}$ is tangent to (the strict transform of) $C_{i-1}$ in $W_{i}$, thus $X_{i}$ vanishes at the intersection point $O_i=C_i\cap C_{i-1}$. We have   $$X_{i}(O_{i})={\Pi_{i+1}}_\star X_{i+1}(O_{i})=0,$$
 hence $\a_{i}=0$.
 \item {\bf Second case} \boldmath$O_i=C_i\cap C_{k}$, $k<i-1$\unboldmath : Then we have $O_{k+1}=C_{k+1}\cap C_k$ and by the previous case,
 \begin{description}
 \item{(1)} $\a_{k+1}=0$, therefore
 $$X_{k+1}={\Pi_{k+2}}_\star X_{k+2}.$$
\item{(2)} The vector field $X_{k+1}$ is tangent to $C_k$, therefore $X_{k+2}$ is tangent to (the strict transform of) $C_k$. 
\end{description}
If $O_{k+2}=C_{k+2}\cap C_k$, we have
$$X_{k+2}(O_{k+2})={\Pi_{k+3}}_\star X_{k+3}(O_{k+2})=0$$
and $\a_{k+2}=0$; by induction we prove $\a_{k+1}=\a_{k+2}=\cdots =0$ till the moment $O_{k+l}$ is not the point $C_{k+l}\cap C_k$ but the point $C_{k+l}\cap C_{k+l-1}$. However if it happens it means that we are in the first case.
 \end{description} 
 
 We have obtained
 \begin{Th} The space of non logarithmic infinitesimal deformations $H^1(U,\T_{\mid U})$ is generated by the $2b_2(S)-l$ cocycles $\mu^i$, $i=0,\ldots,n-1$ and $\t^i$ for indices $i$ such that $O_i$ is at the intersection of two curves. 
 \end{Th}
 The sequence of blowing-ups splits into subsequences 
$$\Bigl(\Pi_{n_1+\cdots+n_\k}\cdots\Pi_{n_1+\cdots+n_\k+l_{\k}-1}\Bigr)\circ\Bigl(\Pi_{n_1+\cdots+n_\k+l_{\k}}\cdots\Pi_{n_1+\cdots+n_\k+n_{\k+1}-1}\Bigr), $$
where $\k=0,\ldots,\rho-1$. The indices wich  correspond to points $O_i$ at the intersection of two curves are 
$$i= n_1+\cdots+n_\k+l_{k},\ldots, n_1+\cdots+n_\k+n_{\k+1}-1,$$
therefore for $\k=0,\ldots,\rho-1$,
$$\a_{n_1+\cdots+n_\k+l_{\k}}=\cdots=\a_{n_1+\cdots+n_\k+n_{\k+1}-1}=0.$$
The equations $(E1)$ become

$$\left\{\begin{array}{lcl}
X_0-{\Pi_1}_\star X_1& = &\displaystyle \a_0\frac{\partial}{\partial u_0}\hspace{10mm}{\rm on} \quad  W_0\\
\hspace{10mm}\vdots&\vdots&\\
X_{l_0-1}-{\Pi_{l_0}}_\star X_{l_0}&=&\displaystyle \a_{l_0-1}\frac{\partial}{\partial u_{l_0-1}} \quad  {\rm on}\quad W_{l_0-1}\\
&&\\
X_{l_0}-{\Pi_{l_0+1}}_\star X_{l_0+1}&=&0 \hspace{10mm} {\rm on} \quad W_{l_0}\\
\hspace{10mm}\vdots&\vdots&\\
X_{n_1-1}-{\Pi_{n_1}}_\star X_{n_1}&=&0 \hspace{10mm} {\rm on}\quad W_{n_1-1}\\
\hspace{10mm}\vdots&\vdots&\\
X_{n_1+\cdots+n_\k}-{\Pi_{n_1+\cdots+n_\k+1}}_\star X_{n_1+\cdots+n_\k+1}& = &\displaystyle \a_{n_1+\cdots+n_\k}\frac{\partial}{\partial u_{n_1+\cdots+n_\k}}\\
&&  \hspace{12mm}  {\rm on} \quad  W_{n_1+\cdots+n_\k}\\
\hspace{10mm}\vdots&\vdots&\\
X_{n_1+\cdots+n_\k+l_{\k}-1}&& \\ 
-{\Pi_{n_1+\cdots+n_\k+l_{\k}}}_\star X_{n_1+\cdots+n_\k+l_{\k}} & = &\displaystyle \a_{n_1+\cdots+n_\k+l_{\k}-1}\frac{\partial}{\partial u_{n_1+\cdots+n_\k+l_{\k}-1}} \\
&& {\rm on} \quad  W_{n_1+\cdots+n_\k+l_{\k}-1}\\
&&\\
X_{n_1+\cdots+n_\k+l_{\k}}&&\\
-{\Pi_{n_1+\cdots+n_\k+l_{\k}+1}}_\star X_{n_1+\cdots+n_\k+l_{\k}+1}&=&0 \hspace{10mm}  {\rm on} \quad W_{n_1+\cdots +n_\k+l_{\k}}\\
&&\\
\hspace{10mm}\vdots&\vdots&\\
X_{n_1+\cdots+n_{\k+1}-1}-{\Pi_{n_1+\cdots+n_{\k+1}}}_\star X_{n_1+\cdots+n_{\k+1}}&=&0 \hspace{10mm}  {\rm on}\quad W_{n_1+\cdots+n_{\k+1}-1}\\
\hspace{10mm}\vdots&\vdots&\\
X_{n_1+\cdots+n_{\rho -1}}-{\Pi_{n_1+\cdots+n_{\rho -1}+1}}_\star X_{n_1+\cdots+n_{\rho -1}+1}& = &\displaystyle \a_{n_1+\cdots+n_{\rho-1}}\frac{\partial}{\partial u_{n_1+\cdots+n_{\rho-1}}}\\
&& \hspace{12mm}  {\rm on} \quad  W_{n_1+\cdots+n_{\rho -1}}\\
\hspace{10mm}\vdots&\vdots&\\
X_{n_1+\cdots+n_{\rho-1}+l_{\rho-1}-1}&& \\ 
-{\Pi_{n_1+\cdots+n_{\rho-1}+l_{\rho-1}}}_\star X_{n_1+\cdots+n_{\rho-1}+l_{\rho-1}} & = &\displaystyle \a_{n_1+\cdots+n_{\rho-1}+l_{\rho-1}-1}\frac{\partial}{\partial u_{n_1+\cdots+n_{\rho-1}+l_{\rho-1}-1}} \\
&& {\rm on} \quad  W_{n_1+\cdots+n_{\rho-1}+l_{\rho-1}-1}\\
X_{n_1+\cdots+n_{\rho-1}+l_{\rho-1}}&&\\
-{\Pi_{n_1+\cdots+n_{\rho-1}+l_{\rho-1}+1}}_\star X_{n_1+\cdots+n_{\rho-1}+l_{\rho-1}+1}&=&0 \hspace{10mm}  {\rm on} \quad W_{n_1+\cdots+n_{\rho-1}+l_{\rho-1}}\\
\hspace{10mm}\vdots&\vdots&\\
X_{n_1+\cdots+n_{\rho}-2}-{\Pi_{n_1+\cdots+n_{\rho}-1}}_\star X_{n_1+\cdots+n_{\rho}-1}&=&0 \hspace{10mm}  {\rm on}\quad W_{n_1+\cdots+n_{\rho}-2}\\
&&\\
X_{n-1}- (\s\Pi_0)_\star X_0 & =& 0 \hspace{10mm} {\rm on}\quad W_{n-1}\\
 \end{array}\right. \leqno{(E2)}$$
 It should be noticed that a block may be reduced to one line, if $l_\k=n_{\k+1}-1$, i.e. if there is in the block only one blowing-up at the intersection of two curves.\\

For $\k=0,\ldots,\rho-1$, the vector fields $X_{n_1+\cdots+n_{\k}+l_{\k}}, \ldots, X_{n_1+\cdots+n_{\k+1}-1}$
 glue together into a vector field that we shall still denote $X_{n_1+\cdots+n_\k+l_{\k}}$. Hence setting
$$\begin{array}{ll}
\P'_0=\P_0\cdots\P_{l_0-1},&  \P''_0=\P_{l_0}\cdots\P_{n_1-1}\\
\hspace{20mm}\vdots&\hspace{20mm}\vdots\\
\P'_{\k}=\P_{n_1+\cdots+n_{\k}}\cdots\P_{n_1+\cdots+n_{\k}+l_{\k}-1},& \P''_{\k}= \P_{n_1+\cdots+n_{\k}+l_{\k}}\cdots\P_{n_1+\cdots+n_{\k+1}-1}\\
\hspace{20mm}\vdots&\hspace{20mm}\vdots\\
  \P'_{\rho-1}=\P_{n_1+\cdots+n_{\rho-1}}\cdots\P_{n_1+\cdots+n_{\rho-1}+l_{\rho-1}-1},&  \P''_{\rho-1}=\P_{n_1+\cdots+n_{\rho-1}+l_{\rho-1}}\cdots\P_{n_1+\cdots+n_{\rho}-1}
  \end{array}$$
  $$  \Pi=\Pi'_0\Pi''_0\cdots\Pi'_{\k}\Pi''_{\k}\cdots  \Pi'_{\rho-1}\Pi''_{\rho-1}.$$
we reduce the system to

$$\left\{\begin{array}{lcl}
X_0-{\Pi_1}_\star X_1& = &\displaystyle \a_0\frac{\partial}{\partial u_0}\hspace{10mm}{\rm on} \quad  W_0\\
\hspace{10mm}\vdots&\vdots&\\
X_{l_0-2}-{\Pi_{l_0-1}}_\star X_{l_0-1}&=&\displaystyle \a_{l_0-2}\frac{\partial}{\partial u_{l_0-2}}  \quad  {\rm on}\quad W_{l_0-2}\\
&&\\
X_{l_0-1}-{\Pi''_0}_\star{\Pi_{n_1}}_\star X_{n_1}&=&\displaystyle \a_{l_0-1}\frac{\partial}{\partial u_{l_0-1}}  \\
&&\hspace{15mm}  {\rm on}\quad W_{l_0-1}\cup\cdots\cup W_{n_1-1}\\
\hspace{10mm}\vdots&\vdots&\\
X_{n_1+\cdots+n_\k}-{\Pi_{n_1+\cdots+n_\k+1}}_\star X_{n_1+\cdots+n_\k+1}& = &\displaystyle \a_{n_1+\cdots+n_\k}\frac{\partial}{\partial u_{n_1+\cdots+n_\k}} \hspace{5mm}  {\rm on} \quad  W_{n_1+\cdots+n_\k}\\
\hspace{10mm}\vdots&\vdots&\\
X_{n_1+\cdots+n_\k+l_{\k}-2}&& \\ 
-{\Pi_{n_1+\cdots+n_{\k}+l_{\k}-1}}_\star X_{n_1+\cdots+n_{\k}+l_{\k}-1} & = &\displaystyle \a_{n_1+\cdots+n_{\k}+l_{\k}-2}\frac{\partial}{\partial u_{n_1+\cdots+n_{\k}+l_{\k}-2}} \\
&& {\rm on} \quad  W_{n_1+\cdots+n_{\k}+l_{\k}-2}\\
&&\\

X_{n_1+\cdots+n_{\k}+l_{\k}-1}&& \\ 
-{\Pi''_{\k}}_\star{\Pi_{n_1+\cdots+n_{\k+1}}}_\star X_{n_1+\cdots+n_{\k+1}} & = &\displaystyle \a_{n_1+\cdots+n_{\k}+l_{\k}-1}\frac{\partial}{\partial u_{n_1+\cdots+n_{\k}+l_{\k}-1}-1} \\
\hspace{64mm}{\rm on}&& W_{n_1+\cdots+n_{\k}+l_{\k}-1} \cup \cdots \cup  W_{n_1+\cdots+n_{\k+1}-1} \\
\hspace{10mm}\vdots&\vdots&\\
X_{n_1+\cdots+n_{\rho -1}}-{\Pi_{n_1+\cdots+n_{\rho -1}+1}}_\star X_{n_1+\cdots+n_{\rho -1}+1}& = &\displaystyle \a_{n_1+\cdots+n_{\rho-1}}\frac{\partial}{\partial u_{n_1+\cdots+n_{\rho-1}}}\\
&& \hspace{12mm}  {\rm on} \quad  W_{n_1+\cdots+n_{\rho -1}}\\
\hspace{10mm}\vdots&\vdots&\\
X_{n_1+\cdots+n_{\rho-1}+l_{\rho-1}-2}&& \\ 
-{\Pi_{n_1+\cdots+n_{\rho-1}+l_{\rho-1}-1}}_\star X_{n_1+\cdots+n_{\rho-1}+l_{\rho-1}-1} & = &\displaystyle \a_{n_1+\cdots+n_{\rho-1}+l_{\rho-1}-2}\frac{\partial}{\partial u_{n_1+\cdots+n_{\rho-1}+l_{\rho-1}-2}} \\
&& {\rm on} \quad  W_{n_1+\cdots+n_{\rho-1}+l_{\rho-1}-2}\\

X_{n_1+\cdots+n_{\rho-1}+l_{\rho-1}-1}-({\Pi''_{\rho-1}}{\s\Pi_0})_\star X_0 & = &\displaystyle \a_{n_1+\cdots+n_{\rho-1}+l_{\rho-1}-1}\frac{\partial}{\partial u_{n_1+\cdots+n_{\rho-1}+l_{\rho-1}-1}} \\

&& {\rm on} \quad  W_{n_1+\cdots+n_{\rho-1}+l_{\rho-1}-1}\cup \cdots \cup W_{n-1}\\
 \end{array}\right. \leqno{(E3)}$$
When $\rho=1$, i.e. when there is only one tree, the linear system reduces to

$$\left\{\begin{array}{lccl}
X_0-{\Pi_1}_\star X_1& = &\displaystyle \a_0\frac{\partial}{\partial u_0} & {\rm over} \;  W_0\\
\hspace{10mm}\vdots&\vdots&\\
X_{l-2}-{\Pi_{l-1}}_\star X_{l-1}&=&\displaystyle \a_{l-2}\frac{\partial}{\partial u_{l-2}} & {\rm over} \;  W_{l-2}\\
&&&\\
X_{l-1}-({\Pi''\circ\s\circ\Pi_{0}})_\star X_{0}&=&\displaystyle \a_{l-1}\frac{\partial}{\partial u_{l-1}}  & {\rm over} \;  W_{l-1}\cup\cdots\cup W_{n-1}
\end{array}\right. \leqno{(E4)}$$
\begin{Cor} \label{relationcocycles} A relation among the cocycles $[\t^i]$ and $[\mu^i]$, $i=0,\ldots,n-1$ contains only $[\t^i]$ in $H^1(S,\T(-Log\  D))$, i.e. indices for which the blown up point $O_i$ is generic.
\end{Cor}

\begin{Cor}[\cite{N1}]\label{baseIH} Let $S$ be a Inoue-Hirzebruch surface with Betti number $b_2(S)=n\ge 1$, then the cocycles $\t^ i$ and $\mu^ i$, $i=0,\ldots,n-1$ define the versal deformation and the versal logarithmic deformation is trivial. Moreover an Inoue-Hirzebruch  surface $S=S_0$ with two cycles of rational curves $\G$ and $\G'$ can be deformed into a Hopf surface with two elliptic curves $\G_u$ and $\G'_u$ blown up respectively $-\G^2$ and $-\G'^2$ times.
\end{Cor}
Proof: In the explicit construction of Inoue-Hirzebruch surfaces \cite{D2}, there is no generic blown up points and $h^1(S,\T)=2n$, hence we have an explicit base of $H^1(S,\T)$ and explicit universal deformation. It is easy to see that any singular point of a cycle may be smoothed for even as well for odd Inoue-Hirzebruch surface.\hfill $\Box$\\

For moduli space of Oeljeklaus-Toma see \cite{OT}.
\begin{Cor}\label{ouvertureOT} Fix any $J$, $M$, $\s$ and consider a large family  $\widehat{\Phi}_{J,\s}:\widehat{\cal S}_{J,\s}\to \widehat{B}_{J}$, then the family is generically versal and the image of the stratum $B_{J,M}$ in the Oeljeklaus-Toma moduli space of surfaces with intersection matrix $M$ contains an open set.
\end{Cor}
Proof: Any large family degenerates to Inoue-Hirzebruch surfaces and at the point $O_{IH}\in B_{J}$ corresponding to this Inoue-Hirzebruch surface, the family is versal. The point $O_{IH}$ is in the closure of any stratum. By openess of the versality it is versal in a neighbourhood  hence on an open set of any stratum. Since the family is generically versal, the image in Oeljeklaus-Toma coarse moduli space contains an open set.\hfill$\Box$ 

\begin{Def} Let $X$ be a complex manifold of dimension $m$. We shall say that $X$ admits a birational structure if there is an atlas $\frak U=(U_i)_{i\in I}$ with charts $\f_i:U_i\to \bb C^m$ such that for each pair $\{i,j\}$ such that $U_i\cap U_j\neq\emptyset$, $\f_j\circ\f_i^{-1}:\f_i(U_{ij})\to \f_j(U_{ij})$ is the restriction of a birational mapping of $\bb P^m(\bb C)$.
\end{Def}
\begin{Cor} Let $M$ be any intersection matrix of a minimal compact complex surface containing a GSS (i.e. of a Kato surface) then the O-T moduli space of such surfaces contains a non empty open set of surfaces admitting a birational structure.
\end{Cor}
Proof: We take $\s=Id$, then the gluing map $\s\circ\Pi$ is birational, then apply Corollary 3.\ref{ouvertureOT}.\hfill$\Box$

\subsection{Infinitesimal logarithmic deformations}
We have seen that a relation is only possible among infinitesimal logarithmic deformation. In fact it can contain neither $\t^i$  when the curve $C_i$ meets two other curves.  In order  to avoid an overflow of notations, we give a complete proof for surfaces with only one tree and we postpone it to the appendix \ref{thehardpart}. The idea of the computation is to work in the first infinitesimal neighbourhood of the maximal divisor. Vanishing of other coefficients should imply to work (if possible) in the successive infinitesimal neighbourhoods.
\begin{Prop}\label{relationlogarithmique} Let $(S,C_0)$ be a marked surface.  If $\sum_{i=0}^{n-1}\a_i[\t^i]=0$ is a relation, then $\a_k=0$ for any index $k$ such that one of the two conditions is fulfilled
\begin{itemize}
\item $O_k$ is the intersection of two rational curves,
\item $O_k$ is a generic point but $C_k$ meets two other curves.
\end{itemize}
 In particular, if the unique regular sequences $r_m$ are reduced to one curve (i.e. $m=1$),  $\{[\t^i],[\mu^i] \mid 0\le i\le n-1\}$ (resp. $\{[\t^i] \mid O_i\ {\rm is\ a\ generic\ point}\}$) is an independant family of $H^1(S,\T)$ (resp. of $H^1(S,\T(-Log D))$) and a base if there is no non trivial global vector fields. 
\end{Prop}
\begin{Rem} {\rm By induction it is possible to show that for any $k<r+s-(p+q)$, a similar Cramer system may by defined and that $\a_k=0$. However, it is not possible to achieve the proof in this way because when $k=r+s-(p+q)$ a new unknown appears. This difficulty is explained by the fact that in general} there is a relation  among the $[\t^i]$'s or  a class vanishes.
\end{Rem}

\section{The hypersurface of non versality}
The constructed families are generically versal. We show in this section that the locus of non-versality is non empty hypersurface. When $b_2(S)=2$ it is the ramification set of the canonical mapping from each stratum to Oeljeklaus-Toma moduli space. It is conjectured that it is a general phenomenon.
\subsection{The generically logarithmically versal family $\cal S_{J,M,\s}$} 

Notice that by \cite{OT}\S 6, $k,t$ and the integers $(m_1,\ldots,m_\rho)$ determine completely the sequence of self-intersections of the rationals curves, i.e. the invariant $\s_n(S)$ and the intersection matrix $M=M(S)$. We have two families of logarithmic deformations, the first one $$\cal S_{J,M,\s}\to B_{J,M}=(\bb C^\star\times \bb C^{l_0-1})\times\cdots\times (\bb C^\star\times \bb C^{l_\kappa -1})\times\cdots\times (\bb C^\star\times\bb C^{l_{\rho-1}-1}),$$
 is generically versal by (3.\ref{logversaldeformation}), the second one 
$$\cal S_{k,\s ,m_1,\ldots,m_\rho}\to U_{k,\s ,m_1,\ldots,m_\rho}=\bb C^\star_{\l}\times (\bb C^{\star})^{\rho-1}\times \bb C^{\e(k,\s ,m_1,\ldots,m_\rho)}$$
  is versal at every point \cite{OT} thm 7.13, therefore
$$\dim B_{J,M}=\dim (\bb C^{\star})^\rho\times \bb C^{\e(k,\s ,m_1,\ldots,m_\rho)},$$
 $$l=\rho+\e(k,\s ,m_1,\ldots,m_\rho),$$
and the bases are equal up to permutation of the factors.\\
\begin{Lem}\label{degrefibreplat} Let $(g_a)_{a\in B_J}$ be a differentiable family of Gauduchon metrics on $\Phi_{J,\s}:\cal S_{J,\s}\to B_J$, $\o_a$ be its associated $(1,1)$ form 
 and let 
$$\deg_{g_a}:H^1(S_a,\cal O^\star)\to \bb R, \quad \deg_{g_a}(L)=\int_{S_a}c_1(L)\wedge \o_a$$
be the degree of a line bundle. Then there is a non vanishing differentiable negative function
$$C:B_J\to \bb R_-^\star$$
such that for any $L^\l\in H^1(S_a,\bb C^\star)\simeq \bb C^\star$,
$$\deg_{g_a}(L^\l)=C(a)\log|\l|.$$
\end{Lem}
Proof:  For any $a\in B_J$ the Lie group morphism $\deg_{g_a}:H^1(S_a,\bb C^\star)\simeq \bb C^\star\to \bb R$ has the form $\deg_{g_a}(L^\l)=C\log|\l|$ where $C\neq 0$ since this morphism is surjective. Besides the family of Gauduchon metric depends differentiably on $a\in B_J$,therefore $C:B_J\to \bb R$ is always positive or always negative. Now, on Enoki surfaces $S_a$, denote by $\G_a$ the topologically trivial cycle of rational curves. Then, by Gauduchon theorem \cite{G}, \cite{LT}, and \cite{DO1}
$$vol(\G_a)=\deg_{g_a}([\G_a])=\deg_{g_a}(L^{t(a)})=C(a)\log|t(a)|$$
where $t(a)=\tr(S_a)$ is the trace of the surface satisfies $0<|t(a)|<1$, therefore $C(a)<0$. Since $C(a)<0$ when $S_a$ is a Enoki surface, $C(a)<0$ everywhere.\hfill$\Box$\\

Remark that a numerically $\bb Q$-anticanonical divisor $D_{-K}$ on a surface $S$ is a solution of a linear system whose matrix is the intersection matrix $M=M(S)$ of $S$. Therefore the index is the least integer $m$ such that $mD_{-K}$ is a divisor and this integer is fixed on any logarithmic family $\Phi_{J,M,\s}:\cal S_{J,M,\s}\to B_{J,M}$.

  \begin{Lem}\label{kappaFavre} Let $F(z_1,z_2)=(\l z_1z_2^\frak s  +P(z_2)+cz_2^{\frac{\frak s  k}{k-1}},z_2^k)$ be a Favre contracting germ associated to a surface of intermediate type $S$. Let $\mu=index(S)\in\bb N^\star$ and $\k\in\bb C^\star$ such that $H^0(S, K_S^{\ot -\mu}\ot L^{\k})\neq 0$. Then 
$$\k=k(S)^{-\mu}\l^{-\mu}.$$
\end{Lem}
Proof: A global section  $\t\in H^0(S,K^{-\mu}\ot L^\k)$ induces a germ  $\t=z_2^\a A(z)\left(\frac{\part}{\part z_1}\wedge\frac{\part}{\part z_2}\right)^{\ot \mu}$  which satisfies  the condition
$$\t(F(z))=\k \Bigl(\det DF(z)\Bigr)^\mu \t(z),$$
where $\a$ is the vanishing order  of $\t$ along $C_{n-1}$ and $A(0)\neq 0$.   Since $\det DF(z)=\l k z_2^{\frak s  +k-1}$, comparison  of lower degree terms gives
$$z^{k\a}A(0)=\k(\l k)^\mu z_2^{\mu(\frak s  +k-1)+\a}A(0)$$ 
hence
$$\left\{\begin{array}{ll}\a(k-1)=\mu(k-1+\frak s  )\\
\\
\k=k(S)^{-\mu}\l^{-\mu}.
\end{array}\right.$$
\hfill$\Box$\\
\begin{Lem} Let $S$ be a minimal complex surface, $\mu$ the index of $S$ and $\k$ such that $H^0(S,K^{\ot -\mu}\ot L^\k)\neq 0$. Then a section of $K^{\ot -\mu}\ot L^\k$ vanishes on all the rational curves in $S$.
\end{Lem}
Proof: Let $D_i$, $i=0,\ldots,n-1$ be the $n$ rational curves in $S$ and suppose that 
$$K^{\ot -\mu}\ot L^\k=\sum_{i=0}^{n-1}k_iD_i.$$
 We have $k_i\ge 0$ for all $i=0,\ldots,n-1$; if one coefficient vanishes,  say  $k_0=0$,  on one hand, since the maximal divisor is connected,
$$c_1(K^{\ot -\mu}\ot L^\k).D_0=\sum_{i=1}^{n-1}k_iD_iD_0>0$$
and on second hand, by  adjunction formula
$$c_1(K^{\ot -\mu}\ot L^\k).D_0=-\mu c_1(K).D_0=\mu(D_0^2+2)\le 0$$
we obtain a contradiction.\hfill$\Box$

\begin{Prop}\label{kappalambda}  Let $\Phi_{J,M,\s}:\cal S_{J,M,\s}\to B_{J,M}$ be a logarithmic family of marked intermediate surfaces with $J\neq\emptyset$ and $\mu$ the common index of the surfaces. Then \\
1) there exists a unique surjective holomorphic function 
 $$\begin{array}{cccc}\k=\k_{J,M,\s}:&B_{J,M}&\to&\bb C^\star\\
 &a&\mapsto&\k(a)
 \end{array}$$
 such that  $H^0(S_a,K_{S_a}^{-\mu}\ot L^{\k(a)})\neq 0$.\\
  2) If  the surfaces admit twisted vector fields  there exists a unique surjective holomorphic function 
 $$\begin{array}{cccc}\l=\l_{J,M,\s}:&B_{J,M}&\to&\bb C^\star\\
 &a&\mapsto&\l(a)
 \end{array}$$
 such that the marked surface $(S_a,C_{0,a})$ is defined by a germ of the form $$F_a(z_1,z_2)=(\l(a) z_1z_2^\frak s  +P_a(z_2),z_2^k).$$
3) The fibers $K_\a:=\{\k=\a\}$ (resp. $\L_\a:=\{\l=\a\}$), $\a\in\bb C^\star$, are closed in $B_J$ hence analytic in $B_{J}\supset B_{J,M}$.\\
4) Let $\overline{B_{J,M}}\subset B_J$ be the closure of $B_{J,M}$ in $B_J$, i.e. the union of $B_{J,M}$ with the smaller strata, then $\k_{J,M,\s}$ extends holomorphically to $\k_{J,M,\s}:\overline{B_{J,M}}\to \bb C$ and $\k_{J,M,\s}^{-1}(0)=\overline{B_{J,M}}\setminus B_{J,M}$.
\end{Prop}
Proof:  1) For $a\in  B_{J,M}$, the complex number $\k(a)$ satisfies $h^0(S_a,K_{S_a}^{-\mu}\ot L^{\k(a)})=1$. It is unique because there is no topologically trivial divisor.\\
We consider a new base space $B_{J,M}\times \bb C^\star$, and let
$$pr_1: B_{J,M}\times \bb C^\star\to B_{J,M},$$
be the  first projection. Let $\cal K\to \cal S_{J,M,\s}$ be the relative canonical line bundle and $\cal L\to \cal S_{J,M,\s}\times\bb C^\star$ be the tautological line bundle such that $\cal L_{a,\tau}$ is the line bundle $L^\tau$ over $S_a$. We consider the family of rank one vector bundles
  $$pr_1^\star\cal K\ot  \cal L \to pr_1^\star\cal S_{J,M,\s}\stackrel{pr_1^\star\Phi_{J,M,\s}}{\longrightarrow} B_{J,M}\times\bb C^\star.$$
Then $\left(pr_1^\star \cal K\ot  \cal L \right)_{(a,\a)}= K_a\ot L^\a$. The set of points 
$$Z=\{(a,\a)\in B_{J,M}\times \bb C^\star \mid h^0(S_a,K_a\ot L^\a)>0\}$$
is an analytic subset. Let 
$$pr:Z\to B_{J,M}$$
 be the restriction to $Z$ of the first projection $pr_1$ over $B_{J,M}$. Then $pr$ is surjective by hypothesis. Each fiber contains only one point. Moreover $pr$ is proper: in fact we consider the closure $\overline{Z}\subset B_{J,M}\times\bb P^1(\bb C)$. By Remmert-Stein theorem, either $\overline{Z}$ is an analytic set in $B_{J,M}\times\bb P^1(\bb C)$ or contains at least one of the hypersurfaces $B_{J,M}\times \{0\}$ or  $B_{J,M}\times\{\infty\}$. But it is impossible because each fiber contains only one point. Therefore $\overline{Z}$ is analytic and $\overline{pr}:\overline{Z}\to B_{J,M}$ is proper hence a ramified covering. Since there is only one sheet, it is the graph of a holomorphic mapping $\k:B_{J,M}\to \bb P^1(\bb C)$. Since for every $a\in B_{J,M}$, $pr^{-1}(a)$ contains exactly one point in $\bb C^\star$, $\k$ has only values in $\bb C^\star$.\\
 Now, $\k$ cannot be constant because $\k=(k\l)^{-\mu}$ and $\l$ is a parameter of a logarithmic versal family, therefore the non-constant mapping $\k:(\bb C^\star)^\rho\times \bb C^{l-\rho}\to \bb C^\star$ is surjective.\\
2) By lemma 4.\ref{kappaFavre}, and \cite{DO1} lemma 4.2.3) there is a numerically canonical divisor if and only if there is a numerically tangent divisor therfore $\k=k^{-1}\l^{-1}$ (see the remark below).\\
3) Consider the hypersurface $\{\k=\a\}\subset B_{J,M}$ for $\a\in\bb C^\star$. The closure $\overline{B_{J,M}}$ of $B_{J,M}$ in $B_J$ is the union of $B_{J,M}$ with lower strata, hence  $\overline{B_{J,M}}\setminus B_{J,M}$ is also a hypersurface. Remmert-Stein theorem shows that $\overline{\{\k=\a\}}$ is analytic or contains an irreducible component of $\overline{B_{J,M}}\setminus B_{J,M}$. However the second possibility is excluded by Grauert semi-continuity theorem because on a whole stratum we would have $H^0(S_a,K^{-\mu}\ot L^\a)\neq 0$ which is impossible because the twisting parameter is not constant. Therefore the slice has an extension. If $\overline{\{\k=\a\}}\cap (\overline{B_{J,M}}\setminus B_{J,M})\neq\emptyset$, the line bundle $K^{-\mu}\ot L^\a$ has a section over ${\cal S_{J,\s}}_{\mid \overline{\{\k=\a\}}}$ hence the zero locus which is the union of all the rational curves by \cite{DO1} would be is a flat family of divisors; however it is impossible because  the configuration changes contradicting flatness (it can be seen that the curve whose self-intersection decreases has a volume which tends to infinity (see \cite{DT1})). Therefore $\overline{\{\k=\a\}}\cap (B_{J}\setminus B_{J,M})=\emptyset$ and each slice is already closed in $B_{J,M}$.\\
4) Since the fibers $K_\a$ are closed in $B_J$, 
$$\lim_{a\to \overline{B_{J,M}}\setminus B_{J,M}}\k_{J,M,\s}(a)=0\ {\rm or}\ \infty.$$
Let $\cal K$ be the relative canonical line bundle, $\t\in  H^0(\cal S_{J,M,\s},\cal K^{-\mu}\ot L^{\k})$ be the flat family of sections over $B_{J,M}$ and $Z$ the associated divisor of zeroes of $\t$. By lemma \thesection.\ref{degrefibreplat},
$$vol(Z_a)=\deg_{g_a}([Z_a])=\deg_{g_a}(K_a^{-\mu}\ot L^{\k(a)})=-\mu\deg_{g_a}(K_a)+C(a)\log|\k(a)|$$
Since $a\mapsto \deg_{g_a}(K_a)$ is differentiable, hence bounded, and $vol(Z_a)>0$, the limit of $\k=\k_{J,M,\s}(a)$ near $\overline{B_{J,M}}\setminus B_{J,M}$ cannot be $\infty$, therefore $\k$ extends continuously and holomorphically.\\
\hfill $\Box$

\begin{Rem} 1) If $index(S)=1$, we have $\l^{-1}=k(S)\k$, i.e. the invariant used here is the inverse of the invariant $\l=\l(S)$ in \cite{DO1}.\\
2) If  $index(S)\neq 1$, $\l(a)$ is defined up to a $(k-1)$-root of unity.
\end{Rem}

\subsection{Existence of relations among the tangent cocycles}
In this section we show that the classes $\{[\t_i]\mid O_i\  {\rm is\ generic}\}$ cannot be linearly independant everywhere, there exist an obstruction.

\begin{Lem}\label{codimaumoinsdeux} Let $\Phi_{J,\s}:\cal S_{J,\s}\to B_J$ any family of minimal surfaces and $B_{J,M}$ any stratum of $B_J$ of intermediate surfaces. Let
$$Z=\{u\in B_J\mid h^0(S_u,\T_u)>0\},$$
be the analytic set of parameters $u\in B_J$ of surfaces with non trivial vector fields.
Then any irreducible component $Z_M$ of $Z$ such that $Z_M\cap B_{J,M}\neq\emptyset$ is contained in $B_{J,M}$ in particular is at least 2-codimensional.
\end{Lem}
Proof: 1) In $B_J$ the Zariski open set of parameters $\{t\neq 0\}$ does not parameterize Inoue surfaces by Proposition \ref{minimalcodim2}.\\
2) Let $M$ such that the surfaces of the stratum $B_{J,M}$ admit twisted vector fields. Then  $Z_M\cap B_{J,M}$ is contained in the hypersurface $\{\l_{M}=1\}$ therefore $Z$ is an analytic set of codimension at least one in a hypersurface, therefore is at least 2-codimensional.\hfill$\Box$

\begin{Th}\label{logversaldeformation} Let $(S,C_0)$ be a minimal marked surface containing a GSS of intermediate type,  with $n=b_2(S)$.  Let $J=I_\infty(C_0)$ and let $\Phi_{J,\s}:\cal S_{J,\s}\to B_J$ be the family of surfaces with GSS associated to $J$ and $\s$. Then, there exists a non empty hypersurface $T_{J,\s}\subset B$ containing $Z=\{u\in B\mid h^0(S_u,\T_u)>0\}$ such that for $u\in B_J\setminus T_{J,\s}$,
\begin{description}
\item{a)}  $\{[\t^i_u],[\mu^i_u] \mid 0\le i\le n-1\}$ is a base of $H^1(S_u,\T_u)$,
\item{b)} $\{[\t^i_u] \mid O_i\ {\rm is\ generic}\}$ is a base of   $H^1(S_u,\Theta_u(-Log\ D_u))$.
\end{description}
Moreover 
\begin{description}
\item{i)} If $T_{J,\s}$ intersects a stratum $B_{J,M}$ then $T_{J,\s}\cap B_{J,M}$ is a hypersurface in $B_{J,M}$,
\item{ii)} $T_{J,\s}$  intersect each stratum $B_{J,M}$ such that the corresponding surfaces admit twisted vector fields and $Z\cap B_{J,M}\subset T_{J,\s}$,
\end{description}
\end{Th}
Proof:  At the point $a=0$ (i.e.  for $a_i=0$ for all $i$), $S_a$ is a Inoue-Hirzebruch surface. By Corollary 3.\ref{baseIH},  the family $\{[\t^i],[\mu^i] \mid 0\le i\le n-1\}$ is a base of $H^1(S_a,\T_a)$ and the family $(\cal S_{J,\s},\Phi_{J,\s},B_J)$ is versal (even universal) at this point, therefore the set of points $T_{J,\s}$ where the family $\{[\t^i],[\mu^i] \mid 0\le i\le n-1\}$ is not effective is at least one codimensional. By openness of the versality the same property holds in a neighbourhood. Since each stratum $B_{J,M,\s}$ has $a=0$ in its closure, $T_{J,\s}\cap B_{J,M,\s}$ is also at least one codimensional. hypersurface in $B_{J,M,\s}$ and we have $i)$. Outside $Z$, $R^1\Pi_\star\T$ is locally free sheaf of rank $2n$, therefore this family is free outside  a hypersurface $T_{J, \s}$. By lemma 4.\ref{codimaumoinsdeux}, $Z$ is of codimension at least 2, therefore by Teleman \cite{Te4}, an irreducible component of $T_{J,\s}$ of codimension one contains $Z$. At a generic point of $a\in B_J$, $S_a$ is a Enoki surface. If $a_i=0$ for exactly one index $i\in J$, we have $\s_n(S_a)=2n+1$ and for this configuration of curves there exists twisted vector fields, therefore by lemma 4.\ref{codimaumoinsdeux}, $T_{J,\s}$ is not empty. This gives ii).
\hfill$\Box$

\begin{Rem} 1) It is possible to prove that  $T_{J,\s}$ does not intersect those strata near Inoue-Hirzebruch surfaces which have only regular sequences $r_1$.\\
2) We shall see that for $\s=Id$, any surfaces with only one tree, $T_{J,\s}$ is a ramification locus of $B_{J,\s}$ over the Oeljeklaus-Toma moduli space \cite{OT}.
\end{Rem}

\subsection{Surfaces with $b_2=2$.\label{subsectionb2=2}}
\subsubsection{Rational curves}
Up to a circular permutation intersection matrix and configuration of the curves $D_0$ and $D_1$ are the following:
\begin{itemize}
\item Surfaces of trace $t\neq 0$: Enoki surfaces and Inoue surfaces,
$$M(S)=\left(\begin{array}{rr}-2&2\\2&-2\end{array}\right),\quad [D_0]=e_0-e_1,\quad [D_1]=e_1-e_0, \quad [D_0]+[D_1]=0. $$

\begin{center}
\includegraphics[width=3cm]{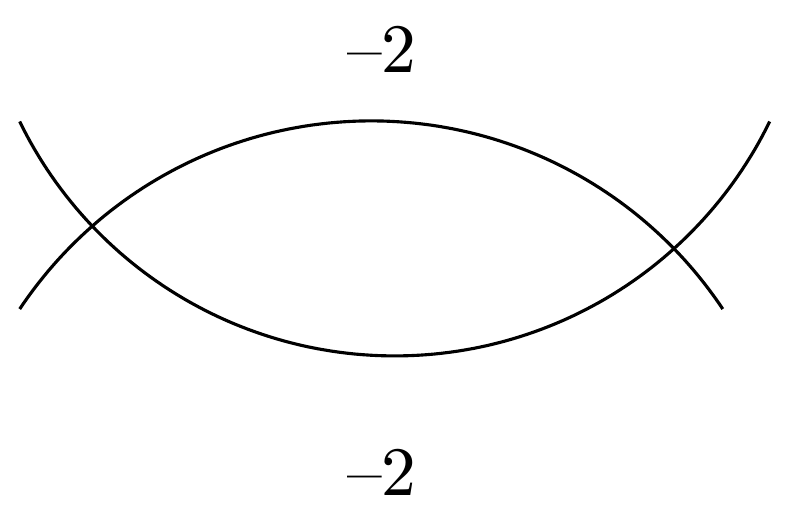}
\end{center}
\item Intermediate surface
$$M(S)=\left(\begin{array}{rr}-1&1\\1&-2\end{array}\right),\quad [D_0]=e_0-e_1-e_0=-e_1, \quad[D_1]=e_1-e_0.$$

\begin{center}
\includegraphics[width=3cm]{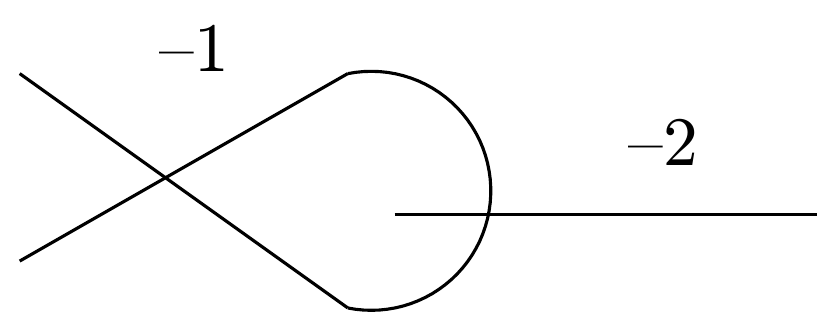}
\end{center}
\item Inoue-Hirzebruch surfaces
$$M(S)=\left(\begin{array}{rr}-4&2\\2&-2\end{array}\right),\quad [D_0]=e_0-e_1-e_0-e_1=-2e_1, \quad [D_1]=e_1-e_0, \quad [D_0]+[D_1]=-e_0,$$
$$M(S)= \left(\begin{array}{rr}-1&0\\0&-1\end{array}\right), \quad [D_0]=e_0-e_1-e_0=-e_1,\quad [D_1]=e_1-e_0-e_1=-e_0.$$
\begin{center}
\includegraphics[width=7cm]{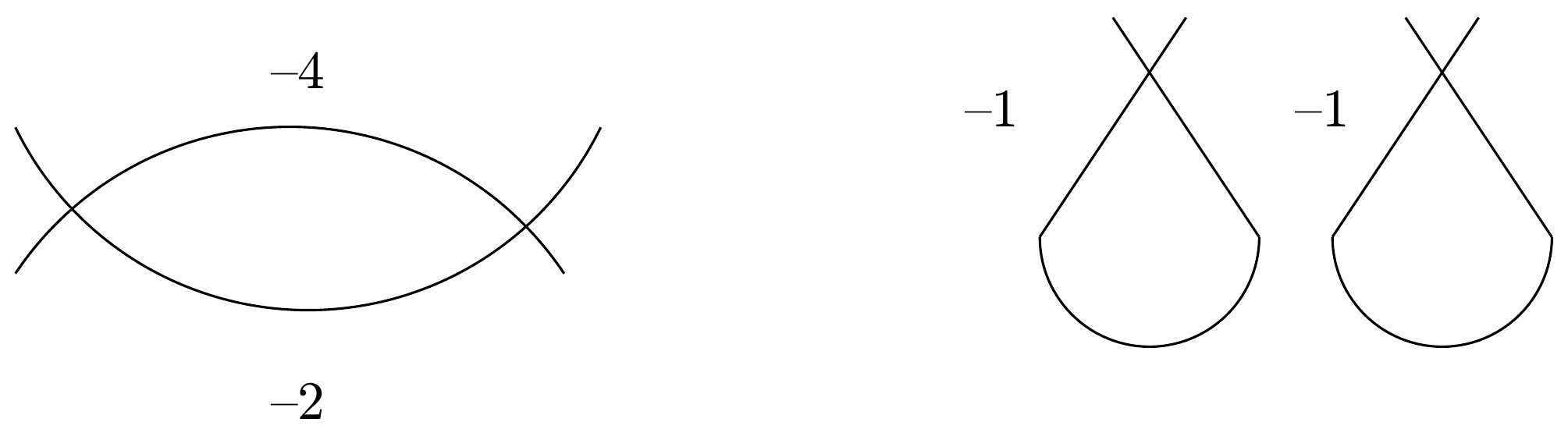}
\end{center}
\end{itemize}

\subsubsection{Intermediate surfaces}
 We consider  intermediate surfaces $S$, since the problem of normal forms is solved for the other cases. There are two curves: one rational curve with a double point $D_1^2=-1$ with one tree $D_0^2=-2$, $D_0D_1=1$. Favre polynomial germs are
$$F_c(z_1,z_2)=(\l z_1z_2+z_2+cz_2^2,z_2^2)$$
where $k=k(S)=2$, $\frak s=1$. Invariant vector fields $\t$ exist if and only if $\l=1$ in which case 
$$\t(z)=\a z_2^{\frak s/(k-1)}\frac{\part}{\part z_1}=\a z_2\frac{\part}{\part z_1}, \quad \a\in\bb C$$
Intermediate surfaces belong to three families, namely for $J=\{0\}$, $J=\{1\}$ and $J=\{0,1\}$.\\
{\bf \boldmath Case $J=\{0\}$}\\
The case $J=\{1\}$ is similar.\\
The family of germs defining surfaces of $\Phi_{J,M,\s}:\cal S_{J,M,\s}\to B_{J,M}$ are 
$$G^J_a(z_1,z_2)=G_a(z_1,z_2)=\bigl(z_2,(z_1+a_1)z_2^2\bigr), \quad a_1\in \bb C^\star, \quad a=(0,a_1)$$
\begin{center}
\includegraphics[width=12cm]{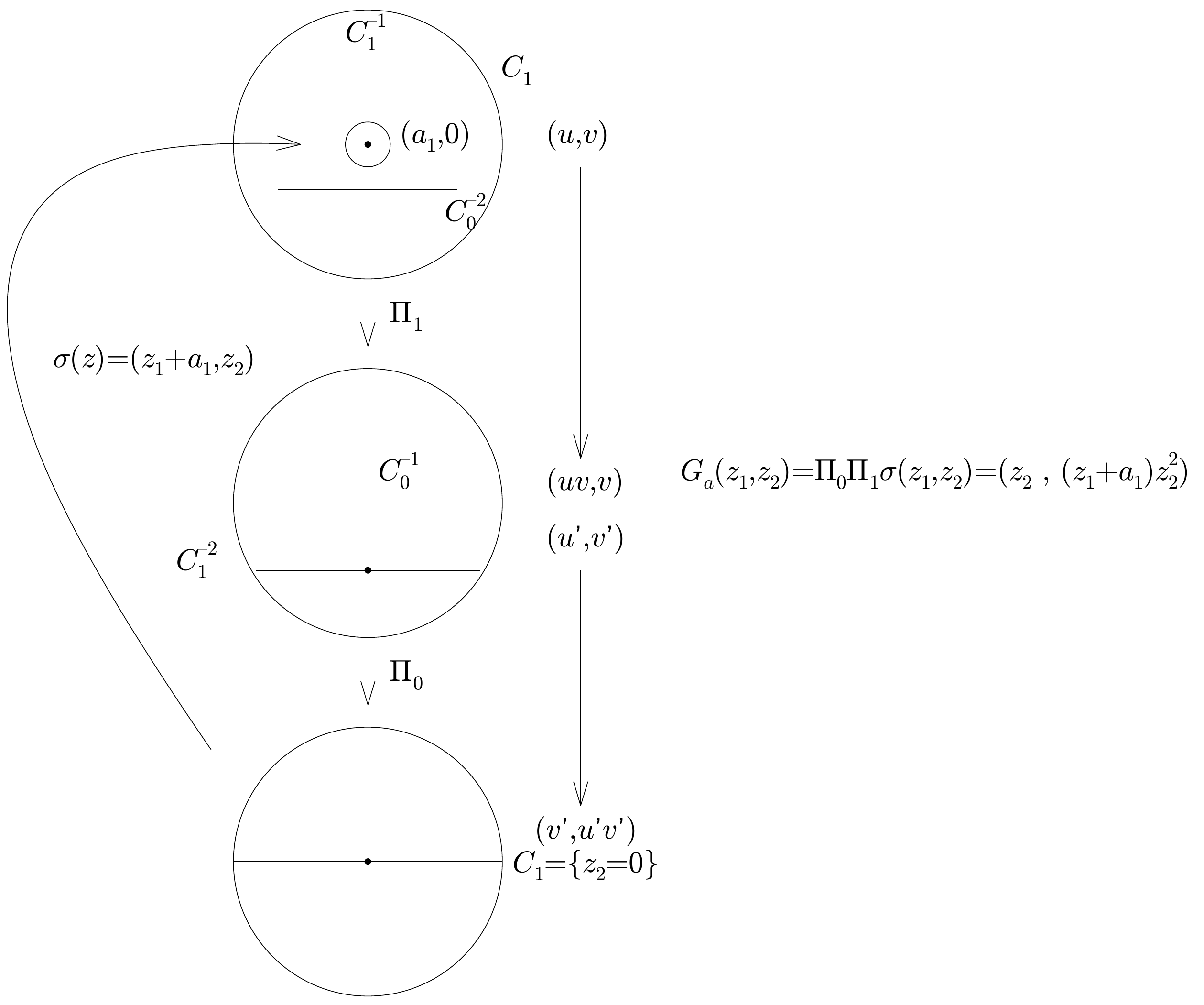}
\end{center}
A germ of isomorphism $\f$ which conjugates $G_a$ and $G_{a'}$ leaves the line $\{z_2=0\}$ invariant, therefore $\f$ has the form $\f(z)=(\f_1(z),Bz_2(1+\t(z))$. A simple computation shows that if $G_a$ and $G_{a'}$ are conjugated then 
$$a_1=\pm a_1'.$$
Besides if we want to determine the twisting parameter $\k$ such that $H^0(S,K^{-1}\otimes L^\k)\neq 0$, we have to solve the equation
$$\mu(G_a(z))=\k \det DG_a(z)\mu(z).$$
Using the relation $D_{-K}=D_\t +D$ of \cite{DO1} or by a direct computation, we know that a section $\mu$ of the twisted anticanonical bundle has to vanish at order two along the cycle, i.e. along $\{z_2=0\}$, therefore $\mu(z)=z_2^2 A(z)\frac{\part}{\part z_1}\wedge \frac{\part}{\part z_2}$ where $A(0)\neq 0$. A straightforward computation shows that
$$\k=-a_1^2.$$
By \cite{DO1} the relation between the twisting parameters $\k$ and  $\wt\l$ chosen so that $H^0(S,\T\otimes L^{\wt\l})\neq 0,$
is $\wt\l=k\k$. Here $k=k(S)=2$, therefore $\wt\l=-2a_1^2$ and 
\begin{itemize} 
\item There is a non-trivial global global vector field if and only if $\wt\l=1$ if and only if 
$$a_1^2=-\frac{1}{2}$$
\item Since $\l=1/\wt\l\in \bb C^\star$ is a parameter of the coarse moduli space, $G_a$ is conjugated to $G_{a'}$ if and only if the corresponding surfaces $S(G_a)$ and $S(G_{a'})$ are isomorphic if and only if $a_1^2=a_1'^2$. In particular the mapping $\bb C^\star =B_{J,M}\to B_{2,1,1}=\bb C^\star$ is 2-sheeted non ramified covering space
\end{itemize}

We are now looking for the missing parameter: we choose $\s(z)=(z_1+\xi z_2+a_1,z_2)$; the infinitesimal deformation is 
$$X(u_1,v_1)=v_1\frac{\part}{\part u_1}.$$
With 
 $$G_{a,\xi}(z_1,z_2)=\bigl(z_2,(z_1+\xi z_2+a_1)z_2^2\bigr)$$
 the same computation gives $\k=-a_1^2$. With $a$ fixed such that $a_1^2=-1/2$ (in order to have a global vector field),   the conjugation relation
 $$\f(G_{a,\xi}(z))=F_c(\f(z))$$
 yields the relations
$$\left\{\begin{array}{cl}
(I)&\f_1\Bigl(z_2,(z_1+\xi z_2+a_1)z_2^2\Bigr)=Bz_2\Bigl[(1+\f_1(z))(1+\t(z))+cBz_2(1+\t(z))^2\Bigr] \\
&\\
(II)&(z_1+\xi z_2+a_1)\Bigl(1+\t\bigl(z_2,(z_1+\xi z_2+a_1)z_2^2\bigr)\Bigr)=B(1+\t(z))^2 
\end{array}\right.$$
we compare the homogeneous parts of the same degree
\begin{itemize}
\item till degree two and homogeneous part $z_1z_2^2$  in $(I)$
\item till degree one and homogeneous part $z_1^2$ in $(II)$
\end{itemize}
A straightforward computation with $a_1^2=-1/2$ yields
 $$c=\xi +2.$$
 Therefore all surfaces with global vector fields are obtained when $\xi$ moves in $\bb C$, and $X$ acts by translation. In particular when $b_2(S)=2$, all surfaces are obtained by simple birational mappings obtained by composition of blowing-ups and an affine map at a suitable place.\\
 We extend the family to Enoki surfaces. The family of marked surfaces $\Phi_{J,\s}:\cal S_{J,\s}\to B_J$ is defined by the family of polynomial germs
$$G_a(z_1,z_2)=\bigl(z_2,z_2^2(z_1+a_1)+a_0z_2\bigr), \quad a=(a_0,a_1)$$
We have
$$\tr(S)=\tr DG_a(0)=a_0$$
therefore $|a_0|<1$. The open set $B_J=\D_{a_0}\times \bb C_{a_1}$ has the following strata
\begin{center}
\includegraphics[width=12cm]{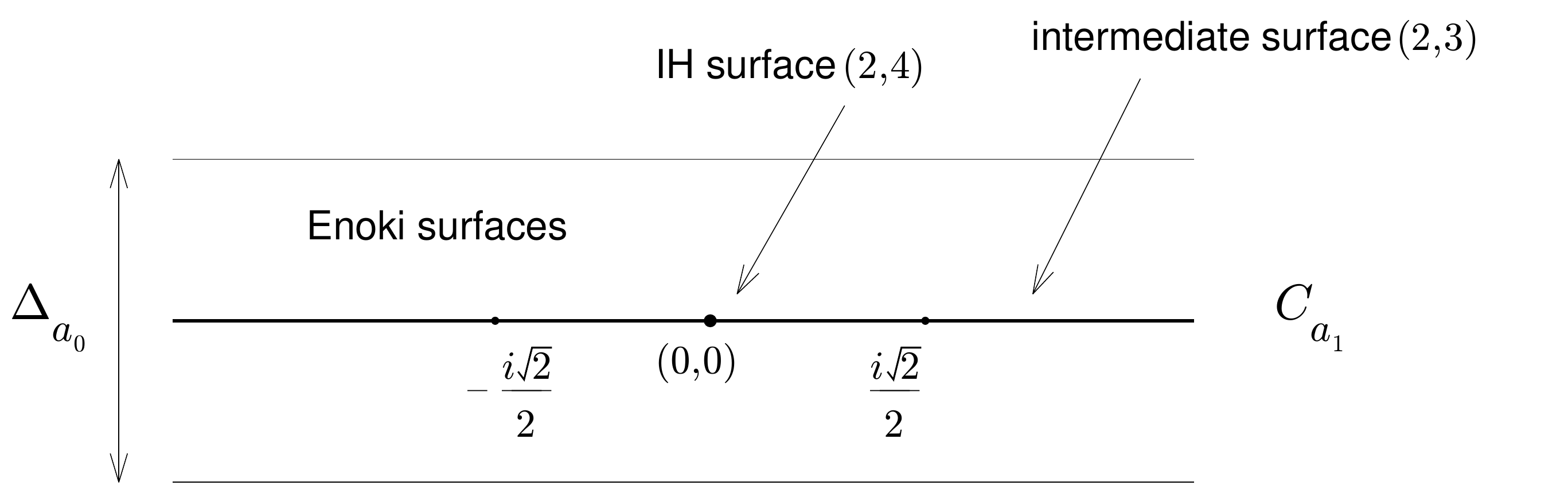}
\end{center}
Notice that for $\f(z_1,z_2)=(-z_1,-z_2)$, 
$$\f\circ G_{(a_0,a_1)}\circ\f^{-1}=G_{(a_0,-a_1)},$$
 therefore there is an involution 
$$i:B_J\to B_J,\quad i(a_0,a_1)=(a_0,-a_1),$$
such that $G_{a}$ and $G_{i(a)}$ give isomorphic surfaces.\\
The mapping from $B_J$ to moduli space is ramified along $\{a_1=0\}$ and $\{a_1=0\}$ is an irreducible component of the  hypersurface $T_{J,\s}$ of the points where the cocycles are not linearly independant. The sheaf of relations is generated by global section by theorem A of Cartan. Let 
$$\a_0(a)[\t_0]+\b_0(a)[\mu_0]+\a_1(a)[\t_1]+\b_1(a)[\mu_1]=0$$
be such a relation. By the same computation as at the beginning of section \ref{section4.2},  
$$\b_0=\b_1=0,$$
therefore we have to solve the system
$$\left\{\begin{array}{lcll}
X_0-\Pi_{1\star}X_1&=&\dps\a_0\frac{\part}{\part u_0}&{\rm at \ the \ point\ }\Pi_1(u_1,v_1)\\
&&&\\
X_1-\s_\star\Pi_{0\star}X_0&=&\dps\a_1\frac{\part}{\part u'_1}&{\rm at \ the \ point\ }\s\Pi_0(u'_0,v'_0)
\end{array}\right.$$
We have
$$D\Pi_1(u_1,v_1)=\left(\begin{array}{cc}v_1&u_1\\0&1\end{array}\right),\quad D(\s\Pi_0)(u'_0,v'_0)=\left(\begin{array}{cc}0&1\\v'_0&u'_0\end{array}\right),$$
Since by Hartogs theorem the vector fields $X_0$ and $X_1$ extend on the whole blown up ball, they are tangent to the exceptional curves and we set
$$X_0=A_0\frac{\part}{\part u'_0}+v'_0B_0\frac{\part}{\part v_0},\quad X_1=A_1\frac{\part}{\part u_1}+v_1B_1\frac{\part}{\part v_1},$$
By straightforward computations similar to those in the appendix we derive that 
$$\a_0(a_0,a_1)=0,$$
for all minimal surfaces, therefore $[\t_0]\neq 0$  on $\D_{a_0}\times \bb C_{a_1}$ (recall that $tr(S_a)=a_0$ and the trace is a holomorphic invariant) .\\
 By proposition 3.\ref{relationlogarithmique},  the four cocycles are independent on the line $\{a_0=0\}$, hence $\a_1(0,a_1)=0$, and the relation reduces to 
$$\a_1(a)[\t_1]=0$$
 with $\a_1(0,a_1)=0$. Therefore $T_{J,\s}\cap \{a_0=0\}=\{(0,\pm \frac{i\sqrt{2}}{2})\}$, $[\t_1]=0$ along $T_{J,\s}\setminus\{a_0=0\}$ but $[\t_1]\neq 0$ at the two points where $\T$ is not locally free. The  mapping from the stratum of Enoki surface to the moduli space of Enoki surfaces is a finite morphism ramified along $T_{J,\s}=\overline{\{a\mid [\t_1](a)=0\}}$, in particular surjective on the moduli space of Enoki surfaces. \\ 

{\bf \boldmath Case $J=\{0,1\}$}\\
With $\s(z)=(z_1+a_1,z_2)$, the family of marked surfaces $\Phi_{J,\s}:\cal S_{J,\s}\to B_J$ is associated to the family of polynomial germs
$$G^J_a(z_1,z_2)=G_a(z_1,z_2)=\bigl(z_2(z_1+a_1), z_2(z_1+a_1)(z_2+a_0)\bigr)$$
$$\det DG_a(z_1,z_2)=z_2^2(z_1+a_1),$$
$$\tr(S_a)=\tr DG_a(0)=a_0a_1,\quad {\rm with}\quad |a_0a_1|<1.$$
There are two lines of intermediate surfaces which meet at the point $(0,0)$ which parametrize the Inoue-Hirzebruch surface with two singular rational curves.
\begin{center}
\includegraphics[width=8cm]{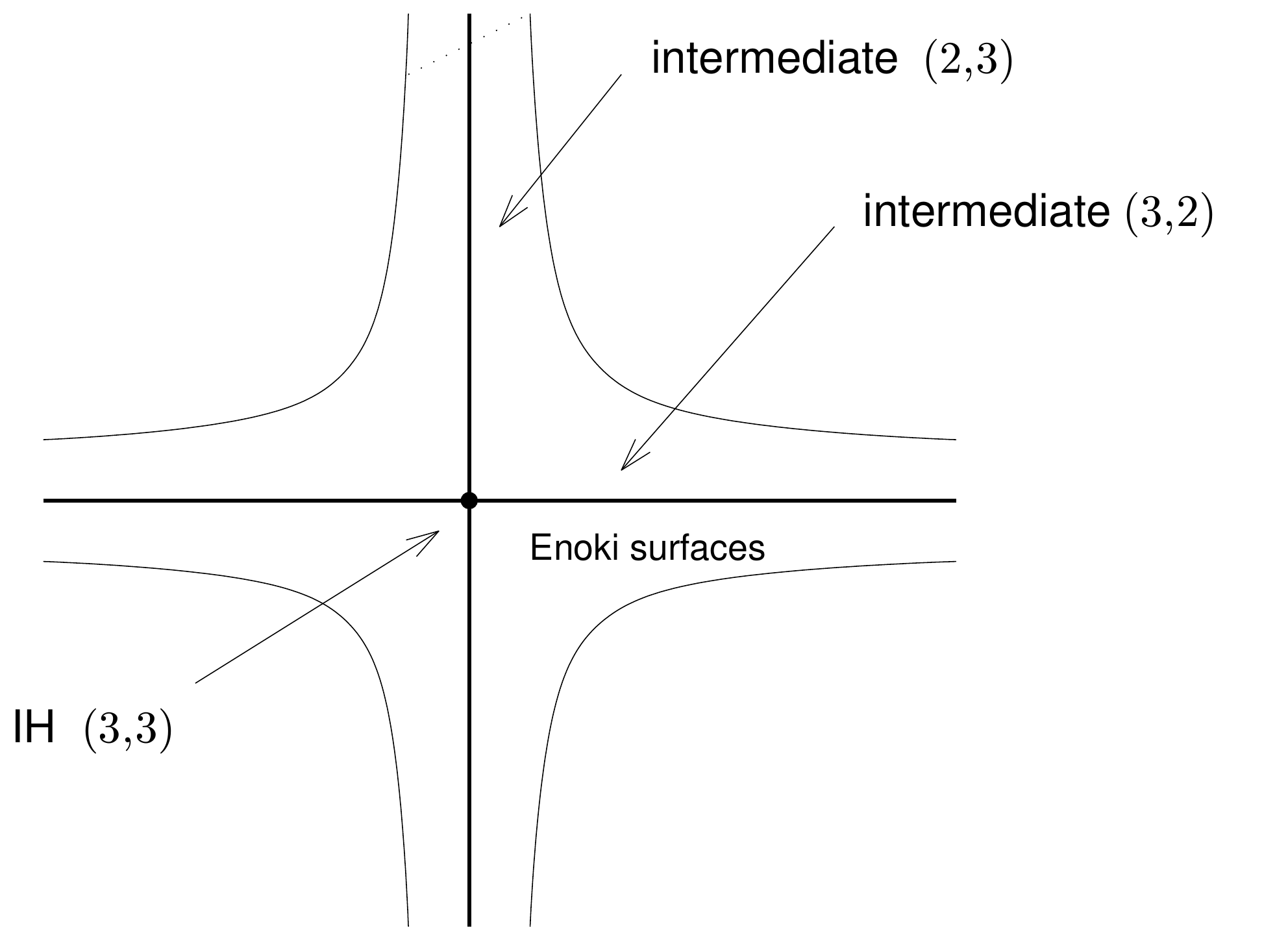}
\end{center}
\begin{itemize}
\item For $a_0=0$, $a_1\neq 0$, $\k=a_1$,
\item For $a_1=0$, $a_0\neq 0$, $\k=a_0$.
\end{itemize}
The involution of the Inoue-Hirzebruch surface which swaps the two cycles induces on the base of the versal family swapping of the two lines of intermediate surfaces.\\
We have obtained
\begin{Th}  Let $F=\Pi\s:(\bb C^2,0)\to (\bb C^2,0)$ be any holomorphic germ, where $\Pi=\Pi_0\Pi_1$ are blowing-ups and $\s$ is any germ of isomorphism. Then $F$ is conjugated to a birational map obtained by the composition of two blowing-ups 
$$(u,v)\mapsto (uv+a,v),\quad (u',v')\mapsto (v',u'v'),$$
and an affine map at a suitable place. If moreover, $S$ is of intermediate type and there is no non-trivial invariant vector field, $F$ is conjugated to the composition of two blowing-ups of the previous types.\\
\end{Th}
\begin{Cor}
Any  minimal surface  with $b_1(S)=1$ and $b_2\le 2$ containing a GSS admits a birational structure.
\end{Cor}

\section{Appendix}
\subsection{On logarithmic deformations of surfaces with GSS, by  Laurent Bruasse} \label{Bruasse}
The results contained in this section is a  not yet published part of the thesis \cite{B}. Notations are those of \cite{Br} and \cite{DO1}.\\

Let $\cal F$ be a reduced foliation on a  compact complex surface $S$. We denote by $T_{\cal F}$ (resp. $N_{\cal F}$) the tangent (resp. normal) line bundle  to $\cal F$.\\
Let $p$ be a singular point of the foliation; in a neighbourhood of $p$ endowed with a coordinate system $(z,w)$ in which $p=(0,0)$, $\cal F$ is defined by a holomorphic vector field
$$\t(z,w)=A(z,w)\frac{\part}{\part z} + B(z,w)\frac{\part}{\part w}.$$
Let $J(z,w)$ be the jacobian matrix of the mapping $(A,B)$. Baum-Bott \cite{BB} and Brunella \cite{Br} have introduced the following two indices:
$$Det(p,\cal F)=Res_{(0,0)}\frac{\det J(z,w)}{A(z,w)B(z,w)}dz\wedge dw$$
$$Tr(p,\cal F)=Res_{(0,0)}\frac{\bigl(\tr J(z,w)\bigr)^2}{A(z,w)B(z,w)}dz\wedge dw$$
where $Res_{(0,0)}$ is the residue at $(0,0)$ (see \cite{GH} p649). We denote by $S(\cal F)$ the singular set of $\cal F$, it is a finite set of points, and let
$$Det\,\cal F:=\sum_{p\in S(\cal F)} Det(p,\cal F), \quad Tr(\cal F):=\sum_{p\in S(\cal F)} Tr(p,\cal F).$$
\begin{Prop}[Baum-Bott formulas, \cite{BB},\cite{Br}] We have
$$Det\,\cal F = c_2(S)-c_1(T_{\cal F}).c_1(S)+c_1(T_{\cal F})^2,$$
$$Tr(\,\cal F) = c_1(S)^2-2c_1(T_{\cal F}).c_1(S)+c_1(T_{\cal F})^2.$$
\end{Prop}
By \cite{DO1}, if $S$ is a minimal compact complex surface with GSS, then
$$Det\,\cal F =n, \quad Tr(\,\cal F) =2n-\s_n(S).$$

\begin{Prop}\label{dimdeffeuilletage} Let $S$ be a minimal surface containing a GSS with $n=b_2(S)\ge 1$ and $tr(S)=0$. If $\cal F$ is a reduced foliation on $S$, then
$$h^1(S,T_{\cal F})=
\left\{\begin{array}{lcl}3n-\s_n(S)&{\rm if}& h^0(S,\T)=0\\3n-\s_n(S)+1&{\rm if}& h^0(S,\T)=1
\end{array}\right.$$
If there is no non-trivial global vector fields this integer is precisely the number of generic blowing-ups.
\end{Prop}
Proof: By Riemann-Roch formula
$$\begin{array}{lcl}
h^0(T_{\cal F})-h^1(T_{\cal F})+h^2(T_{\cal F}) & = & \chi(S)+\frac{1}{2}\Bigl( c_1(T_{\cal F})^2- c_1(T_{\cal F})c_1(K)\Bigr)\\
&&\\
&=&c_1(T_{\cal F})^2 = \s_n(S)-3n
\end{array}$$
since 
\begin{itemize}
\item by the first Baum-Bott formula 
$c_1(T_{\cal F}).c_1(S)=c_1(T_{\cal F})^2$ and
\item  by the second and the previous observation
$c_1(T_{\cal F})^2=-Tr(\cal F)+c_1(S)^2=\s_n(S)-3n$.
\end{itemize}
Suppose first that $S$ is of intermediate type. Two cases occur
$$h^0(T_{\cal F})=\left\{\begin{array}{cl}0&{\rm if}\ h^0(S,\T)=0\\ 1&{\rm if\ not}
\end{array}\right.$$
Moreover, by Serre duality $h^2(T_{\cal F})=h^0(K\ot T_{\cal F}^\star)$ and
$$c_1(K).\bigl(c_1(K)-c_1(T_{\cal F})\bigr)=c_1(S)^2+c_1(S).c_1(T_{\cal F})=-n+(\s_n(S)-3n)=-4n+\s_n(S)<0$$
 Let $e_i$, $i=0,\ldots,n-1$ be the Donaldson classes in $H^2(S,\bb Z)$ which trivialize the negative intersection form. In $H^2(S,\bb Z)$, $c_1(K)=\sum_{i=0}^{n-1} e_i$ and $c_1(K)-c_1(T_{\cal F})=\sum_{i=0}^{n-1}a_ie_i$. Since 
$c_1(K).\bigl(c_1(K)-c_1(T_{\cal F})\bigr)<0$ we have $\sum_i a_i>0$ therefore $h^0(K\ot T_{\cal F}^\star)=0$ by \cite{N2} Lemma (2.3). \\
If $S$ is a Inoue-Hirzebruch surface, there are two foliations, each defined by a twisted vector field $\t\in H^0(S,\T\ot L^{\l})$, with $\l$ an irrational quadratic number (see \cite{DO1}), hence $T_{\cal F}=L^{1/\l}$. We have $-K=D$ or $-K^{\ot 2}=2D$ and there is no topologically trivial divisor, therefore 
$$h^0(L^{1/\l})=0,\quad{\rm and}\quad h^2(L^{1/\l})=h^0(K\ot L^\l)=0.$$
We conclude by Riemann-Roch theorem that
$$h^1(T_{\cal F})=h^1(L^{1/\l})=0=3n-\s_n(S)$$
which is the annouced result.\hfill$\Box$\\

We have a canonical injection $0\stackrel{i}{\to} T_{\cal F}\to \T(-Log\, D)$. The aim of the following proposition is to compare logarithmic deformations and deformations which respect the foliation:
\begin{Prop} There exists  an exact sequence of sheaves of $\cal O_S$-Modules 
$$0\to T_{\cal F}\stackrel{i}{\to} \T(-Log\, D)\to N_{\cal F}\otimes \cal O(-D)\to 0.\leqno{(\spadesuit)}$$
\end{Prop}
Proof: Let $\cal U=(U_i)$ be a finite  covering by open sets endowed with holomorphic 1-forms $\omega_i$ defining the foliation $\cal F$. On each open set $U_i$ we consider the morphism
$$\begin{array}{cccc}j:&\T(-Log\,D)_{\mid U_i}&\longrightarrow& N_{\cal F}\ot \cal O_S(-D)_{\mid U_i}\\
&\t&\longmapsto&\omega_i(\t)
\end{array}$$
Since $\t$ is tangent to $D$, $\omega_i(\t)$ vanishes on $D$, therefore the morphism is well defined on $U_i$. Moreover, by definition, the normal bundle $N_{\cal F}$ is defined by the cocycle $(g_{ij})_{ij}=(\omega_i/\omega_j)_{ij}\in H^1(\cal U,\cal O^\star)$, therefore $j$ is well defined on $S$ and its kernel is clearly $\im i$. It remains to check that $j$ is surjective: outside $D$ it is obvious since the foliation has singular points only at the intersection of two curves and we have the exact sequence
$$0\to T_{\cal F}\stackrel{i}{\to} \T(-Log\, D)_{\mid S\setminus D}=\T_{\mid S\setminus D}\to N_{\cal F}\to 0.$$
Let $x\in D$, $f_x\in N_{\cal F,x}\otimes \cal O(-D)_x$ and $U$ an open neighbourhood of $x$ on which $f$ is defined.
\begin{itemize}
\item If $x$ is not at the intersection of two curves, let $(z,w)$ be a coordinate system in which $D=\{z=0\}$ and $\cal F$ defined by $\omega=a(z,w)dz+zb(z,w)dw$. Since $f$ vanishes on $D$, $f=zg$. Let $\t=z\a(z,w)\frac{\part}{\part z}+\b(z,w)\frac{\part}{\part w}$ be a logarithmic vector field. We have to find $\a$ and $\b$ such that
$$f(z,w)=zg(z,w)=\omega(\t)=z(a\a+b\b)$$
i.e. $g\in (a,b)$. The are solutions because $x$ is not a singular point of the foliation hence $a$ is invertible at $x$.
\item If $x$ is at the intersection of two curves, 
$$\omega=wa\,dz + zb\,dw,  \quad f=zwg\quad{\rm and}\quad \t=z\a(z,w)\frac{\part}{\part z}+w\b(z,w)\frac{\part}{\part w}.$$
 We have to solve $g=a\a+b\b$. By \cite{K} p171 (see also \cite{DO1} p1528), the order of $\t$ is one, hence $a$ or $b$ is invertible and $g\in (a,b)$.
\end{itemize}
\hfill$\Box$\\

Let $S$ be a minimal compact complex surface containing a GSS with $n=b_2(S)\ge 1$ and $tr(S)=0$. If $S$ is not a Inoue-Hirzebruch surface then $S$ admits a unique holomorphic foliation $\cal F$ \cite{DO1} p1540 given by a $d$-closed section of $H^0(S,\Omega^1(Log\, D)\ot L^{k(S)})$. If $S$ is a Inoue-Hirzebruch surface, it admits exactly two foliations defined by twisted vector fields.\\
The exact sequence ($\spadesuit$) yields
$$0\to H^1(S,T_{\cal F})\to H^1(S,\T(Log\,D))\to H^1(S,N_{\cal F}\ot \cal O(-D)) \to H^2(S,T_{\cal F})$$
In fact, if $S$ is not a Inoue-Hirzebruch surface, $\cal F$ is unique, thence $\Omega^1$ contains a unique non-trivial coherent subsheaf which is  $\cal O(-D)\ot L^{1/k}$. As $N_{\cal F}^\star$ is another, $N_{\cal F}=\cal O(D)\ot L^k$ and $H^0(S, N_{\cal F}\ot \cal O(-D))=H^0(S,L^{k})=0$ because $k\neq 1$. We have also $h^2(S,N_{\cal F}\ot \cal O(-D))=h^2(S,L^k)=h^0(S,K\ot L^{1/k})=0$. By Riemann-Roch theorem, $h^1(S,N_{\cal F}\ot \cal O(-D))=0$ and we obtain the isomorphism
$$0\to H^1(S,T_{\cal F})\to H^1(S,\T(Log\,D))\to 0$$
If $S$ is a Inoue-Hirzebruch surface $N_{\cal F}=\cal O(D)\ot L^\l$ where $\l$ is an irrationnal  number and we have the same conclusion.\\
With (\ref{dimdeffeuilletage}) we have proved:

\begin{Th} Let $S$ be a minimal compact complex surface containing a GSS with $n=b_2(S)\ge 1$ and $tr(S)=0$. Then:
$$h^1(S,\T(Log\, D))=h^1(S,T_{\cal F})=\left\{ \begin{array}{lcl}
3n-\s_n(S)&{\rm if}&h^0(S,\T)=0\\
3n-\s_n(S)+1&{\rm if}&h^0(S,\T)=1
\end{array}\right.$$
In particular any logarithmic deformation keeps the foliation.
\end{Th}
\begin{Rem} If $tr(S)\neq 0$, the theorem remains true by \cite{DK}.
\end{Rem}

\subsection{Infinitesimal logarithmic deformations: the hard part}\label{thehardpart}
We give in this section the proof of proposition 3.\ref{relationlogarithmique}.\\
Since $\s(0)=O_{n-1}$ is the intersection of two transversal rational curves which are contracted by $F$, there is  a conjugation by a linear map $\f$ (in particular birational) such that $\f^{-1}F\f=F'=\Pi'\s'$, satisfies
$$\frac{\part\s'_1}{\part z_2}(0)=\frac{\part\s'_2}{\part z_1}(0)=0.\leqno{(S)}$$
 It means that $\s'^{-1}(C_{n-1})$ is tangent to $z_2=0$ and the other curve is tangent to $z_1=0$, therefore their strict transforms meet the exceptional curve $C_0$ respectively at $\{u'=v'=0\}$ and $\{u=v=0\}$.\\
 Therefore in the following computations we shall suppose that the condition $(S)$ is satisfied.
Let $\Pi''=\Pi_l\circ\Pi_{l+1}\circ\cdots\circ\Pi_{n-1}$ be the composition of blowing-ups at the intersection of two curves and of $\Pi_l$, then it is the composition of mappings $(u,v)\mapsto (uv,v)$ or $(u',v')\mapsto (v',u'v')$, and of $\Pi_l(u',v')=(v'+a_{l-1},u'v')$, hence 
$$\Pi''(x,y)=(x^py^q+a_{l-1},x^ry^s)$$
 where $\left(\begin{array}{cc}p&q\\r&s\end{array}\right)$ is the composition of matrices $\left(\begin{array}{cc}1&1\\0&1\end{array}\right)$ or $\left(\begin{array}{cc}0&1\\1&1\end{array}\right)$, the last one being of the second type, therefore  
 $$\det \left(\begin{array}{cc}p&q\\r&s\end{array}\right)=\pm 1.$$
We have
$$\Pi''\s \Pi_0(u_0,v_0)=\Pi''\bigl(\s_1 (u_0v_0,v_0),\s_2(u_0v_0,v_0)\bigr)= \bigl(\s_1^p\s_2^q(u_0v_0,v_0)+a_{l-1},\s_1^r\s_2^s(u_0v_0,v_0)\bigr).$$

{\bf First case:} there are at least two singular sequences, then
$$1\le p\leq r,\quad 1\le q\le s,\quad p+q<r+s.$$

  The jacobian is
$$D(\Pi''\s \Pi_0)(u_0,v_0) = \left(
\begin{array}{cc}
v_0\s_1^{p-1}\s_2^{q-1}(u_0v_0,v_0)P(u_0,v_0) &  \s_1^{p-1}\s_2^{q-1}(u_0v_0,v_0)Q(u_0,v_0) \\
 &\\
 v_0\s_1^{r-1}\s_2^{s-1}(u_0v_0,v_0)R(u_0,v_0) & \s_1^{r-1}\s_2^{s-1} (u_0v_0,v_0)S(u_0,v_0)\end{array}\right)$$
where 
$$\left\{\begin{array}{l}
P(u,v)= p\s_2(uv,v)\partial_1\s_1(uv,v)+q\s_1(uv,v)\partial_1\s_2(uv,v),\\
\\
Q(u,v)= p\s_2(uv,v)\Bigl(u\partial_1\s_1(uv,v)+\partial_2\s_1(uv,v)\Bigr) + q\s_1(uv,v)\Bigl(u\partial_1\s_2(uv,v)+\partial_2\s_2(uv,v)\Bigr)\\
\\
R(u,v)=r\s_2(uv,v)\partial_1\s_1(uv,v)+s\s_1(uv,v)\partial_1\s_2(uv,v)\\
\\
S(u,v)= r \s_2(uv,v)\Bigl( u\partial_1\s_1(uv,v) + \partial_2\s_1(uv,v)\Bigr) + s\s_1(uv,v) \Bigl( u\partial_1\s_2(uv,v)+\partial_2\s_2(uv,v)\Bigr)
\end{array}\right.$$

For $i=1,\ldots,l-1$ we have also
$$D\Pi_i(u_i,v_i)=\left(\begin{array}{cc}v_i&u_i\\ 0&1\end{array}\right)$$

In the local chart $(u_i,v_i)$ containing $O_i$, for $i=0,\ldots,l-1$, $X_i$ is tangent to $C_i=\{v_i=0\}$, hence we have
$$X_i(u_i,v_i)=\left(\begin{array}{c}A_i(u_i,v_i)\\ \\v_iB_i(u_i,v_i)\end{array}\right).$$

For $i=0,\ldots,l-2$, we have at the point $$(u_i,v_i)=\Pi_{i+1}(u_{i+1},v_{i+1})=(u_{i+1}v_{i+1}+a_{i},v_{i+1}),$$
$$\left(\begin{array}{c}A_i(u_{i+1}v_{i+1}+a_{i},v_{i+1})\\ \\v_{i+1}B_i(u_{i+1}v_{i+1}+a_{i},v_{i+1})\end{array}\right)- \left(\begin{array}{cc} v_{i+1}&u_{i+1}\\ &\\ 0&1\end{array}\right) \left(\begin{array}{c}A_{i+1}(u_{i+1},v_{i+1})\\  \\v_{i+1}B_{i+1}(u_{i+1},v_{i+1})\end{array}\right)= 
\left(\begin{array}{c}\a_i\\ \\0\end{array}\right)$$

For $i=l-1$, at the point 
$$\begin{array}{lcl}
(u_{l-1},v_{l-1})&=& \Pi''\circ\s\circ\Pi_0(u_0,v_0)= \Pi''\bigl(\s_1(u_0v_0,v_0), \s_2(u_0v_0,v_0))\bigr)\\
&&\\
&=&
\bigl(\s_1^p\s_2^q(u_0v_0,v_0)+a_{l-1},\s_1^r\s_2^s(u_0v_0,v_0)\bigr),
\end{array}$$

$$
\left(\begin{array}{c}A_{l-1}(\Pi''\s\Pi_0(u_0,v_0))\\ \\
\s_1^r\s_2^s(u_0v_0,v_0)B_{l-1}(\Pi''\s\Pi_0(u_0,v_0))\end{array}\right)- D(\Pi''\s\Pi_0)(u_0,v_0)\left(\begin{array}{c}A_0(u_0,v_0)\\  \\v_0B_0(u_0,v_0)\end{array}\right) 
= 
\left(\begin{array}{c} \a_{l-1}\\ \\0\end{array}\right)$$
Equivalently, we obtain

For $i=0,\ldots, l-2$,
$$A_i(u_{i+1}v_{i+1}+a_i,v_{i+1}) - v_{i+1}\left\{ A_{i+1}(u_{i+1},v_{i+1})+u_{i+1}B_{i+1}(u_{i+1},v_{i+1})\right\} = \a_i  \leqno{(I_i)}$$   
$$B_i(u_{i+1}v_{i+1}+a_i,v_{i+1}) -   B_{i+1}(u_{i+1},v_{i+1}) =0
\leqno{(II_i)}$$

For $i=l-1$, omitting subscripts,
$$A_{l-1}( \Pi''\circ\s\circ\Pi_0(u,v))
 - v\s_1^{p-1}\s_2^{q-1}(uv,v)\Bigl\{ P(u,v)A_0(u,v) +Q(u,v)B_0(u,v)\Bigr\}  = \a_{l-1} \leqno{(I_{l-1})}$$

$$\s_1\s_2(uv,v)B_{l-1}( \Pi''\circ\s\circ\Pi_0(u,v)) 
- v\Bigl\{R(u,v)A_0(u,v) +S(u,v)B_0(u,v)\Bigr\}
=0 \leqno{(II_{l-1})}$$

 For $i=0,\ldots,l-1$ and for $v_{i+1}=0$  the equations $(I_i)$ yield,
$$A_i(a_i,0)=\a_i.\leqno{(1)}\label{1}$$
We put $u_i=t_i+a_i$, $i=0,\ldots,l-1$,
$$A_i(u_i,v_i)=A_i(a_i,0)+A'_i(t_i,v_i)=A_i(a_i,0)+\sum_{j+k>0}a^i_{j,k}t_i^jv_i^k,$$ 
$$B_i(u_i,v_i)=B_i(a_i,0)+B'_i(t_i,v_i)=B_i(a_i,0)+ \sum_{j+k>0}b^i_{j,k}t_i^jv_i^k.$$
For $i=0,\ldots,l-2$, equations $(II_{i})$ give
$$B:=B_0(a_0,0)=\cdots =B_{l-1}(a_{l-1},0),\leqno{(2)}$$
$$B'_1(t_1,0)=\cdots=B'_{l-1}(t_{l-1},0)=0.\leqno{(3)}$$

Replacing $A_i$ and $B_i$ by their expressions we have by (2),

For $i=0,\ldots, l-2$,
$$\begin{array}{l}
A'_i\bigl((t_{i+1}+a_{i+1})v_{i+1},v_{i+1}\bigr) - v_{i+1}\Bigl\{ A_{i+1}(a_{i+1},0)+ A'_{i+1}(t_{i+1},v_{i+1}) \\
\hspace{40mm} + (t_{i+1}+a_{i+1})\bigl[B + B'_{i+1}(t_{i+1},v_{i+1})\bigr]\Bigr\} =0
\end{array}\leqno{(I'_i)}$$

$$B'_i\bigl((t_{i+1}+a_{i+1})v_{i+1},v_{i+1}\bigr) -   B'_{i+1}(t_{i+1},v_{i+1}) =0
\leqno{(II'_i)}$$

$$\left\{\begin{array}{l}
A'_{l-1}(\s_1^p\s_2^q(uv,v),\s_1^r\s_2^s(uv,v)) \\ \\
 -v\s_1^{p-1}\s_2^{q-1}(uv,v)\Bigl\{ P(u,v)\bigl[A_0(a_0,0)+A'_0(t,v)\bigr]
 +  Q(u,v)\bigl[B+B'_0(t,v)\bigr]\Bigr\} = 0\end{array}\right.\leqno{(I'_{l-1})}$$

$$\left\{ \begin{array}{l}
\s_1\s_2(uv,v)\bigl[B+B'_{l-1}(\s_1^p\s_2^q(uv,v),\s_1^r\s_2^s(uv,v))\bigr] \\ \\
\hspace{23mm}- v\Bigl\{R(u,v)\bigl[A_0(a_0,0)+A'_0(t,v)\bigr] 
  + S(u,v)\bigl[B+B'_0(t,v)\bigr]\Bigr\}= 0
\end{array}\right. \leqno{(II'_{l-1})}$$

Now, from the equations $I'_i$ and $II'_i$, $0\le i\le l-1$, we show that some terms vanish. In fact:
For $i=0,\ldots, l-2$, we divide $(I'_i)$  by $v_{i+1}$, we set $v_{i+1}=0$,  we apply  $(3)$, and we compare linear terms:
$$a^0_{10}-a^1_{10}=\cdots = a^{l-2}_{10}-a^{l-1}_{10}=B. \leqno{(4)}$$
Equations (4) give 
$$a^0_{10}-a^{l-1}_{10}-(l-1)B=0\leqno{(5)}$$

For $i=0,\ldots,l-2$, we divide $(II'_i)$ by $v_{i+1}$, we set $v_{i+1}=0$,  we apply  $(3)$:
$$b^0_{10}a_1+b^0_{01}-b^1_{01}=0,\quad  b^1_{01}=\cdots=b^{l-1}_{01}.\leqno{(6)}$$
$$b^1_{10}=\cdots =b^{l-1}_{10}=0.\leqno{(7)}$$

Dividing $(I'_{l-1})$ by $v^2\s_1^{p-1}\s_2^{q-1}(uv,v)$, setting $v=0$, recalling that $\partial_1\s_2(0)=\part_2\s_1(0)=0$ and cancelling the factor $\partial_1\s_1(0)\partial_2\s_2(0)\neq 0$, we obtain 
$$a_{1,0}^{l-1}(t+a_0) - \Bigl\{ p\bigl[A_0(a_0,0)+A'_0(t,0)\bigr]
+ (p+q)(t+a_0)\bigl[B+B'_0(t,0)\bigr] \Bigr\} =0
\leqno{(8)}$$
Constant part of (8) is
$$a_0 a_{10}^{l-1} - \left\{ pA_0(a_0,0)+(p+q)a_0B\right\}=0 \leqno{(9_c)}$$
Linear part of (8) is
$$ pa^0_{10}-a^{l-1}_{10}  + (p+q)B + (p+q)a_0 b^0_{10}=0 \leqno{(9_l)}$$
Dividing $(II'_{l-1})$ by $v^2$, setting $v=0$ and cancelling the factor term $\partial_1\s_1(0)\partial_2\s_2(0)\neq 0$, we obtain
$$(t+a_0)B - r\bigl[A_0(a_0,0)+A'_0(t,0)\bigr] 
 - (r+s)(t+a_0)\bigl[B+B'_0(t,0)\bigr] =0
\leqno{(10)}$$
Constant part of (10) is
$$ rA_0(a_0,0)+ (r+s-1)a_0B =0 \leqno{(11_c)}$$
Linear part of (10) is
$$ra^0_{10} + (r+s-1)B +(r+s)a_0b^0_{10} =0 \leqno{(11_l)}$$
The determinant of the linear system (5), (9$_c$), (9$_l$), (11$_c$) and (11$_l$) with unknowns $a^0_{10}$, $a^{l-1}_{10}$, $A_0(a_0,0)$, $B$ and $b^0_{10}$ is 
$$\begin{array}{lcl}\D&=&\left| \begin{array}{ccccc} 1&-1&0&-(l-1)&0\\&&&&\\
0&a_0&-p&-(p+q)a_0&0\\&&&&\\
p&-1&0&(p+q)&(p+q)a_0\\&&&&\\
0&0&r&(r+s-1)a_0&0\\&&&&\\
r&0&0&(r+s-1)&(r+s)a_0\end{array}\right|\\
&&\\
&&\\
& = & a_0^2(ps-qr)\Bigl\{(ps-qr)+1-(p+s)-rl\Bigr\}\neq 0
\end{array}$$
 
 Therefore, by (4), (6) and (7)
$$a^0_{10}=\cdots =a^{l-1}_{10}=B=0, \quad b^0_{01}=\cdots=b^{l-1}_{01}\quad {\rm and} \quad 
b^0_{10}=\cdots =b^{l-1}_{10}=0.\leqno{(12)}$$
Moreover we obtain 
$$\begin{array}{|c|}\hline \overset{\ }{\underset{\  }{\a_0=A_0(a_0,0)= 0}}. \\\hline\end{array}\leqno{(13)}$$

{\bf Second case:} there is only one singular seqence $s_m$, $m\ge 1$, then 
$$\left(\begin{array}{cc}p&q\\r&s\end{array}\right)=\left(\begin{array}{cc}0&1\\1&1\end{array}\right)\left(\begin{array}{cc}1&1\\0&1\end{array}\right)^{m-1}=\left(\begin{array}{cc}0&1\\1&m\end{array}\right)$$
$$\Pi''\s\Pi_0(u_0,v_0)=\bigl(\s_2(u_0v_0,v_0)+a_{l-1},\s_1\s_2^m(u_0v_0,v_0)\bigr)$$
$$D\bigl(\Pi''\s\Pi_0\Bigr)(u,v)=\left(
\begin{array}{cc}vP(u,v)&Q(u,v)\\&\\v\s_2^{m-1}(uv,v)R(u,v)&\s_2^{m-1}(uv,v)S(u,v)\end{array}
\right)$$
where
$$\left\{\begin{array}{lcl}P(u,v)&=&\part_1\s_2(uv,v)\\
&&\\
Q(u,v)&=&u\part_1\s_2(uv,v)+\part_2\s_2(uv,v)\\
&&\\
R(u,v)&=&\s_2\part_1\s_1(uv,v)+m\s_1\part_1\s_2(uv,v)\\
&&\\
S(u,v)&=&\s_2(uv,v)\Bigl(u\part_1\s_1(uv,v)+\part_2\s_1(uv,v)\Bigr)\\
&&+m\s_1(uv,v)\Bigl(u\part_1\s_2(uv,v)+\part_2\s_2(uv,v)\Bigl)
\end{array}\right.$$
The new equations are
$$\left\{\begin{array}{ll}(I'_{l-1})&A'_{l-1}(\s_2(uv,v),\s_1\s_2^m(uv,v))\\
&\hspace{25mm}-v\Bigl\{P(u,v)[A_0(a_0,0)+A'_0(t,v)]+Q(u,v)[B+B'_0(t,v)]\Bigr\}=0\\
&\\
(II'_{l-1})&\s_1\s_2(uv,v)\Bigl(B+B'_{l-1}(\s_2(uv,v),\s_1\s_2^m(uv,v))\Bigr)\\
&\hspace{25mm}-v\Bigl\{R(u,v)[A_0(a_0,0)+A'_0(t,v)]+S(u,v)[B+B'_0(t,v)]\Bigr\}=0
\end{array}\right.$$

The end of the proof  follows the same lines than in the first case. Details are left to the reader.\\

\begin{Rem} By induction it is possible to show that for any $k<r+s-(p+q)$, a similar Cramer system may by defined and that $\a_k=0$. However, it is not possible to achieve the proof in this way because when $k=r+s-(p+q)$ a new unknown appears. This difficulty is explained by the fact that  {\bf there is a relation} among the $\t^i$'s over an hypersurface in $B_J$.
\end{Rem}

Georges Dloussky,\\
 Aix-Marseille University (AMU), department of Mathematics, 39 rue F.Joliot-Curie, Marseille 13013,\\
 Laboratoire d'Analyse, Topologie et Probabilités (LATP), Unité Mixte de Recherche (UMR) 7353 CNRS\\
  georges.dloussky@univ-amu.fr

\end{document}